%% file: main.tex
\documentclass[11pt]{article}{}
\pdfoutput=1

\newcommand{\doctitle}{Numerical analysis of non-local calculus on finite weighted graphs, with application to reduced-order modelling of dynamical systems}
\newcommand{\docauthor}{M. Duschenes, S. Srivastava, \& K. Garikipati}
\newcommand{\docaffil}{Applied Physics, Mathematics, Mechanical Engineering, and Michigan Institute for Computational Discovery \& Engineering\\ University of Michigan}

\newcommand{\docheader}{Numerical analysis ... dynamical systems - \docauthor}
\newcommand{\docfooter}{University of Michigan}

\newcommand{\pdftitle}{Numerical analysis ... dynamical systems}

\usepackage{main}

\usepackage{xpatch}
\xapptocmd\appendices{%
  \crefalias{section}{appendix}%
}{}{\PatchFailed}
\crefname{section}{Sec.}{Secs.}
\Crefname{section}{Sec.}{Secs.}
\crefname{equation}{Eq.}{Eqs.}
\Crefname{equation}{Eq.}{Eqs.}
\creflabelformat{equation}{#2(#1)#3}
\crefname{figure}{Fig.}{Figs.}
\Crefname{figure}{Fig.}{Figs.}
\creflabelformat{figure}{#2#1#3}
\crefname{table}{Tab.}{Tabs.}
\Crefname{table}{Tab.}{Tabs.}
\creflabelformat{table}{#2#1#3}
\Crefname{appendix}{App.}{Apps.}
\crefname{appendix}{App.}{Apps.}
\Crefname{si}{Supp Inf.}{Supp Inf.}
\crefname{si}{Supp Inf.}{Supp Inf.}
\newcommand{\appsecname}{Supplementary Information}

\begin{document}

\pretitle{\begin{center}\vskip -80pt}%
\title{\Large\doctitle}
\posttitle{\end{center}}
\preauthor{\begin{center} \vskip -10pt}
\author[1,2,4]{M. Duschenes}
\author[1,4]{S. Srivastava}
\author[1,2,3,4]{K. Garikipati \thanks{Corresponding author at: Department of Mechanical Engineering, University of Michigan, United States. \emph{E-mail address}: krishna@umich.edu (K. Garikipati).}}
\affil[1]{Department of Mechanical Engineering, University of Michigan, United States}
\affil[2]{Applied Physics, University of Michigan, United States}
\affil[3]{Department of Mathematics, University of Michigan, United States}
\affil[4]{Michigan Institute for Computational Discovery \& Engineering, University of Michigan, United States}
\postauthor{\end{center} \vskip -20pt}
\predate{\begin{center} \vskip -0pt}
\date{} 
\postdate{\end{center} \vskip -40pt}%
\maketitle
\begin{abstract}
\pagestyle{abstract}
\noindent We present an approach to reduced-order modelling that builds off recent graph-theoretic work for representation, exploration, and analysis of computed states of physical systems (Banerjee \etal, \emph{Comp. Meth. App. Mech. Eng.}, \textbf{351}, 501-530, 2019). We extend a non-local calculus on finite weighted graphs to build such models by exploiting polynomial expansions and Taylor series. In the general framework for non-local calculus on graphs, the graph edge weights are intricately linked to the embedding of the graph, and consequently to the definition of the derivatives. In a previous communication (Duschenes and Garikipati, arXiv:2105.01740), we have shown that radially symmetric, continuous edge weights derived from, for example Gaussian functions, yield inconsistent results in the resulting non-local derivatives when compared against the corresponding local, differential derivative definitions. Taking inspiration from finite difference methods, we algorithmically compute edge weights, considering the embedding of the local neighborhood of each graph vertex. Given this procedure, we ensure the consistency of the non-local derivatives in this setting, a crucial requirement for numerical applications. We show that we can achieve any desired orders of accuracy of derivatives, in a chosen number of dimensions without symmetry assumptions in the underlying data. Finally, we present two example applications of extracting reduced-order models using this non-local calculus, in the form of ordinary differential equations from parabolic partial differential equations of progressively greater complexity. 
\end{abstract}


\newpage
\clearpage
\singlespacing
\pagenumbering{arabic}
\pagestyle{document}


\newpage

\section{Introduction} \label{sec:intro}
Of interest in this work is the derivation of reduced-order models that are derived from high dimensional solutions to physical systems. Our approach falls within the wide scope of dynamical systems including projection, coarse-graining, dimensionality reduction, and sub-manifold extraction among other methods. It is based on a graph theoretic approach for representation, exploration and analysis of computed states of physical systems that we have proposed recently. \cite{Banerjee2019} This approach draws from the recognition that the complexity of high-dimensional computations is typically distilled into quantities of interest, which we refer to as \emph{states}, for the purpose of design, decision-making and high-throughput computing for optimization. Reduced-order models derived on these states, but which themselves exist on a low-dimensional manifold, have efficiency of representation. Since the original states are computed with high-fidelity, the reduced-order models ``inherit'' some of this accuracy. Their derivation on a principled basis with controlled approximations also renders them susceptible to analysis.

The states that we consider are functionals extracted from the high-dimensional solutions. Examples abound in engineering science and include the lift, drag or pressure difference in computational fluid dynamics, measures of load, deformation and total energies in solid mechanics, phase volume fractions and free energies in materials physics, and cross sections in nuclear physics. \cite{Kochunas2020} These states undergo \emph{transitions} as some \emph{parameter} of the system is varied, such as its physical, possibly time-like, parameters, boundary or initial conditions. 

\begin{figure}[ht]
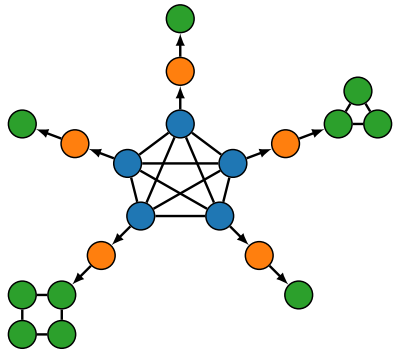

	\centering
	\tikzfig{0.7}{figures/clique} 
	\caption{Example graph, with a central clique, and branching clique subgraphs. Vertex coloring may represent various local attributes.}
	\label{fig:cliquegraph}
\end{figure}

\begin{figure}[ht]
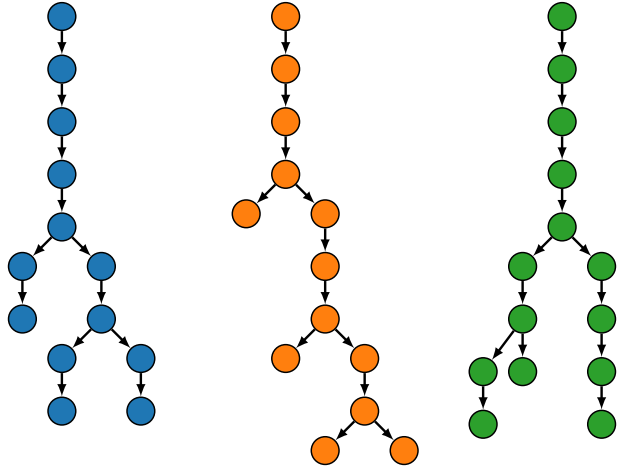

	\centering
	\tikzfig{0.7}{figures/tree} 
	\caption{Example directed tree graphs with branches.}
	\label{fig:treegraph}
\end{figure}

The states, of which there is typically a small number $\sim \mathcal{O}(10)$ relevant to a system, form a low-dimensional vector of quantities of interest. We have shown\cite{Banerjee2019} that by defining the states to be vertices on a graph, and the transitions to be edges between the vertices, a one-to-one correspondence is uncovered between the properties of the computed physical system and the elements of graph theory. Weighted edges between vertices can then be assigned, or found through graph theoretic principles, encoding relationships or transitions between states, and the magnitude of such correlations. Reversible linear systems lead to undirected, fully connected, clique graphs (\cref{fig:cliquegraph}), whereas dissipative dynamical systems are represented by directed trees. Path-dependence, arising in dissipative dynamical systems, is manifested in the acyclicity of trees (\cref{fig:treegraph}). Notions of centrality reveal insights to the relations between states and path traversal properties are induced on graphs by equilibrium or dynamical transitions on the physical system.

Here, we extend the graph theoretic framework to reduced-order modelling based on computed states and transitions of physical systems. The graph of states (vertices) and transitions (edges) is a discrete manifold describing the physical system. On this manifold we seek to define differential equations and functional representations for components of the state vector. A central example in this communication is a reduced-order model for the evolution of a material system in the form of an ordinary differential equation for the phase fraction driven by itself and other components of the state vector such as the free energy and strains. This requires computation of derivatives on these states: of the free energy, with respect to the phase fractions, yielding global chemical potentials, or with respect to the strains, yielding global stress-like quantities. An alternative form is a direct functional representation for the free energy with phase fractions and strains as arguments. In this case, we seek functional representations as approximations to the Taylor series, thus also needing the above derivative computation. 

We have adopted a discrete non-local calculus on weighted graphs for these derivatives on states. This formalism, previously applied to image analysis \cite{Gilboa2008,Desquesnes2013,Hein2007} begins with definition of a vector space on the graph, followed by the gradient operator, inner products and partial derivatives. This non-local calculus is internally consistent, with adjointness relations and a divergence theorem resulting. Its relation to classical differential calculus is determined by the graph weights. Through the inner products, the weights form kernels determining the support of integral operators on the graph manifold, and crucially controlling the behavior of partial derivatives relative to their classical definitions. This is of central importance to numerics on the reduced-order models.

With functions and non-local partial derivatives in hand, algebraic and non-local differential operators can be introduced and viewed as bases in the reduced-order representation. Their evaluation at the vertex data points opens the door to a stepwise regression approach for choosing the best reduced-order model. Parsimony of representation can be dialed in by using thresholds on the stepwise regression.

In this communication we consider first-order dynamical processes in materials that are models of microstructure formation. The generating partial differential equations are gradient flows including the Allen-Cahn and Cahn-Hilliard phase field models, and elasticity. Further examples of physical systems will be explored in future communications.


The non-local calculus on finite weighted graphs is laid out in \cref{sec:graphtheory}, illustrative physical systems and discussions are presented in \cref{sec:physicalsystems}, and conclusions appear in \cref{sec:concl}. Detailed proofs and implementation details appear in the Appendices. 


\section{The graph theoretic approach} \label{sec:graphtheory}
The major formalism used in this work relies on methods from graph theory. \cite{Banerjee2019,West2001,Newman2010} A graph, denoted by $G = (V,E)$, consists of a set of vertices, $x_{} \in V$, with $\abs{V} = n$. The vertices are connected by a set of edges, $e \in E$, where each edge $e = (x_{},y_{})$ is a pair of vertices. A vertex is denoted as being the part of an edge by $x_{} \in e$. The graph can also have global and local attributes on the vertices and edges, as we describe below. Of importance are edge weights $w$, that are generally functions of the local vertex attributes.

We make precise the notions outlined in the Introduction. A state is a low-dimensional vector of dimension $p+s$. It includes \emph{parameters} such as physical constants, or initial and boundary conditions of a generating PDE, that can be collected into a low-dimensional vector $x \in \reals[p]$. In what follows, it proves convenient to interchangeably refer to the vertices $x$ as these state vectors, since they fulfill this role on the discrete manifold that is the graph. Additionally, functionals, $u \in \reals[s]$, are obtained from a high-fidelity computation of dimension $P$, possibly a finite-dimensional version of a PDE, where $p+s \ll P$. Each state $\{x_{},u(x)\}$ is represented as a vertex $x \in V$ on the graph. Transitions between these states, including changes in system parameters, or steps of the numerical solver, allow for connections to be made between vertices on the graph. These transitions therefore can be thought of as inducing edges $e = (x_{},y_{})$ between states $x_{}$ and $y_{}$.

\subsection{Non-local calculus on finite, weighted graphs} \label{sec:graphtheory_nonlocalcalculus}
As outlined in the Introduction, we seek to develop reduced-order models in the form of differential equations for the evolution of some component of the state vector on the graph, or as functional representations. In the latter case, we explore forms that are motivated by the Taylor Series of local, differential calculus. In both cases, we need to first define a calculus on the discrete manifold that is the graph. We adopt a discrete, non-local calculus on finite, weighted graphs. \cite{Hein2007,Gilboa2008,Elmoataz2008,Desquesnes2013,Lozes2014}

Gilboa \etal define a discrete calculus, consisting of non-local operators, based on differences between states $x_{},y_{} \in V$, and an edge weight $w(x_{},y_{})$.\cite{Gilboa2008} The edge weights $w$ can be defined in a number of ways, and $w(x_{},y_{}) = w(x_{} - y_{})$ is natural for our purposes. As will be discussed, several constraints are placed on the forms of the weight functions, to ensure that the non-local calculus converges to local differential operators in certain limits. 

To fix ideas, we consider, without loss of generality of the development, that $s = 1$. Then, the $u(x_{})$ are scalar functionals at a single vertex. Vectors $v(x_{},y_{})$ are functions on vertex pairs $x,y$ or the edge $e = (x,y)$. Accordingly, we have a vector space $\mathnormal{L}$ such that there is a mapping $V\times V \mapsto \mathnormal{L}$ or $E\mapsto \mathnormal{L}$. Integration on the discrete manifold $G$ is a sum over vertices $x \in V$, and is relevant to operations such as inner products, which we define to be normalized by the size of the graph, $\abs{V} = n$. We draw attention to the distinction between the state vectors $\{x,u\} \in \reals[p+s]$ represented as vertices $x \in V$, and the edge vectors $v(x_{},y_{}) \in \mathnormal{L}$. These are distinct vector spaces, and are not copies of each other. We also note that for this treatment, there are no self edges and so vectors $v(x_{},y_{})$ at vertex $x$ are defined for $y \neq x$.

\subsubsection{Non-local calculus definitions} \label{sec:graphtheory_nonlocalcalculus_definitions}
The non-local gradient operator $\gradient[][][{u}] = \gradient[][{w}][[{u}]](x_{},y_{}) $, with respect to the weight $w$, is defined as the vector functional of scalars $u$ at state $x_{}$. It represents the vector of weighted differences with \textit{all other vertices} $y \in V$:
\begin{equation}
	 \gradient[][{w}][[{u}]](x_{},y_{}) \equiv [u(y_{})-u(x_{})]\sqrt{w(x_{},y_{})},\quad \forall ~ y \neq x \in V.
\end{equation}
The inner product between scalars $u_{\alpha}$ and $u_{\beta}$ is defined as the sum over the vertices:
\begin{align}
	\langle u_{\alpha},u_{\beta}\rangle \equiv&~ \frac{1}{n}\sum_{x\in V} u_{\alpha}(x_{})u_{\beta}(x_{}),
\end{align}
and the contraction between vectors $v_{\alpha}$ and $v_{\beta}$ at vertex $x$ is defined as the sum over disparate vertices:
\begin{align}
	[\dotproduct[{v_{\alpha}}][{v_{\beta}}]](x_{}) \equiv&~ \frac{1}{n-1}\sum_{y \in V\setminus \{x\}} v_{\alpha}(x_{},y_{})v_{\beta}(x_{},y_{}).
\end{align}
Inner products, and norms between pairs of vectors are likewise defined. An adjoint relation is obtained between the gradient and the divergence operators, as is a divergence theorem. This leads to a self-adjoint Laplacian. We do not reproduce these results here, but direct the interested reader to other work. \cite{Hein2007,Gilboa2008,Elmoataz2008,Desquesnes2013,Lozes2014} We work with the partial derivatives of $u(x_{})$ with respect to $x^{\mu}$, the $\ith[\mu]$ component of the state $x$. We denote the partial derivatives of local, differential calculus with $\partial$, and non-local partial derivatives with $\delta$. The latter are obtained as the contraction between the function gradient and the specifically defined unit vector,
\begin{align}
	\difference[1]{u(x_{})}{x^{\mu}} =&~ [\dotproduct[{\gradient[][{w}][{[{u}](x_{})}]}][{\hat{x}^{\mu}_{}}]] \label{eq:diff_1} \\
	=&~ \frac{1}{n-1}\sum_{y \in V\setminus \{x\}} [u(y_{}) - u(x_{})](y^{\mu}_{} - x^{\mu}_{})w(x_{},y_{}) \label{eq:diff_1_sum}, 
\end{align}
where the unit vector at state $x$ in the $\mu$ direction is defined as a vector gradient between states $x$ and $y$:
\begin{equation}
	\hat{x}^{\mu}(x_{},y_{}) \equiv \gradient[][{w}][[x^{\mu}]](x_{},y_{}) = [y^{\mu}_{} - x^{\mu}_{}]\sqrt{w(x_{},y_{})}.
\end{equation}
The weights $w(x_{},y_{})$ are chosen such that the unit vectors are normalized:
\begin{equation}
 	[\dotproduct[{\hat{x}^{\mu}}][{\hat{x}^{\nu}}]] = \delta^{\mu\nu}, \label{eq:unitnorm}
\end{equation}
where $\delta^{\mu\nu}$ is the Kronecker delta. As will be discussed, this constraint together with the inner products also being normalized by the volume of the space, leads to other important properties including consistency of these non-local derivatives with their counterparts in differential calculus.

\noindent Higher order derivatives also can be computed as approximations to $\nderivative[k+1]{u(x)}{x^{\mu_{0}}}{x^{\mu_{k}}}$ using an extension of \cref{eq:diff_1} and can be defined recursively:
\begin{align}
\ndifference[k+1]{u(x_{})}{x^{\mu_{0}}}{x^{\mu_{k}}} ~=&~ 
[\dotproduct[{\gradient[][{w}][{[ 
	{\dotproduct[{\gradient[][{w}][{[
	{\cdots[\dotproduct[{\gradient[][{w}][{[{\dotproduct[{\gradient[][{w}][{[u(x_{})]}]}][{\hat{x}^{\mu_{0}}}] }]}]}][{\hat{x}^{\mu_{1}}}]]\cdots}]}]}][{\hat{x}^{\mu_{k-1}}}]}]}]}][{\hat{x}^{\mu_{k}}}]] \label{eq:diff_p}
\end{align}
It will now be shown, that under specific constraints on the weight functions, the above partial derivatives are, to leading order, the corresponding partial derivatives of differential calculus. However certain properties of differential calculus, such as the commutativity of mixed partial derivatives, do not automatically hold to all orders of approximation with all possible choices of edge weights. 

\subsection{Weights defined on neighborhoods of a vertex}
\label{sec:localweights}
\subsubsection{Local weight definitions for first order derivatives in one dimension}\label{sec:nonlocalerror_localderivative}
Within this non-local calculus formalism for the partial derivatives, we are free to choose the edge weight definitions $w$, and are solely constrained by desiring a certain order of accuracy between the non-local derivatives and their differential derivative counterparts. For example if the normalization condition in \cref{eq:unitnorm} is imposed, then derivatives may have up to first order accuracy. We desire to define non-local derivatives with any order of accuracy for any order of derivative, to allow for general application, and propose the following weight definitions with finite support over a local neighborhood of vertices around a vertex of interest. 

\noindent Given a point $\widetilde{x}$ and a neighborhood $\mathcal{N}(\widetilde{x})$ of nearby points, we seek a certain order of accuracy between the non-local and differential calculus derivatives, and find the $w(\widetilde{x},x)$ with support over the neighborhood $\mathcal{N}(\widetilde{x})$ that yields this order of accuracy. This neighborhood is defined as a set of vertices within a subset of the whole graph of vertices $\mathcal{N}(\widetilde{x}) \subseteq \widetilde{V} \subseteq V$ and has size $d = d(\widetilde{x}) = \abs{\mathcal{N}(\widetilde{x})}$. To quantify the accuracy of our derivative definitions, we consider various components of the errors. We first define  monomials for the $d = \abs{\mathcal{N}(\widetilde{x})}$ data points included in the neighborhood in terms of
\begin{align}
	z = z(\widetilde{x}) =&~ x - \widetilde{x} \in \reals[d \times p] \label{eq:zmonomial}.
\end{align}
For brevity we will use the notation of the sum over these monomials for a fixed $\widetilde{x}$ solely with the neighborhoods $\sum_{x \in \mathcal{N}(\widetilde{x})} \equiv \sum_{\mathcal{N}(\widetilde{x})}$. 

\noindent Of central interest is the error of non-local $l$-order derivatives at a local point $\widetilde{x} \in \widetilde{V}$
\begin{align}
	\varepsilon_{l}(\widetilde{x}) \equiv \unindifference[l]{u(\widetilde{x})}{x} - \uninderivative[l]{u(\widetilde{x})}{x}, \label{eq:error_local_derivative}
\end{align}
as well as the global derivative error over all local points
\begin{align}
	\varepsilon_{l} \equiv \frac{1}{\abs{\widetilde{V}}}\sum_{\widetilde{x} \in \widetilde{V}} \varepsilon_{l}(\widetilde{x}). \label{eq:error_global_derivative}
\end{align}
The generalization to higher dimensions for errors $\varepsilon_l^{\mu\nu\cdots}$ of derivatives along dimensions $\mu,\nu,\cdots$ is straightforward. We also introduce the local error of a model $u(x|\widetilde{x})$ for a function $u(x)$ at a local point $x \in V \setminus \widetilde{V}$, based around a point $\widetilde{x} \in \widetilde{V}$
\begin{align}
	e(x|\widetilde{x}) \equiv u(x|\widetilde{x}) - u(x), \label{eq:error_local_model}
\end{align}
and the global model error over all local points
\begin{align}
	e \equiv \frac{1}{\abs{V \setminus \widetilde{V}}}\sum_{x \in V \setminus \widetilde{V}} e(x|\widetilde{x}). \label{eq:error_global_model}
\end{align}
As will be discussed, in general we choose the model points $\widetilde{x}$ for a given evaluation point $x$ as the point $\widetilde{x} \in \widetilde{V}$ that is closest to $x$, with respect to a chosen metric $\norm{\cdot}$ over the $p$ dimensional space. Given the dependence of the weights $w(z)$, we can define these errors in terms of this distance and the components of $z$ have norm $\norm{z_i} \geq h$, $i = 1 \cdots d$, for a length scale $h$ that depends on the neighborhood. The aim of our analysis is to show the convergence properties of these errors as functions of the length scale $h$, and in particular in the limit of a continuous graph as $h \to 0$.

\noindent The analysis will not impose any symmetries and each $w(x-\widetilde{x}) ~\forall x \in \mathcal{N}(\widetilde{x})$ will be found to ensure each derivative in the graph has a specified order of accuracy. For $l = \{1,\dots,k\}$ order derivatives, by assuming an adequately continuous function $u(x)$ and expanding the definition of the non-local derivatives in \cref{eq:localderivative_a} as a Taylor series about the model point $\widetilde{x}$, we find there are $q_l$ constraints to be imposed on the derivative weights if we desire that the error scales with length scale $h$ as a power $r_l$. These $q_l$ constraints come from requiring the lower order $l^{\prime} < l$ terms in the Taylor series of the error to be identically zero, ensuring that the leading order error term is of order $r_l$ and so
\begin{align}
	\unindifference[l]{u(\widetilde{x})}{x} =&~ \uninderivative[l]{u(\widetilde{x})}{x} + \mathcal{O}(h^{r_{l}}).
\end{align}
In $p$ dimensions, the mixed derivatives for a given total order of derivative $l$, and $q_l$ can be shown in \cref{app:nonlocalerror_number_terms} to grow exponentially with $p$ and $r$, and there are
\begin{align}
	q =&~ q(p,r) = \binom{p+r}{r} - 1 \label{eq:pqr_constraints}
\end{align}
unique constraints when accounting for mixed derivatives commuting. Depending on the definitions of the derivative neighborhoods via $\mathcal{N}({\widetilde{x}})$ and the behavior of this error $\varepsilon_l(\widetilde{x})$, there will emerge a relationship found between derivative order $l$, the model order $k$, the number of dimensions $p$, the number of constraints $q_l$, and the final scaling $r_{l}$.

\noindent We will first derive the weight definitions in $p=1$ dimensions for first order derivatives, and then describe how the definitions generalize to higher order derivatives and higher dimensions. Given a point $\widetilde{x}$ and a neighborhood $\mathcal{N}(\widetilde{x})$, in this $p=1$ dimensional setting, we may evaluate non-local derivatives of the form
\begin{align}
	\difference[1]{u(\widetilde{x})}{x} =&~ \frac{1}{\abs{\mathcal{N}(\widetilde{x})}} \sum_{x \in \mathcal{N}(\widetilde{x})} \left( u(x) - u(\widetilde{x})\right)(x - \widetilde{x})w(\widetilde{x},x). \label{eq:localderivative}
\end{align}

\noindent Guided by stencils that are generated by discretization methods, we define weights to be of the form
\begin{align}
	w(\widetilde{x},x) =&~ \frac{\abs{\mathcal{N}(\widetilde{x})}}{(x - \widetilde{x})^2}a(x-\widetilde{x}) \label{eq:localweightfunction}
\end{align}
where $a = a(x - \widetilde{x})$ are (dimensionless) reduced weight functions that are independent of vertex spacing, $\norm{x - \widetilde{x}}$ and not necessarily symmetric in $x$ and $\widetilde{x}$. The derivatives therefore take the form: 
\begin{align}
	\difference[1]{u(\widetilde{x})}{x} =&~ \sum_{x \in \mathcal{N}(\widetilde{x})} \frac{u(x) - u(\widetilde{x})}{(x - \widetilde{x})}a(x-\widetilde{x}). \label{eq:localderivative_a}
\end{align}
The reduced weights are strictly functions of the monomials in $z = x - \widetilde{x} \sim h$, and represent the weight of each point $x$ in the stencil of points around $\widetilde{x}$. For example, a symmetric difference two-point stencil in $p=1$ dimensions would have $a(x-\widetilde{x}) = \frac{1}{2}$. 

\noindent In previous works involving non-local calculus on graphs, \cite{Gilboa2008,Elmoataz2008} continuous weights $w(\widetilde{x},x)$ are chosen, that have support over the entire graph. These weights generally are chosen to decay away from $\widetilde{x}$, such as a Gaussian weight $w \sim e^{-\abs{x-\widetilde{x}}^2}$, with the conjecture that faster than polynomial decay will lead to the non-local derivatives converging to the differential derivatives. We have conducted detailed analyses and numerical studies with these continuous weights \cite{Duschenes2021}, and have shown that in fact weights with non-local support over neighborhoods $\mathcal{N}$ such that $\abs{\mathcal{N}(\widetilde{x})} \approx n$ (the dimension of $G$), lead to a constant error between the derivative definitions, for any form of the weights. Local neighborhoods of weights with strictly finite support, $\abs{\mathcal{N}(\widetilde{x})} \ll n$ are therefore required for rigorous convergence for any distribution of vertices on the graph.

\begin{figure}[hpt]
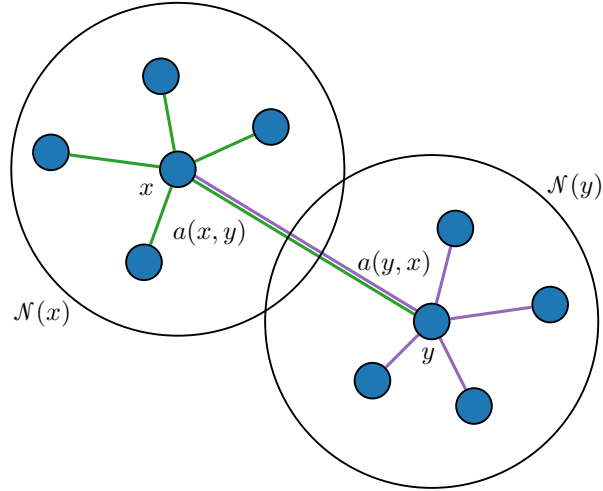

\centering
\tikzfig{0.9}{figures/localweights}
\caption{Local neighborhoods of adjacent vertices $x$ and $y$, where the absence of symmetry in the weight functions for the edge weight between $x$ and $y$ is shown.}
\label{fig:localweights}
\end{figure}

\noindent Depending on the local neighborhoods, as per \cref{fig:localweights}, the reduced weights are not necessarily symmetric: $a(x-y) \neq a(y-x)$, $x,y \in \widetilde{V}$. The graph of vertices can therefore be thought of as being potentially directed and multi-edged, where pairs of vertices $x$ and $y$ may have an edge $a(y-x)$ from $x$ to $y$, and an edge $a(x-y)$ from $y$ to $x$. Depending on the state vectors of the graph, such as a uniformly spaced mesh of data, or various boundary conditions on the data, there may be additional symmetries that yield a symmetric, undirected graph where $a(y-x) = a(x-y)$. 

\noindent \textbf{Remark}: As an aside, a modified $k$-order Taylor series model may be constructed using a training dataset as a functional representation for $u(x) \approx u_k(x|\widetilde{x})$, where a different Taylor series expansion is developed around each possible point $\widetilde{x}$. This modified Taylor series has an identical form to that of the standard Taylor series, with non-local derivatives replacing differential derivatives, and there being additional linear coefficients associated with each term that are fit with regression approaches to account for the error in the non-local derivatives. A complete error analysis for this modified Taylor series approach can be referred to in \cref{app:error_analysis_model}, including an analysis of the number of constraints necessary to guarantee an order of accuracy $r$ in \cref{app:error_analysis_local}, an analysis of the $l$-order non-local derivatives' convergence in $p$ dimensions in \cref{app:nonlocalerror_localderivative_p_higher}, and an analysis of the convergence of this modified Taylor series to the standard Taylor series in \cref{app:error_analysis_coef} in $p=1$ dimensions and \cref{app:error_analysis_coef_p} in general $p$ dimensions. These analyses reveal the expected order of accuracy: a $\ith[k]$ order (modified) Taylor series has error scaling as $h^{k+1}$.

\noindent To show this scaling relation of the derivative error, we may expand the differences of functions in the non-local derivative definitions in a Taylor series about $\widetilde{x}$:
\begin{align}
	\difference[1]{u(\widetilde{x})}{x} =&~ \sum_{s=0}^{\infty}\frac{1}{(s+1)!}\uninderivative[s+1]{u(\widetilde{x})}{x} \sum_{x \in \mathcal{N}(\widetilde{x})} (x-\widetilde{x})^{s}a(x - \widetilde{x}).
\end{align}
where the sum over $\mathcal{N}(\widetilde{x})$ comes from the weights and \cref{eq:localderivative_a}. The non-local derivative Taylor series expansion can be rewritten as a function of the monomials of $z$ in \cref{eq:zmonomial}
\begin{align}
	\difference[1]{u(\widetilde{x})}{x} =&~ \sum_{s=0}^{\infty}\frac{1}{(s+1)!}\uninderivative[s+1]{u(\widetilde{x})}{x} \sum_{\mathcal{N}(\widetilde{x})} z(\widetilde{x})^{s}a(z(\widetilde{x})),
\end{align}
which can be written succinctly as a product of two infinite dimensional matrices
\begin{align}
	\difference[1]{u(\widetilde{x})}{x} =&~ d_{\infty}^T V_{\infty}^T a .
\end{align}
Here, $V_{\infty} = V_{\infty}(z) \in \reals[d \times \infty]$ is the Vandermonde-like matrix of powers of $z$ with elements
\begin{align}
	V_{\infty_l}^T =&~ z^{l-1} \in \reals[d],
\end{align}
$d_{\infty} = d_{\infty}(\widetilde{x}) \in \reals[\infty]$ is the vector of derivatives with elements
\begin{align}
	d_{\infty_l} =&~ \frac{1}{l!}\uninderivative[l]{u(\widetilde{x})}{x} \in \reals[],
\end{align}
and $a = a(z) \in \reals[d]$ is the vector of reduced weights with elements
\begin{align}
	a_{l} = a(z_l).
\end{align}
Here we have used base $1$ indexing to be consistent with previous definitions of these matrices.

\noindent If we desire that the non-local derivatives are $r_{1} = r$ order accurate: 
\begin{align}
	\difference[1]{u(\widetilde{x})}{x} =&~ \derivative[1]{u(\widetilde{x})}{x} + \sum_{s=r}^{\infty}\frac{1}{(s+1)!}\uninderivative[s+1]{u(\widetilde{x})}{x} \sum_{\mathcal{N}(\widetilde{x})} z(\widetilde{x})^{s}a(z(\widetilde{x})), \\
	=&~ \derivative[1]{u(\widetilde{x})}{x} + \mathcal{O}(z^{r}),
\end{align}
then the weights can be found from solving the linear system of equations of the first $r$ moments:
\begin{align}
	\sum_{\mathcal{N}(\widetilde{x})} z(\widetilde{x})^{s}a(z(\widetilde{x})) = \delta_{0s} \quad s = \{0,\dots,r-1\} \label{eq:localweightconstraint_moments}.
\end{align}
In $p=1$ dimensions, the number of constraints equals the number of lower order derivative terms, since there are no mixed derivatives and one derivative for each order of derivative. Therefore first order derivatives have 
\begin{align}
	q_1 =&~ r_1 = r
\end{align}
constraints.

\noindent This can be written as a linear problem, given we partition the matrices as
\begin{align}
	V_{\infty} =&~ \left[\begin{array}{ccc} V & \bar{V} \end{array} \right],
	\intertext{where $V = \left[1~Z\right] \in \reals[d \times r],Z \in \reals[d \times r-1], \bar{V} \in \reals[d \times \infty]$ are Vandermonde matrices,}
	V =&~	\left[ 
		\begin{array}{ccccc} 
		1 & z & z^2 & \cdots & z^{r-1}
		\end{array} 
		\right],\\
	Z =&~	\left[ 
		\begin{array}{cccc} 
		z & z^2 & \cdots & z^{r-1}
		\end{array} 
		\right],\\
	\intertext{and}
	d_{\infty} =&~ \left[\begin{array}{c} d \\ \bar{d} \end{array} \right],
\end{align}
where $d = \left[d_1 ~ \hat{d}^T \right]^T \in \reals[r]$, $\hat{d} \in \reals[r-1]$, and $\bar{d} \in \reals[\infty]$.
The linear problem to be solved is
\begin{align}
	d^TV^Ta =&~ d_1 \\
	\intertext{which can be written as a linear combination of the derivatives}
	d_1(1^Ta - 1) ~+&~ \hat{d}^TZ^Ta = 0.
	\label{eq:lincombderivatives}
\end{align}
The weights must yield consistent derivatives for all functions $u$, and therefore for all $d\in \reals[r]$. The problem can thus be simplified to
\begin{align}
	V^Ta =&~ e_1 \label{eq:linear_localweights}
\end{align}
where $e_1 \in \reals[r]$ has elements $e_{1_l} = \delta_{1l}$.

\noindent When $p = 1$, the dimension of the matrix $d = r$ and the size of the local neighborhood, the number of vertices in $\mathcal{N}(\widetilde{x})$, equals the order of accuracy, the Vandermonde matrix is full rank,\cite{Turner1966,Pugliese2000} with the known ordinary least squares pseudo-inverse solution
\begin{align}
	a = (VV^T)^{-1}Ve_1.
\end{align}

\noindent Given this linear problem for $a(x-\widetilde{x})$ in \cref{eq:linear_localweights}, the weight functions for all expansion points $\widetilde{x}$ for the modified Taylor series models can be found using \cref{eq:localweightfunction}. The weights are then used in the original non-local calculus definitions in \cref{eq:localderivative}, yielding non-local first derivatives that are $r$ order accurate. This procedure can be repeated for higher order derivatives, and as will be shown in \cref{app:nonlocalerror_localderivative_higher} of the \appsecname, distinct sets of edge weights will be found to ensure that each derivative at each expansion point has the desired order of accuracy.


\subsubsection{Local weight definitions for higher order derivatives in higher dimensions} \label{sec:nonlocalerror_localderivative_higher}
Given the definitions of the first derivatives in \cref{eq:localderivative} for $p=1$, the dimensionality of $x$, we may take several approaches to the form of the non-local $l>1$ higher derivatives, particularly in $p>1$ dimensions. We choose to write higher order derivatives recursively from the first derivatives
\begin{align}
	\difference[1]{u(\widetilde{x})}{x^{\mu}} =&~ \sum_{\mathcal{N}^{{\mu}}(\widetilde{x})} \frac{u(x) - u(\widetilde{x})}{x^{\mu} - \widetilde{x}^{\mu}} {a^{\mu}}(x-\widetilde{x}),\quad \mu \in \{0,\dots,p-1\}. \label{eq:localderivative_p}
\end{align}
Here $\mathcal{N}^\mu$ is the neighborhood used for the $x^\mu$ coordinate, $\mu \in \{0,\dots,p-1\}$. From analysis that can be found in \cref{app:nonlocalerror_localderivative_p_higher} of the \appsecname, we find that due to this recursive nature of the derivatives, higher order derivatives will involve nested sums over different neighborhoods $\mathcal{N}(\widetilde{x}),\mathcal{N}(\widetilde{x}^{\prime})$. An important distinction is that unless there is strict isotropy and homogeneity in the graph, the sets of monomials 
\begin{align}
	\{z(\widetilde{x},x) : x \in \mathcal{N}(\widetilde{x})\} \neq&~ \{z(\widetilde{x}^{\prime},x) : x \in \mathcal{N}(\widetilde{x}^{\prime})\} \quad ~\forall~ \widetilde{x} \neq \widetilde{x}^{\prime} \in \widetilde{V}, 
\end{align}
and so different neighborhoods will not be equivalent about their base points. Performing an identical analysis of expanding the definitions of the non-local derivatives in a standard Taylor series about the model points $\widetilde{x}$, lower order error terms in higher order derivatives will not cancel in the case of general unstructured data, leading to the order of accuracy of higher order derivatives depending on the order of accuracy of all lower derivatives. Therefore we must choose our first order derivative weights to be adequately accurate to retain a desired order of accuracy of higher order derivatives.

\noindent Therefore for each subsequent derivative along a dimension, the weights are used from the fixed set of computed weights of the first derivatives $\{a^{\mu}\}$ and there is a single set of $p$ weights for all derivatives at a given point. The weights obey the general constraints for all $s = \{0,\dots,r-1\}$:
\begin{align}
	\sum_{\mathcal{N}^{{\mu}}(\widetilde{x})} \frac{z^{\mu_0\cdots\mu_{s}}}{z^{\mu}}{a^{\mu}}(z) =&~ \delta_{0s}\delta^{\mu_0\mu} \cdots \delta^{\mu_{s}\mu}
	. \label{eq:localweightconstraint_moments_p_main}
\end{align}
and the Vandermonde matrix of constraints is generalized into a multi-dimensional version with all polynomials of $z^{\mu_0\cdots\mu_{s}}$. 

\noindent Given this approach to the weights with fixed order of first derivative accuracy $r$ that are computed from $q(p,r)$ constraints in \cref{eq:pqr_constraints}, higher order derivatives for all dimensions will scale with decreasing accuracy as
\begin{align}
	\unindifference[l]{u(\widetilde{x})}{x} =&~ \uninderivative[l]{u(\widetilde{x})}{x} + \mathcal{O}(z^{r+1-l}) \\
	\intertext{and so}
	r_{l} =&~ r + 1 - l \label{eq:derivative_scaling_local}.
\end{align}
Therefore $r \geq k$ should be chosen to ensure all derivatives have a non-zero order of accuracy for a fixed set of weights. The $\ith[l]$ order partial derivatives in higher dimensions therefore commute up to an $r+1-l$ order of accuracy, and this proposed stencil-based non-local calculus provides a rigorous basis of derivative operators to compute reduced-order models on unstructured data.

\noindent The neighborhood for estimation of each derivative is determined in such a way that the multi-dimensional Vandermonde matrix $V$ in $p$ dimensions remains full rank and non-singular. The selection is done by preemptively sorting points in terms of some defined metric, in this work, the minimum Euclidean distance $\norm{x-\widetilde{x}}$, and conditionally adding points to the neighborhood if they continue to allow a definition of the pseudo-inverse of $V$. Points are added until the size of the neighborhood is equal to the number of constraints and the neighborhood points satisfy the constraints given by \cref{eq:localweightconstraint_moments_p_main}.

\subsection{Numerical error of modified Taylor series model in $p$ dimensions}\label{sec:nonlocalerror_localderivative_results}
\begin{figure}[hpt]
\centering
\begin{subfigure}[t]{0.49\textwidth}
	\centering
	\includegraphics[width=\textwidth]{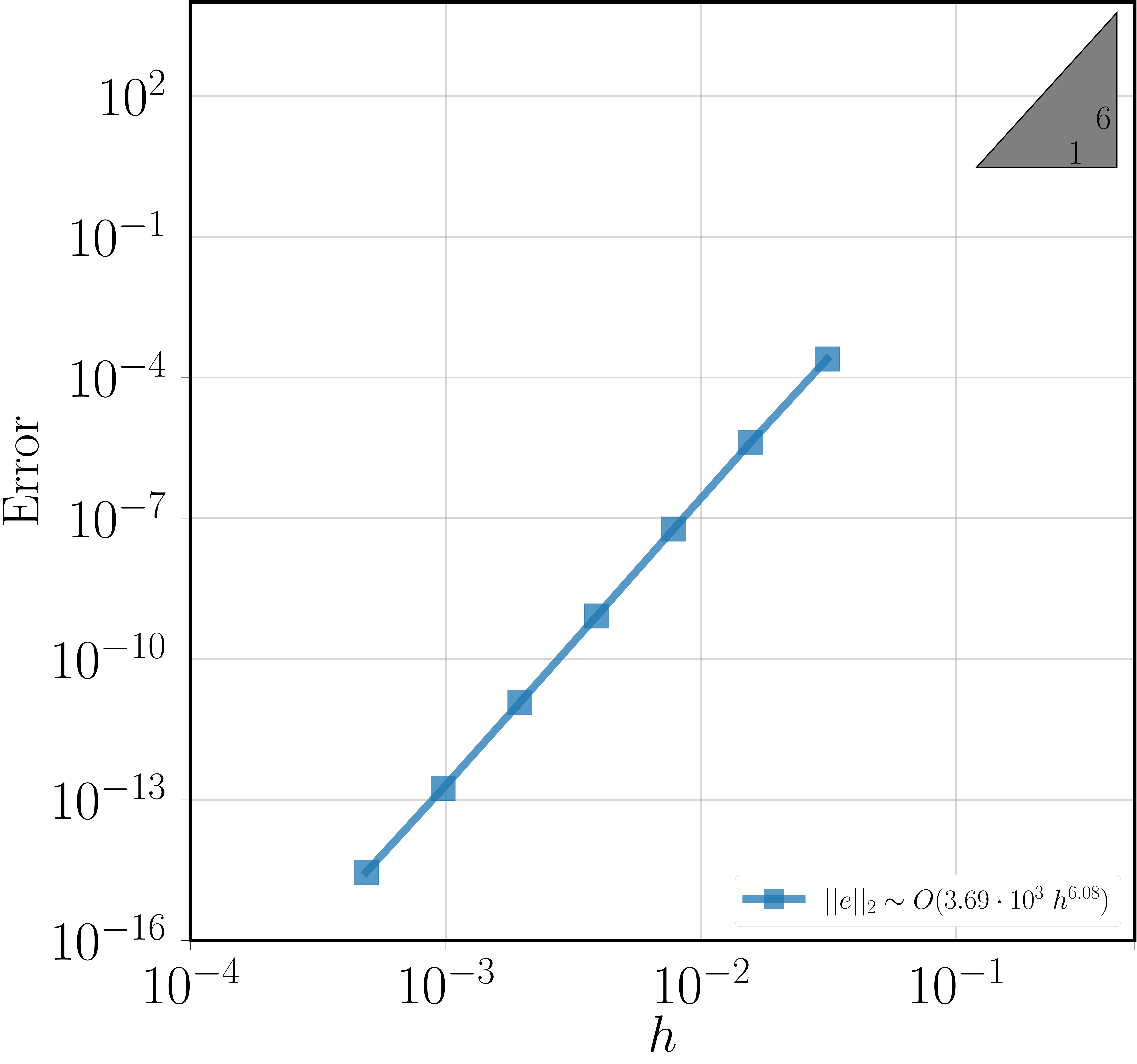}
	\subcaption{Global model error $e$.}
	\label{fig:errorscaling_1_error_global}
\end{subfigure}
\hfill
\begin{subfigure}[t]{0.49\textwidth}
	\centering
	\includegraphics[width=\textwidth]{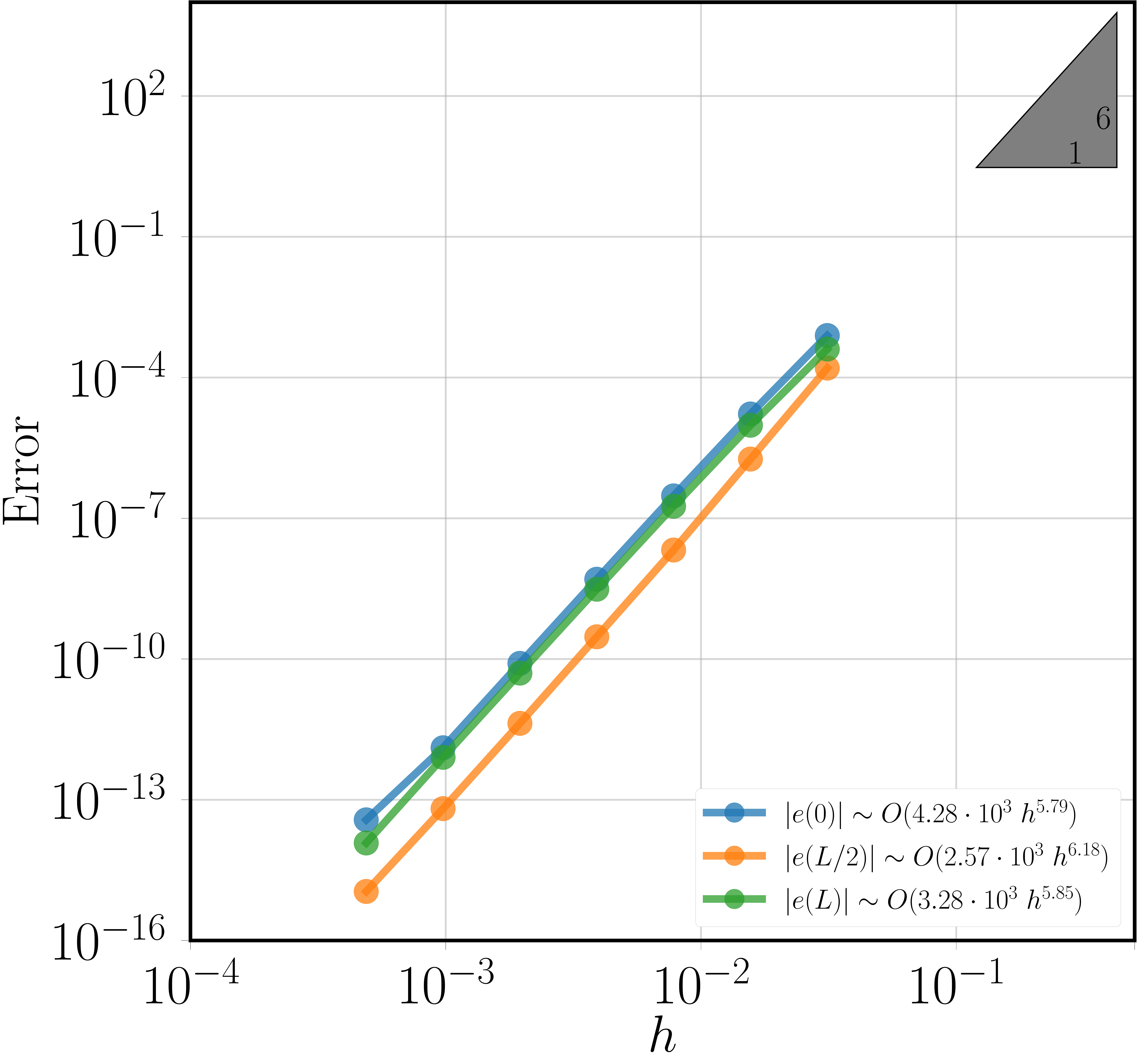}
	\subcaption{Local model error $e(x)$ at points in $(0,L)\times(0,L)$.}
	\label{fig:errorscaling_1_error_local}
\end{subfigure}\\
\vspace{1cm}
\begin{subfigure}[t]{0.49\textwidth}
	\centering
	\includegraphics[width=\textwidth]{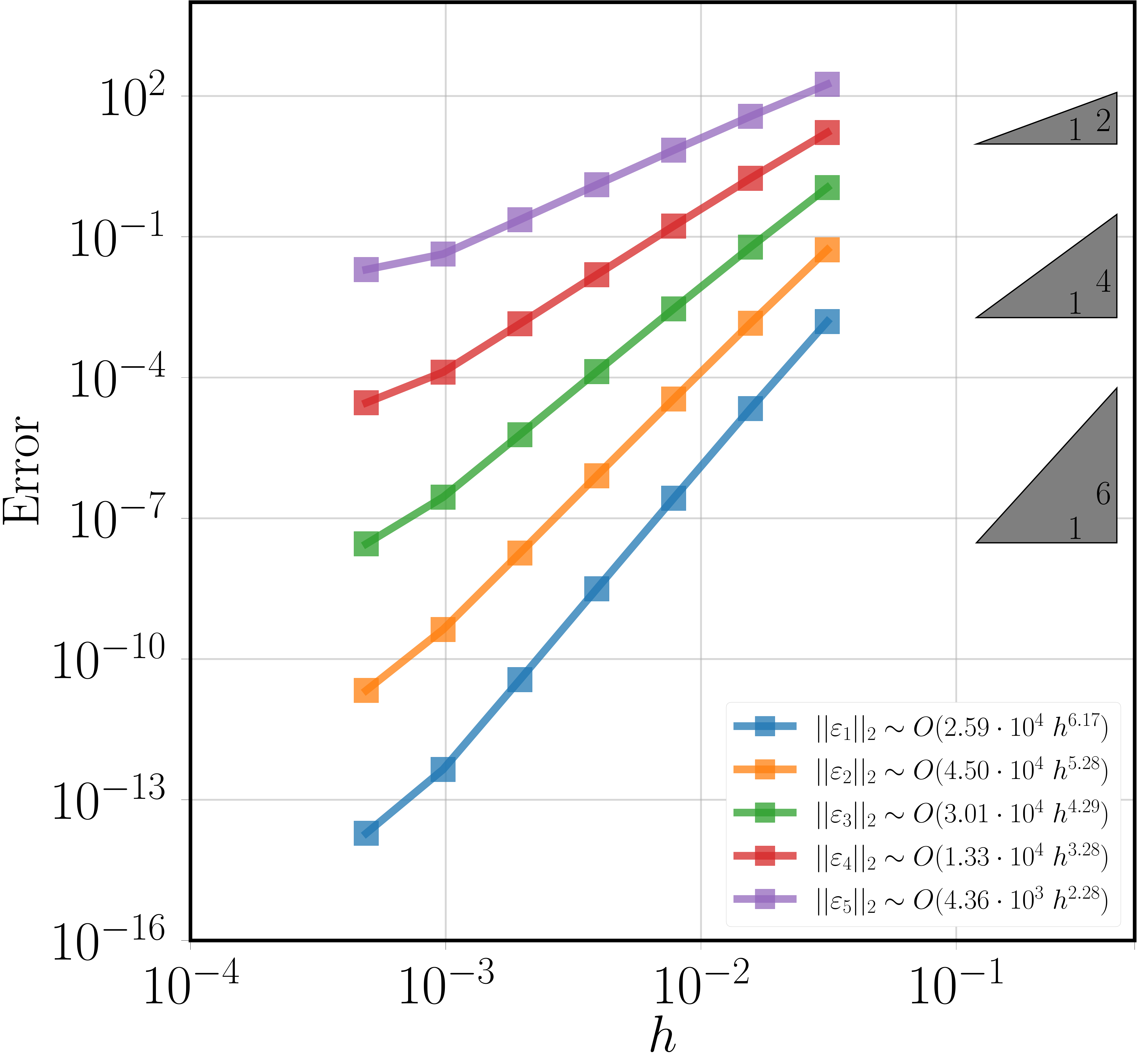}
	\subcaption{Global derivative error $\varepsilon_{l}$, $l = \{1 \cdots k\}$ for derivatives of order $l$ with respect to $x$.}
	\label{fig:errorscaling_1_derivativeerror_global}
\end{subfigure}
\hfill
\begin{subfigure}[t]{0.49\textwidth}
	\centering
	\includegraphics[width=\textwidth]{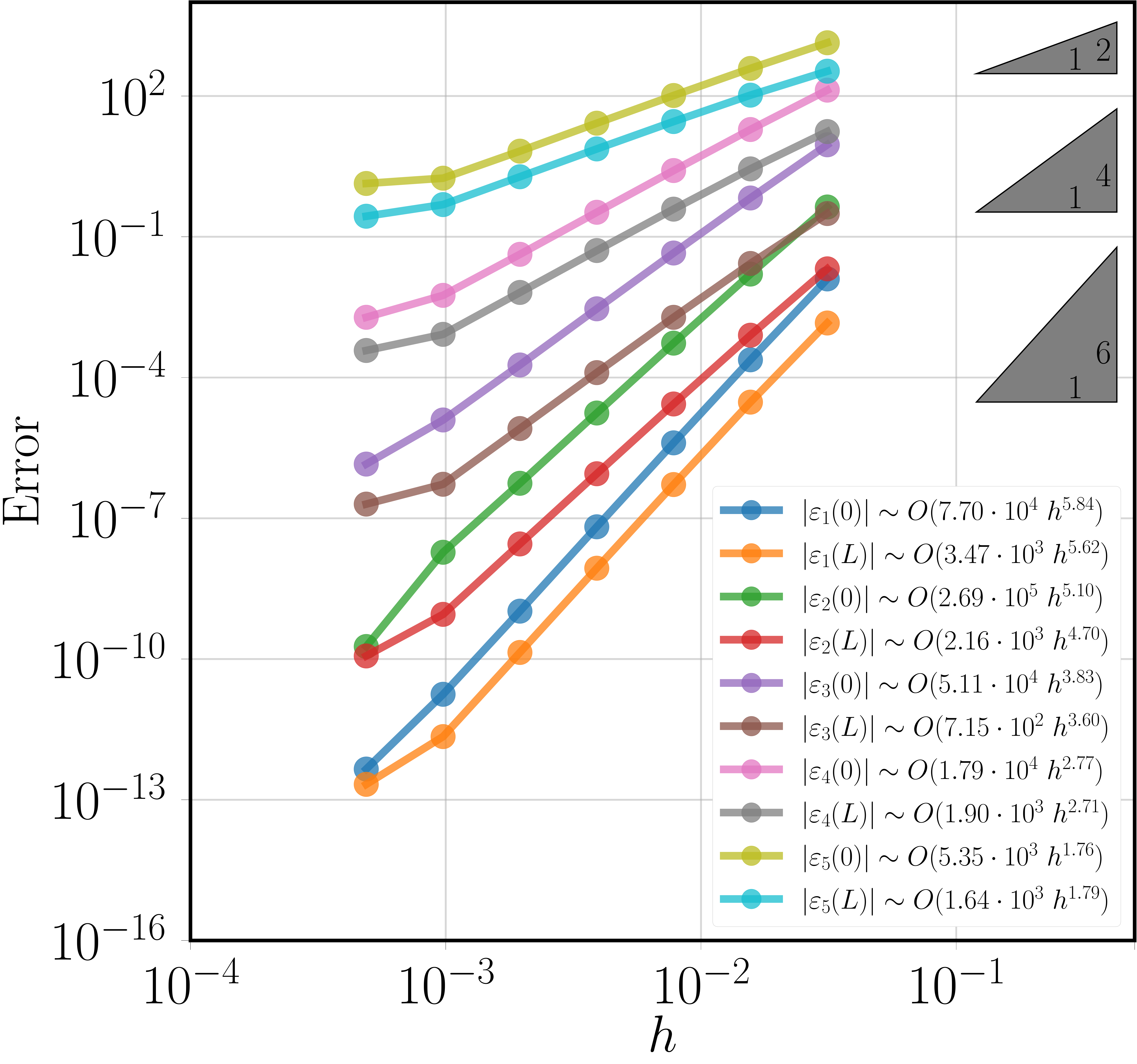}
	\subcaption{Local derivative error $\varepsilon_{l}(x)$, $l = \{1 \cdots k\}$ for derivatives of order $l$ with respect to $x$.}
	\label{fig:errorscaling_1_derivativeerror_local}
\end{subfigure}
\caption{Error scaling as a function of data length scale $h$, of a $r=k+1$ accurate, $k=5$-order Taylor series model for a $K=8$ order polynomial $u(x) = \sum_{\sum s \leq K}\alpha_{s} x^{s}$ in $p=1$ dimensions. Polynomial coefficients $\alpha_{s} \sim U[-1,1]$ are sampled from a uniform distribution and all unique commuting Taylor series derivatives are included in the model up to $k$ order terms. Inset triangles indicate theoretical slopes of scaling of fits with length scale $h$, where model error $e \sim \mathcal{O}(h^{k+1})$ and $\ith[l]$ order derivative error $\varepsilon_l \sim \mathcal{O}(h^{r+1-l})$.}
\label{fig:errorscaling_1}
\end{figure}

\begin{figure}[hpt]
\centering
\begin{subfigure}[t]{0.49\textwidth}
	\centering
	\includegraphics[width=\textwidth]{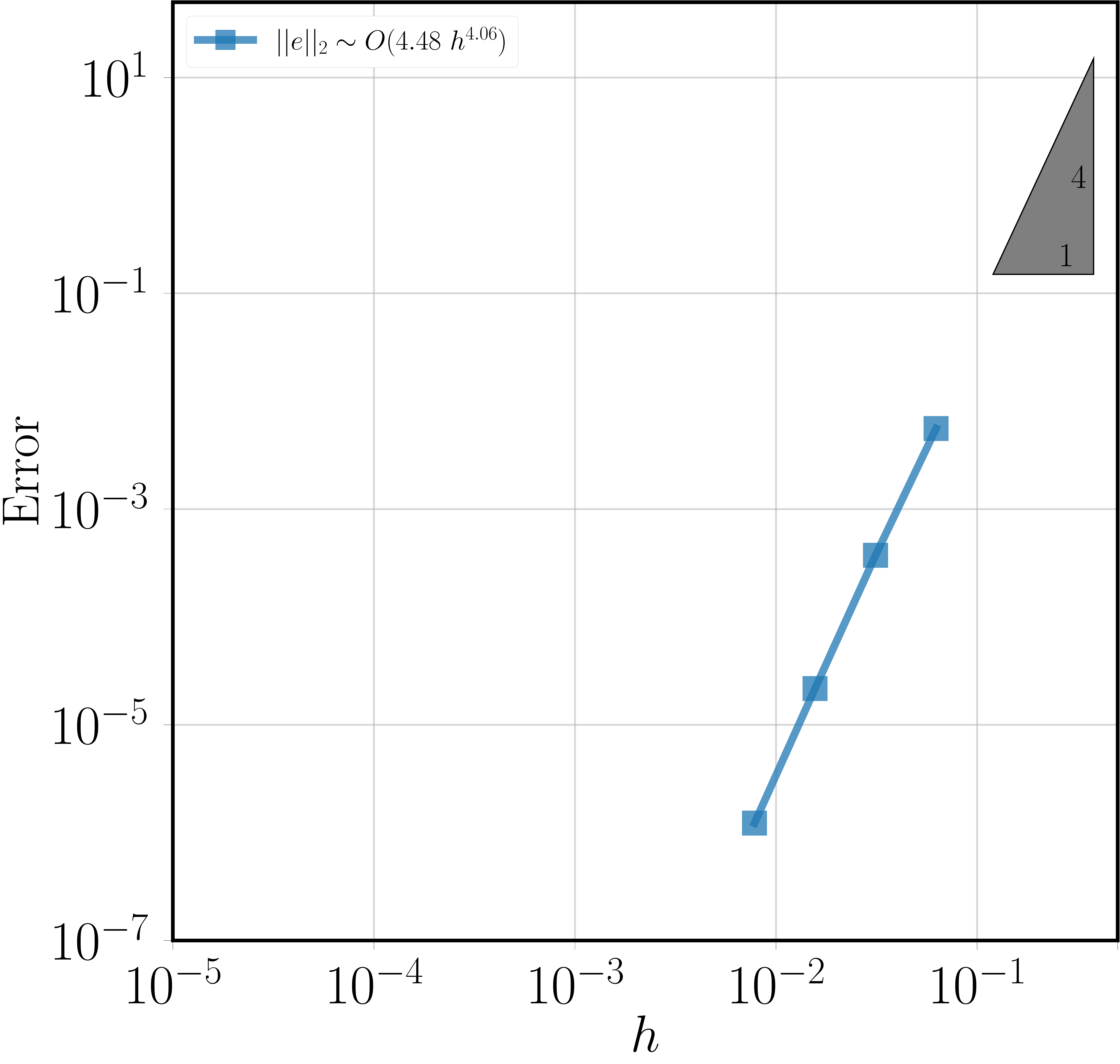}
	\subcaption{Global model error $e$.}
	\label{fig:errorscaling_2_error_global}
\end{subfigure}
\hfill
\begin{subfigure}[t]{0.49\textwidth}
	\centering
	\includegraphics[width=\textwidth]{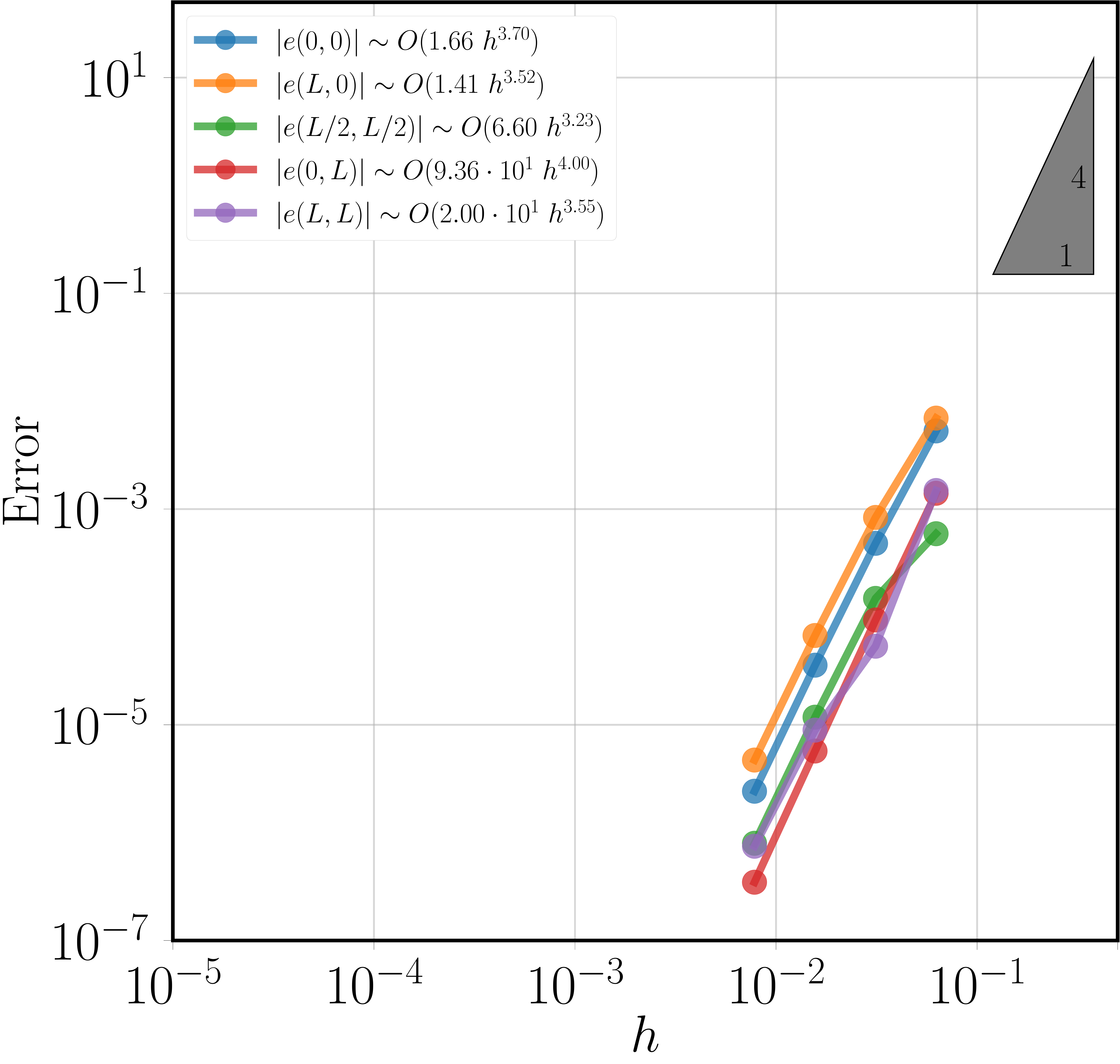}
	\subcaption{Local model error $e(x)$ at points in $(0,L)\times(0,L)$.}
	\label{fig:errorscaling_2_error_local}
\end{subfigure}\\
\vspace{1cm}
\begin{subfigure}[t]{0.49\textwidth}
	\centering
	\includegraphics[width=\textwidth]{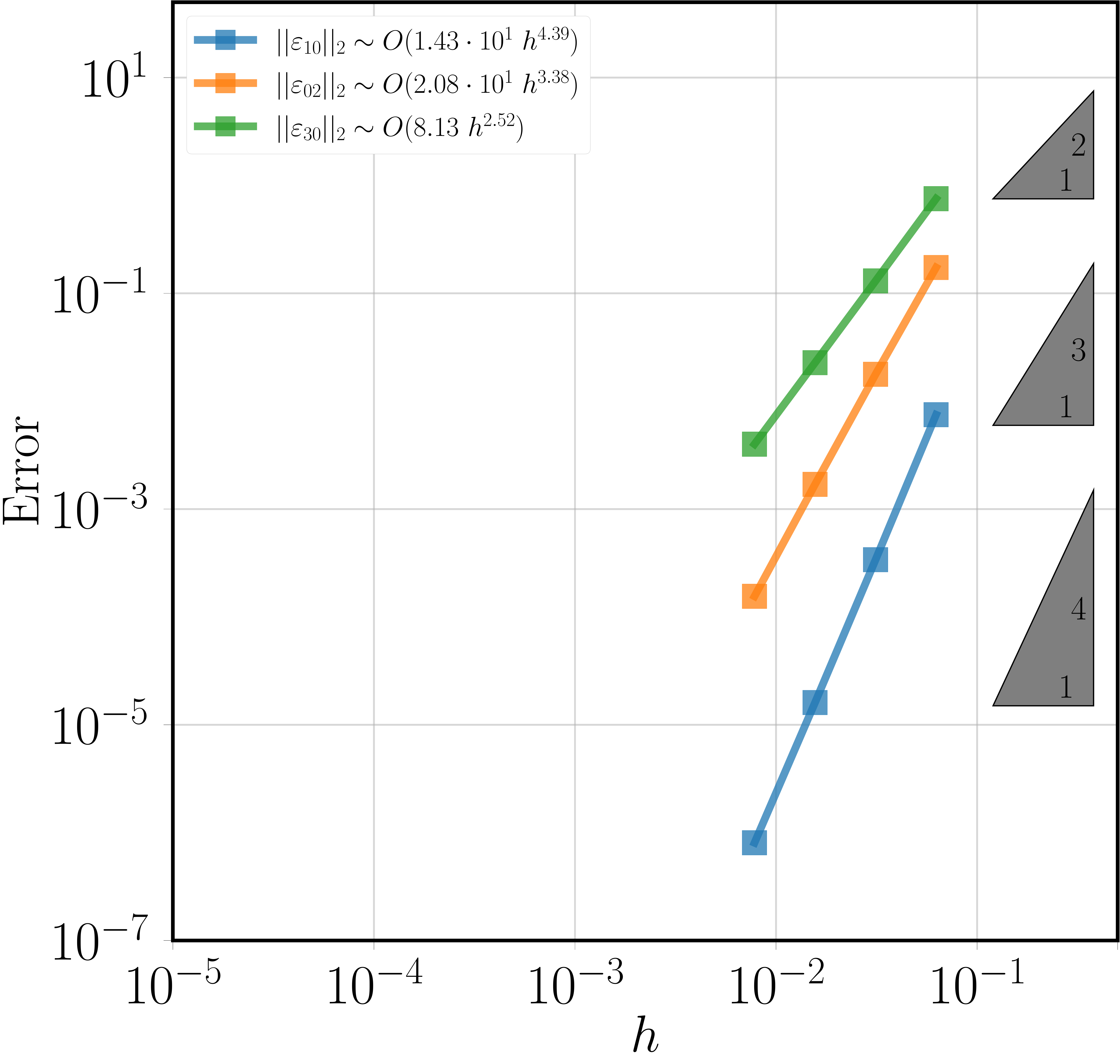}
	\subcaption{Global derivative error $\varepsilon_{lm}$, $l+m = \{1 \cdots k\}$ for derivatives of order $l,m$ with respect to $x = \{x^{0}, x^{1}\}$.}
	\label{fig:errorscaling_2_derivativeerror_global}
\end{subfigure}
\hfill
\begin{subfigure}[t]{0.49\textwidth}
	\centering
	\includegraphics[width=\textwidth]{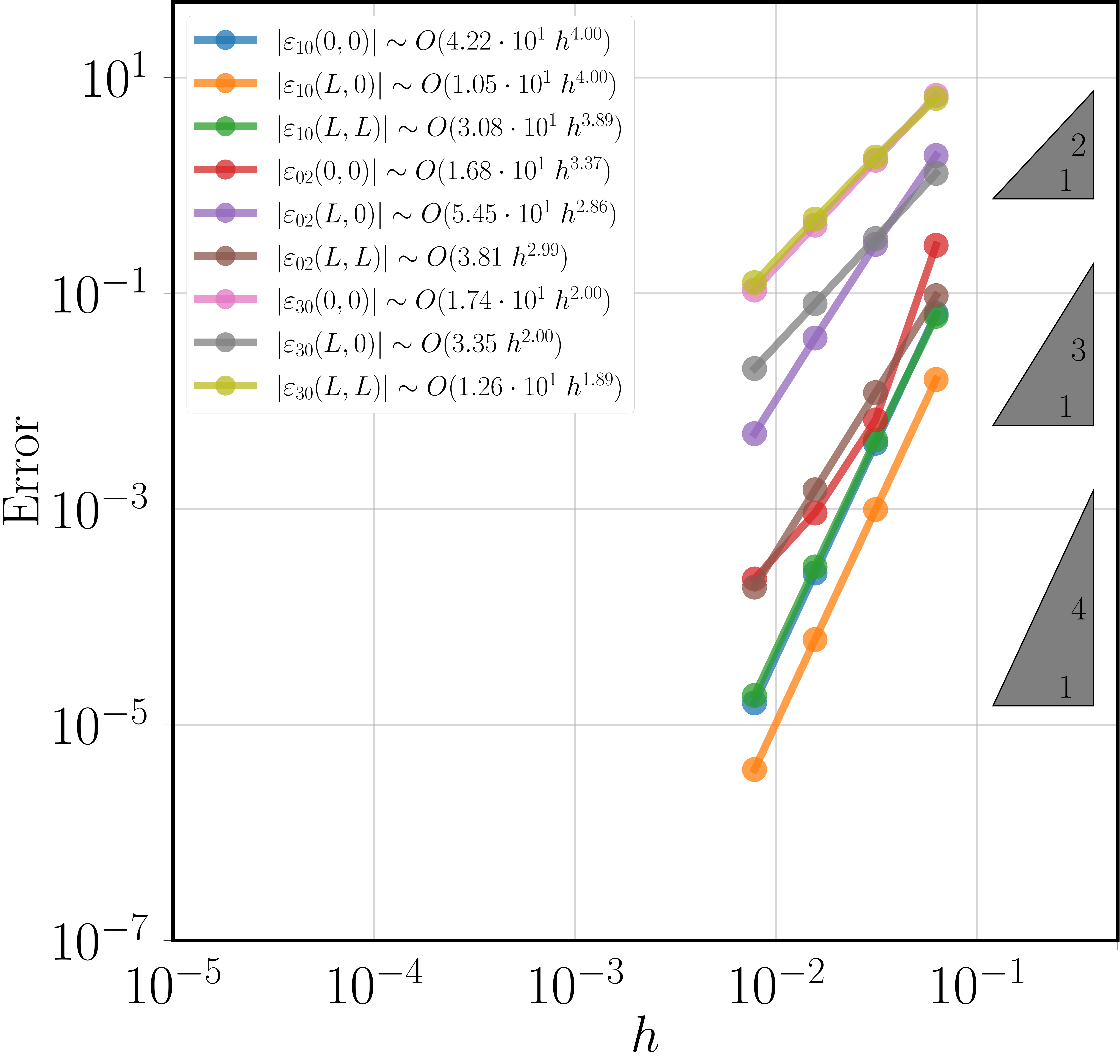}
	\subcaption{Local derivative error $\varepsilon_{lm}(x)$, $l+m = \{1 \cdots k\}$ for derivatives of order $l,m$ with respect to $x = \{x^{0}, x^{1}\}$.}
	\label{fig:errorscaling_2_derivativeerror_local}
\end{subfigure}
\caption{Error scaling as a function of data length scale $h$, of a $r=k+1$ accurate, $k=3$-order Taylor series model for a $K=6$ order polynomial $u(x^{0}, x^{1}) = \sum_{l+m \leq K}\alpha_{lm} {x^{0}}^{l}{x^{1}}^{m}$ in $p=2$ dimensions. Polynomial coefficients $\alpha_{lm} \sim U[-1,1]$ are sampled from a uniform distribution and all unique commuting Taylor series derivatives are included in the model up to $k$ order terms. Inset triangles indicate theoretical slopes of scaling of fits with length scale $h$, where model error $e \sim \mathcal{O}(h^{k+1})$ and $\ith[l,m]$ order derivative error $\varepsilon_{lm} \sim \mathcal{O}(h^{r+1-l-m})$.}
\label{fig:errorscaling_2}
\end{figure}

We now present in in \cref{fig:errorscaling_1,fig:errorscaling_2} numerical results for $p=1,2$ dimensional error scaling of the local and global model error for a $k$-order modified Taylor series, as per \cref{eq:error_local_model,eq:error_global_model}, as well as the local and global error of non-local derivatives, as per \cref{eq:error_local_derivative,eq:error_global_derivative}. As predicted by the analysis in \cref{app:error_analysis_model}, the plots confirm that this $k$-order Taylor series model, with $n = \mathcal{O}(1/h)$ local models and $r = k+1$-order accurate local stencils, has non-local $\ith[l]$ order derivative error of $\varepsilon_l = \mathcal{O}(h^{r+1-l})$, and local and global model error of $e = \mathcal{O}(h^{k+1})$.

\section{Physical systems of interest} \label{sec:physicalsystems}

We apply the non-local calculus on graphs to obtain reduced-order models for two physical systems of first-order, therefore dissipative, dynamics. In addition to computing derivatives of states as laid out at length in \cref{sec:graphtheory}, we use regression methods to select combinations of non-local derivative and algebraic operators on the states.

\subsection{Allen-Cahn dynamics}
\label{sec:physicalsystems_diffusion}

\begin{figure}[hpt]
\centering
\begin{subfigure}[t]{0.32\textwidth}
	\centering
	\includegraphics[width=\textwidth]{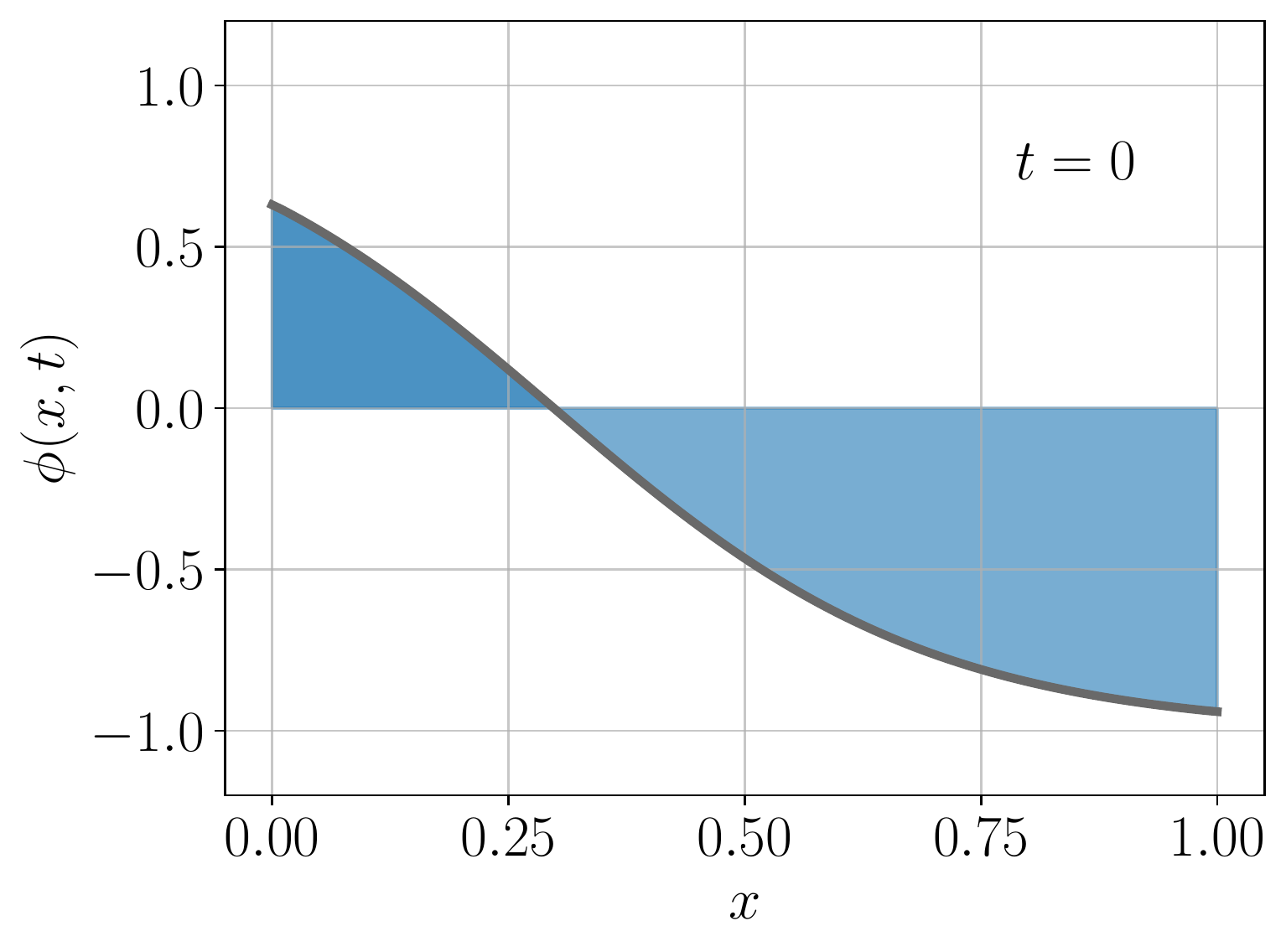}
	\subcaption{Initial condition.}
	\label{fig:AllenCahn1D_init}
\end{subfigure}
\hfill
\begin{subfigure}[t]{0.32\textwidth}
	\centering
	\includegraphics[width=\textwidth]{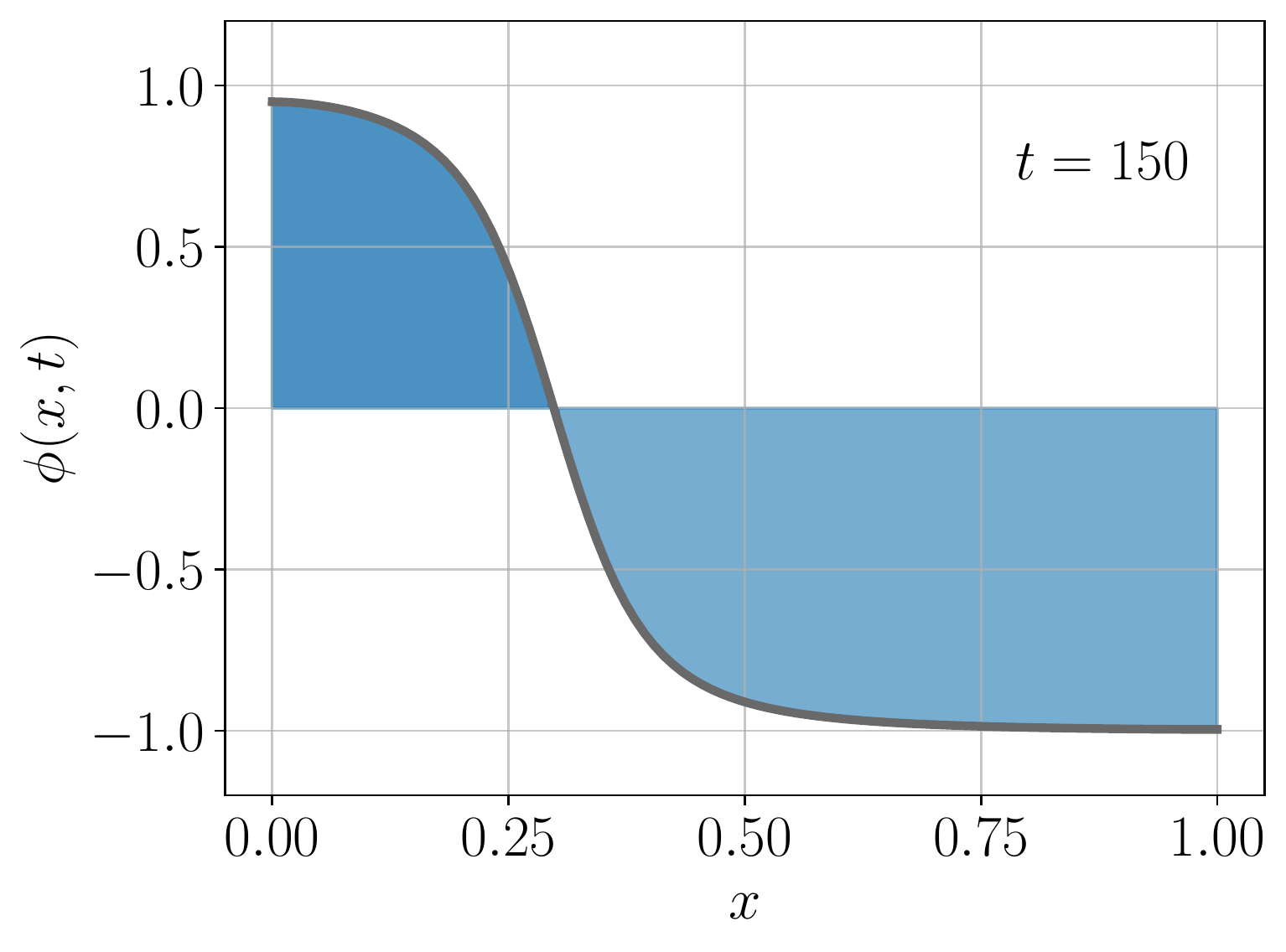}
	\subcaption{Solution at intermediate $t$.}
	\label{fig:AllenCahn1D_neq}
\end{subfigure}
\begin{subfigure}[t]{0.32\textwidth}
	\centering
	\includegraphics[width=\textwidth]{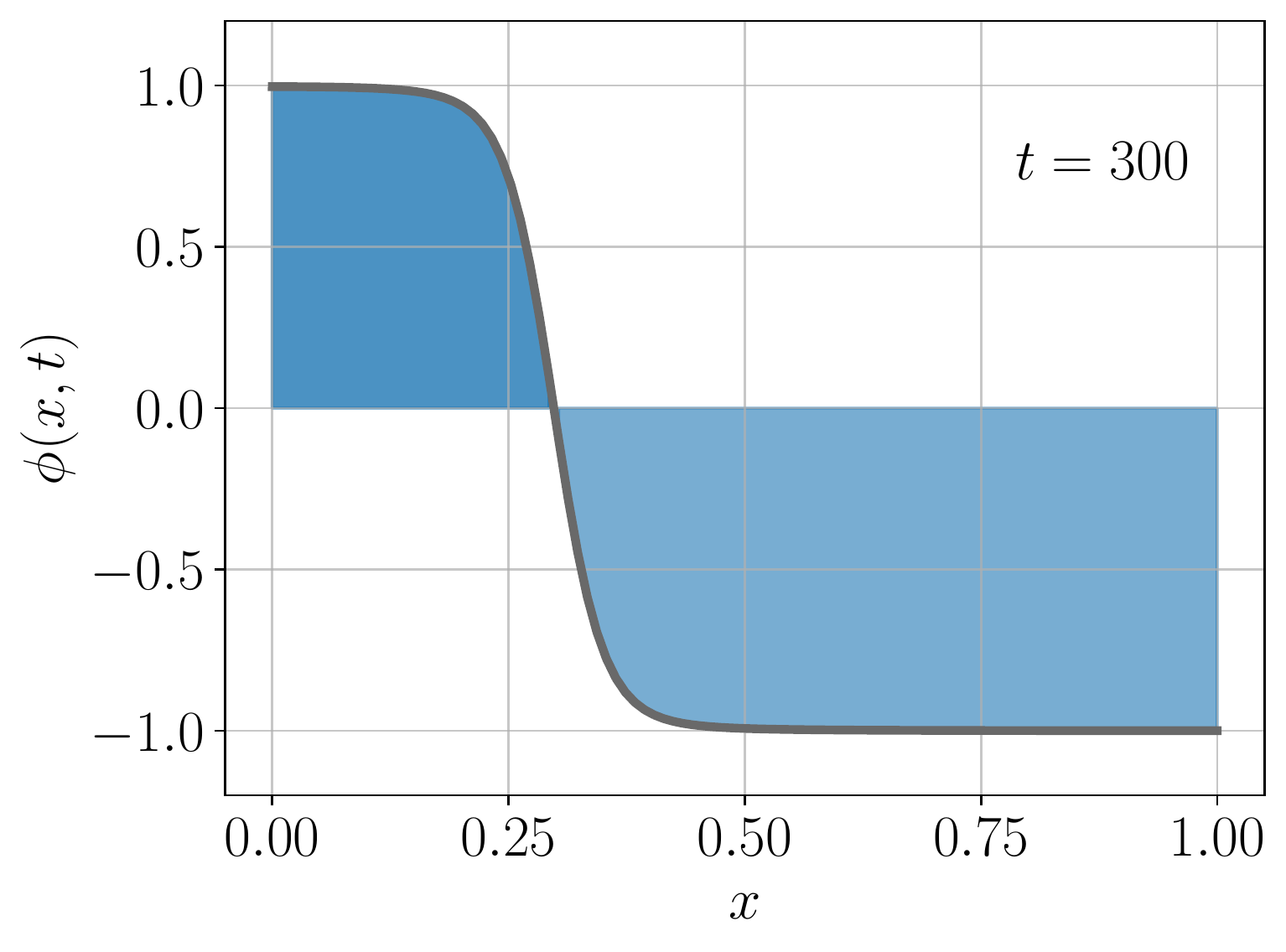}
	\subcaption{Near equilibrium solution.}
	\label{fig:AllenCahn1D_eq}
\end{subfigure}
\caption{Field evolution of 1D Allen-Cahn dynamics with $M_{\phi} = 10^{-3}$ and $\lambda = 1$ at $0$, $150$, and $300$ time steps. A Backward-Euler scheme is used with a time step of $\Delta t = 10^{-2}$.}
\label{fig:AllenCahn1D}
\end{figure}

Consider a field $\phi = \phi(x,t): \Omega \times [0,T] \mapsto \mathbb{R}$, governed by first order dynamics driven by gradient flow:
\begin{align}
	\derivative[1]{\phi}{t} =&~ -M_{\phi}\unidifference[1]{\psi}{\phi}, \quad \text{in}\; \Omega \times [0,T], \label{eq:CH}\\
	\nabla\phi\cdot\boldsymbol{n} &=~ 0\quad \text{on}\; \partial\Omega
	\label{eq:CHbcs}\\
	\phi(x,0) &= \phi_0(x) \label{eq:CHics}
\end{align}
Here, the free energy density 
\begin{align}
	\psi =&~ f(\phi) + \frac{1}{2}\lambda\abs{\gradient[][]{\phi}}^2
\end{align}
includes $f$, an algebraic Landau energy density of the form
\begin{align}
	f(\phi) =&~ (1-\phi^2)^2,
	\label{eq:landau}
\end{align}
with wells at $\phi = \pm 1$. The gradient energy $\lambda\abs{\gradient[][]{\phi}}^2$, with $\lambda > 0$ penalizes sharp transitions between the positive and negative phases $\Omega_\pm \subset \Omega$, which are defined by
\begin{equation}
x \in 
\begin{cases}
    \Omega_+ & \text{if }\phi(x) \ge 0\\
    \Omega_- & \text{if }\phi(x) < 0
\end{cases}
\end{equation}
The kinetics are controlled by the local mobility, $M_{\phi} \ge 0$. \Cref{eq:CH,eq:CHbcs} constitute the Allen-Cahn equation..\cite{Allen1979} An example of the system dynamics in 1D appears in \cref{fig:AllenCahn1D}

\noindent Of interest here are states of the system, which are functionals of $\phi$, such as the volume averaged total energy
\begin{align}
	\Psi[\phi] =&~ \int\limits_{\Omega} \psi(x,t)~\mathrm{d}V = \int\limits_{\Omega} \left( f(\phi) + \frac{1}{2}\lambda\abs{\gradient[][]{\phi}}^2\right)~\mathrm{d}V
	\label{eq:totalenergy}
\end{align}
as well as others obtained, using the indicator function
\begin{equation}
    I(\phi) = \begin{cases}
                1 & \text{if }\phi \ge 0\\
                0 & \text{if }\phi <0,
              \end{cases}    
              \label{eq:indicfn}
\end{equation}
on the Landau energy density and its derivative
\begin{align}
	F_\pm[\phi] &= 
	\begin{cases}
	& \int\limits_{\Omega} ~ I(\phi)~ f~ \mathrm{d}V,\\
	& \int\limits_{\Omega} ~ \left(1-I(\phi)\right)~ f~ \mathrm{d}V,
	\end{cases}
	\\
	F^{\prime}_\pm[\phi] &= 
	\begin{cases}
	& \int\limits_{\Omega} ~ I(\phi)~ f^{\prime}~ \mathrm{d}V,\\
	& \int\limits_{\Omega} ~ \left(1-I(\phi)\right)~ f^{\prime}~ \mathrm{d}V,
	\end{cases}
\end{align}
powers of the phase field and its gradients
\begin{align}
	\varphi_{k_\pm}[\phi] =&~ 
	\begin{cases}
	& \int\limits_{\Omega} ~ I(\phi)~ \phi^k~ \mathrm{d}V\\
	& \int\limits_{\Omega} ~ \left(1-I(\phi)\right)~ \phi^k~ \mathrm{d}V
	\end{cases}, \\
	\varphi_{\scriptscriptstyle{\gradient[k_{\pm}][][]}}[\phi] =&~
	\begin{cases}
	& \int\limits_{\Omega}~ I(\phi) ~\gradient[k][]{\phi}~ \mathrm{d}V\\
	& \int\limits_{\Omega}~ \left(1-I(\phi)\right) ~\gradient[k][]{\phi}~ \mathrm{d}V
	\end{cases}, \\
	\varphi_{\scriptscriptstyle{\gradient[][k_{\pm}][]}}[\phi] =&~ 
	\begin{cases}
	& \int\limits_{\Omega}~ I(\phi) ~\vert\gradient[][]{\phi}\vert^k~ \mathrm{d}V\\
	& \int\limits_{\Omega}~ \left(1-I(\phi)\right) ~\vert\gradient[][]{\phi}\vert^k~ \mathrm{d}V
	\end{cases}{},\\
	\varphi_{\scriptscriptstyle{k},\scriptscriptstyle{\nabla^l}_{\pm}} =&~ 
	\begin{cases}
	& \int\limits_\Omega I(\phi)~ \phi^k~ \nabla^l\phi~ \mathrm{d}V\\
	& \int\limits_\Omega \left(1-I(\phi)\right)~ \phi^k~ \nabla^l\phi~ \mathrm{d}V
	\end{cases},
	\label{eq:moments}
\end{align}

\noindent The local chemical potential is obtained by computing the variational derivative
\begin{align}
\frac{\mathrm{d}}{\mathrm{d}\varepsilon}\Psi[\phi + \varepsilon\chi]\Big\vert_{\varepsilon = 0} =&~ \int\limits_{\Omega} \difference[1]{\psi}{\phi}~ \frac{\mathrm{d}}{\mathrm{d}\varepsilon}(\phi+\varepsilon\chi)\Big\vert_{\varepsilon = 0}~\mathrm{d}V \nonumber \\
	=&~ \int\limits_\Omega \chi \left(f^{\prime} - \lambda \gradient[2][]{\phi}\right) ~\mathrm{d}V + \lambda\int\limits_{\partial\Omega}\chi \nabla\phi\cdot\boldsymbol{n}\mathrm{d}S.
\end{align}
Applying the boundary conditions in \cref{eq:CHbcs} the local chemical potential is 
\begin{equation}
    \mu_\phi = \difference[1]{\psi}{\phi} := f^\prime - \lambda \gradient[2][]{\phi}.
    \label{eq:localchempot}
\end{equation}
We note that at equilibrium the Euler-Lagrange relation is $\mu_\phi = f^\prime - \lambda \gradient[2][]{\phi} = 0$. \\
\noindent These dynamics dissipate the total free energy since:
\begin{align}
	\frac{\mathrm{d}\Psi}{\mathrm{d}t} = \int\limits_{\Omega} \difference[1]{\psi}{\phi}&~ \derivative[1]{\phi}{t}~\mathrm{d}V = -\int\limits_{\Omega} M_{\phi}\mu_\phi^2~\mathrm{d}V\\ 
	\intertext{and therefore}
	 \frac{\mathrm{d}\Psi}{\mathrm{d}t} &\le 0.
	 \label{eq:dissip}
\end{align}

\subsubsection{A reduced-order model of gradient flow}

\noindent We now seek a reduced-order model in terms of the states defined in \crefrange{eq:totalenergy}{eq:moments}. We consider the volume fraction of the positive phase, $\varphi_{} \equiv \varphi_{1^+} = \int\limits_\Omega I(\phi) \phi ~\mathrm{d}V$. A reduced-order model can be derived:
\begin{align}
	\derivative[1]{\varphi_{}}{t} =&~ \int\limits_{\Omega} (I(\phi) + {I}^{\prime}(\phi) \phi)\derivative[1]{\phi}{t}~ \mathrm{d}V.
\end{align}
Noting from \cref{eq:indicfn} that $I^\prime(\phi) = \pm\delta[\phi;0]$, it follows that the second term in parentheses vanishes. Then, substituting \cref{eq:CH,eq:landau,eq:localchempot} leads to
\begin{alignat}{2}
    \frac{\partial \varphi_{}}{\partial t} &= - &&\int\limits_\Omega I(\phi)~ M_\phi \frac{\delta\Psi}{\delta \phi}~ \mathrm{d}V\label{eq:rom1}\\ 
    &= &&\int\limits_\Omega I(\phi)~ M_\phi\left((4\phi^2-4\phi^4) + \lambda\phi\nabla^2\phi \right) ~\mathrm{d}V\label{eq:rom2}
\end{alignat}
We note that \cref{eq:rom2} defines an exact evolution equation for $\varphi_{}$. However, given the gradient flow form of the integrand in \cref{eq:rom1} it is of interest to explore its preservation in a global gradient flow model analogous to \cref{eq:CH}, and of the form $\partial\varphi_{}/\partial t \sim \delta\Psi/\delta\varphi_{}$, where the gradient $\delta\Psi/\delta\varphi_{}$ can be thought of as a global chemical potential, defined following the non-local calculus in \crefrange{sec:graphtheory_nonlocalcalculus}{sec:localweights}, and discrepancy terms are added:
\begin{align}
	\derivative[1]{\varphi_{}}{t} =&~ -M_{\varphi} \difference[1]{\Psi}{\varphi_{}} + \mathcal{E}_{\varphi},\label{eq:romdiffeq}\\
	\varphi_{}(0) =& \int\limits_\Omega I(\phi(x,0))\phi(x,0) ~\mathrm{d}V\label{eq:romic}
\end{align}
 A kinetic parameter $M_{\varphi}$ has been introduced, and $\mathcal{E}_{\varphi}$ is the discrepancy in the assumed global gradient flow model.

\noindent Guided by \cref{eq:rom2} we propose model forms
\begin{align}
	M_{\varphi} =&~ M_{\varphi}(\varphi_{\scriptscriptstyle{k}_{\pm}},\varphi_{\scriptscriptstyle{\gradient[k_{\pm}][][]}},\varphi_{\scriptscriptstyle{\gradient[][k_{\pm}][]}},\dots),
	\intertext{and}
	\mathcal{E}_{\varphi} =&~ \mathcal{E}_{\varphi}(\varphi_{\scriptscriptstyle{k}_{\pm}},\varphi_{\scriptscriptstyle{\gradient[k_{\pm}][][]}},\varphi_{\scriptscriptstyle{\gradient[][k_{\pm}][]}},\dots).	
\end{align}
We also propose a functional form for the total free energy of
\begin{align}
	\Psi =&~ \Psi(\varphi_{\scriptscriptstyle{k}_{\pm}},\varphi_{\scriptscriptstyle{\gradient[k_{\pm}][][]}},\varphi_{\scriptscriptstyle{\gradient[][k_{\pm}][]}},\dots).	
\end{align}

\begin{figure}[ht]
\centering
\includegraphics[width=0.5\textwidth]{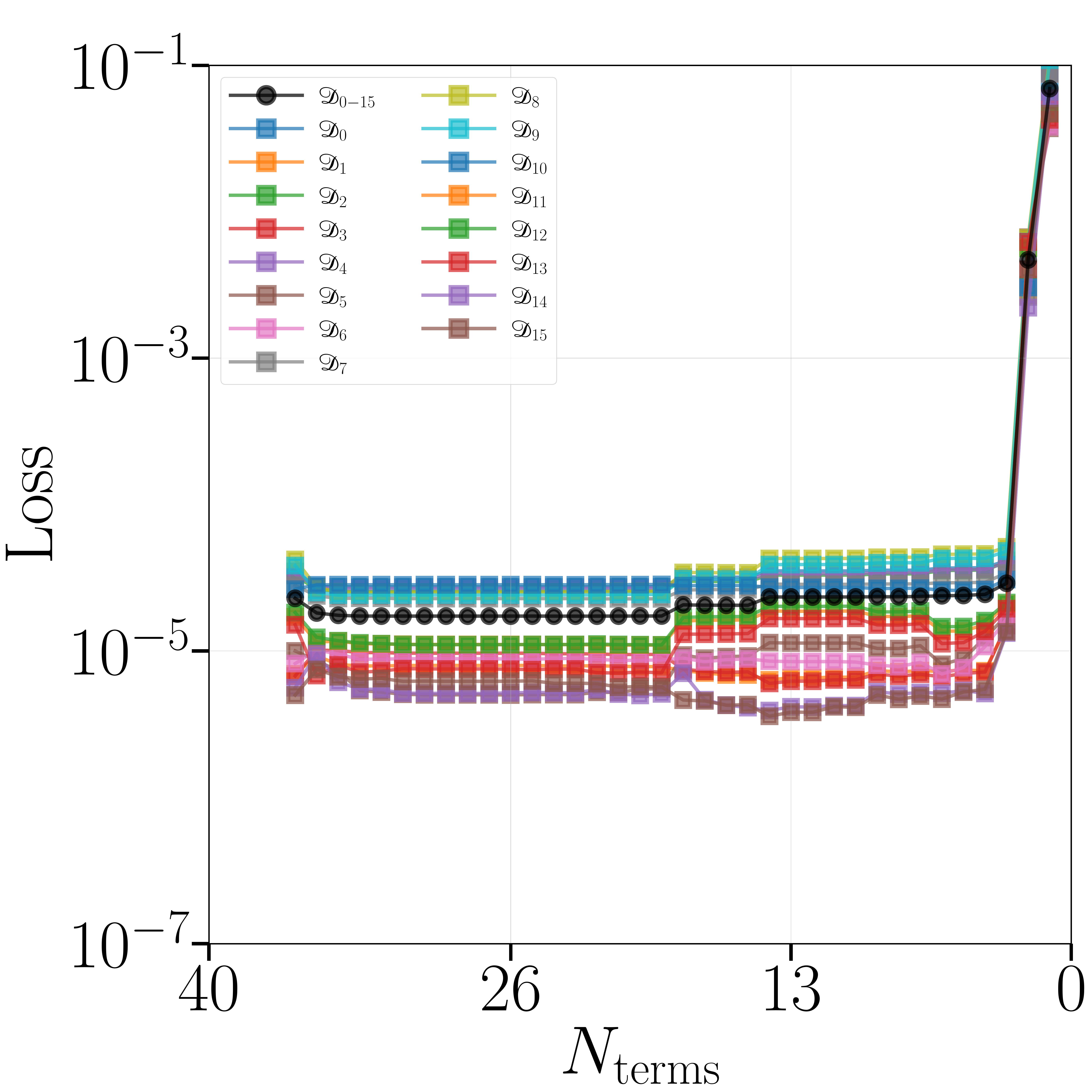}
\captionsetup{width=1\linewidth}		
\caption{Loss curves for the reduced-order model of Allen-Cahn dynamics with increasing parsimony.}
\label{fig:diffusion_loss}	
\end{figure}

\begin{figure}[ht]
	\centering
	\includegraphics[width=0.5\textwidth]{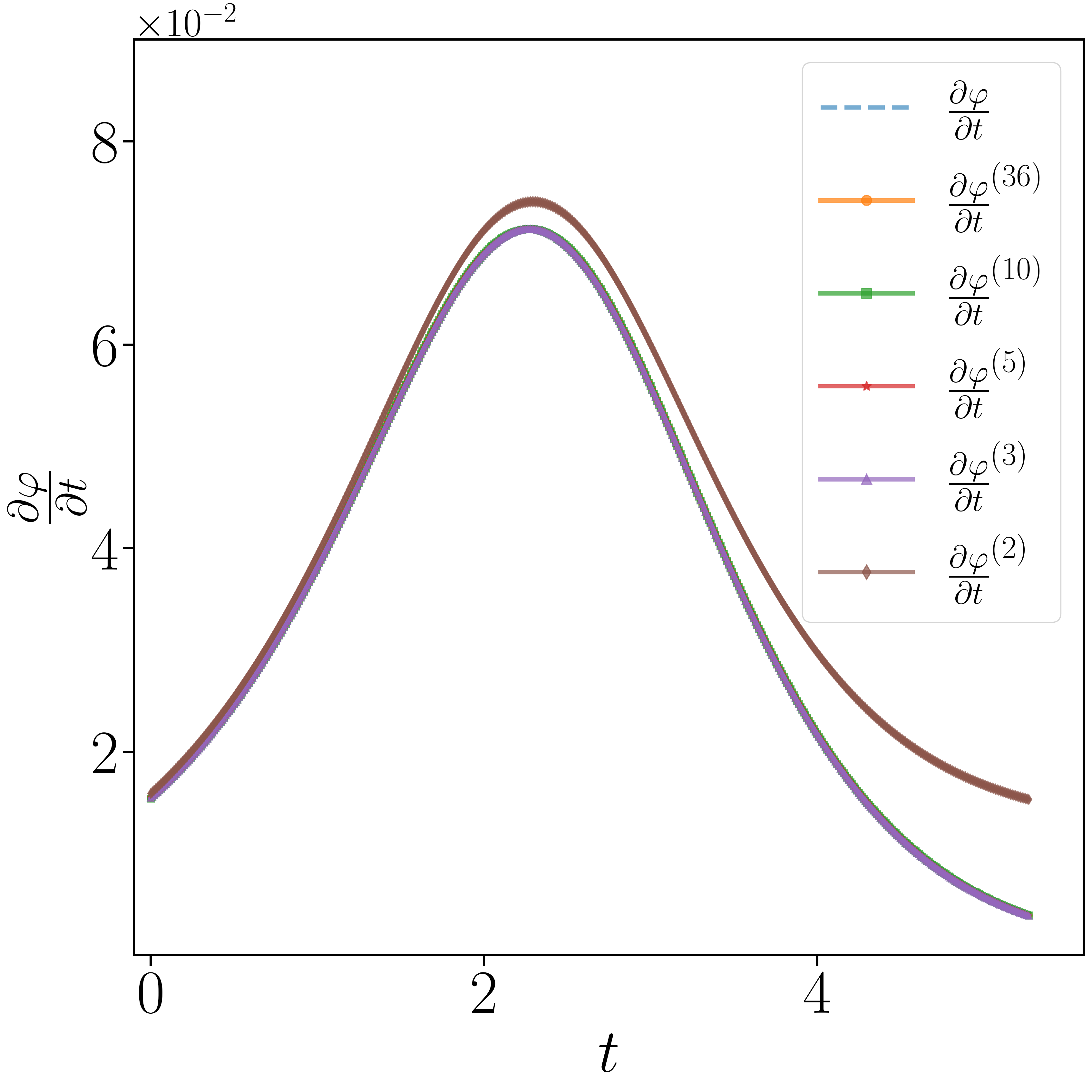}
	\captionsetup{width=1\linewidth}		
	\caption{Fitted global dynamics for $\mathcal{D}_{5}$ as the model becomes more parsimonious. The dashed blue curve is the first order derivative data. The number of operators in each reduced-order model appears in superscript parentheses. The data are matched by all models except for the two-term model.}
	\label{fig:diffusion_bestfit}
\end{figure}

\noindent High fidelity data is obtained for different choices of material parameters $M_{\phi},\lambda$, as well as initial conditions. Each simulation yields states defined as in \crefrange{eq:totalenergy}{eq:moments} at each time instant. They are vertices of a directed graph with edges whose sense is dictated by time. We call attention to the directedness of the graph arising from the first-order nature of these dynamics that impose time irreversibility due to dissipation as we have observed previously.\cite{Banerjee2019} As an example, the data presented in \cref{fig:AllenCahn1D} produces a non-branching tree with $387$ vertices corresponding to an equal number of time instants making up a trajectory. Please refer to \cref{fig:treegraph} for examples of directed branching trees. The full dataset consists of $N = 16$ such trajectories, represented as $\mathcal{D}_j = \{\varphi_{k_{\pm}},\varphi_{\scriptscriptstyle{\gradient[2_{\pm}][][]}}, \Psi,\unidifference[1]{\Psi}{\varphi},\dots\}^{(j)}$. Given a basis of operators, for example,
\begin{align}
	v = \{\varphi_{k_{+}},\varphi_{\scriptscriptstyle{\gradient[2][][]}}&~,F_{+},F^{\prime}_{+},\dots\}, \\
\intertext{and expanding the kinetic parameter and discrepancy:}
	M_{\varphi} =~ \sum_{\alpha} \gamma^{M \alpha} v_{\alpha}, \quad&\quad \mathcal{E}_{\varphi} =~ \sum_{\beta} \gamma^{\beta} v_{\beta} 
\end{align}
we can fit a model with the $\gamma$ coefficients. We choose this global model basis to include $q=36$ terms. The second and fourth powers, and product of linear and Laplacian terms are suggested by the exact model in \cref{eq:rom2}. Others account for global effects from the boundaries and non-local correlations that can  be exactly resolved only by spatially integrating the PDE in \crefrange{eq:CH}{eq:CHics}. For a global, reduced-order model to match the data from $N$ trajectories, terms are needed that reflect finite domain effects. In general, this translates to requiring more terms than in the local PDE model. Ridge regression and cross validation across the trajectories prevents overfitting, and yields the results in \cref{fig:diffusion_loss,fig:diffusion_bestfit}.

\noindent We can see from the loss curves in \cref{fig:diffusion_loss} that there is a sharp Pareto front at around $3$ terms in the model, where the maximum model has $36$ terms. Recalling that the exact reduced-order model, \cref{eq:rom2} has three terms, we note that the inferred ones also have low losses with three operators. We also see that the fitted curves in \cref{fig:diffusion_bestfit} match quite well visually, and only with less than $3$ terms does the reduced model differ significantly from the exact dynamics. The resulting models for $10$ to $1$ operators are

\input{./figures/model_diffusion_local.tex}

The Pareto front with a sharp increase in loss occurs exactly when the mobility term goes to zero and the discrepancy term dominates the model. This suggests the global diffusion model \crefrange{eq:romdiffeq}{eq:romic} is a sound ansatz, and is only being partially masked by overfitting from the persisting discrepancy terms. Future work should investigate other ansatz, and how well this mixture of physics informed terms, with discrepancy terms affects the effectiveness of the reduced-order models.


\subsection{Microstructures in a gradient-regularized model of non-convex elasticity} \label{sec:physicalsystems_microstructures}
We apply the graph theoretic approach to studying reduced-order models for the mechano-chemical response of solids that undergo phase transformations mediated by composition (chemical) and strain (elastic) variables. This physical system is driven by a free energy density function that, naturally, is parameterized by composition  and strain. The underlying functional form of the free energy density is non-convex in both these quantities. Microstructures develop as phases and symmetry-breaking structural variants arise corresponding to negative eigenvalues of the Hessian of the free energy density in the strain-composition space. These models are ill-posed, however; a condition arising from the absence of penalization on inter-phase and inter-variant interfaces and manifesting in PDE solutions that are pathologically discretization-dependent. Well-posedness and penalization of interfaces are restored by free energy density functionals $\psi = \psi_{\textrm{hom}} + \psi_{\textrm{grad}}$, dependent on composition and strain gradients. When decoupled, the chemical component of this problem causes spinodal decomposition and Ostwald ripening, described by the Cahn-Hilliard equation. \cite{Cahn1958} The non-linear elasticity problem is described by strain gradient elasticity, of which many forms have been proposed.

Our treatment follows Toupin. \cite{Toupin1962} and has appeared as decoupled, gradient-regularization of non-convex, non-linear elasticity \cite{Rudraraju2014,Teichertetal2017}, as well as mechano-chemical spinodal decomposition \cite{Rudraraju2016,Sagiyamaetal2016,SagiyamaGarikipati2017b} The free energy density functional $\psi$ depends on the fields $\{c,\vec{E},\gradient[][]{c},\gradient[][]{\vec{E}}\}$ where $c$ is the composition and $\vec{E}$ is the Green-Lagrange strain tensor.

\begin{figure}[ht]
	\centering
	\includegraphics[width=0.6\textwidth]{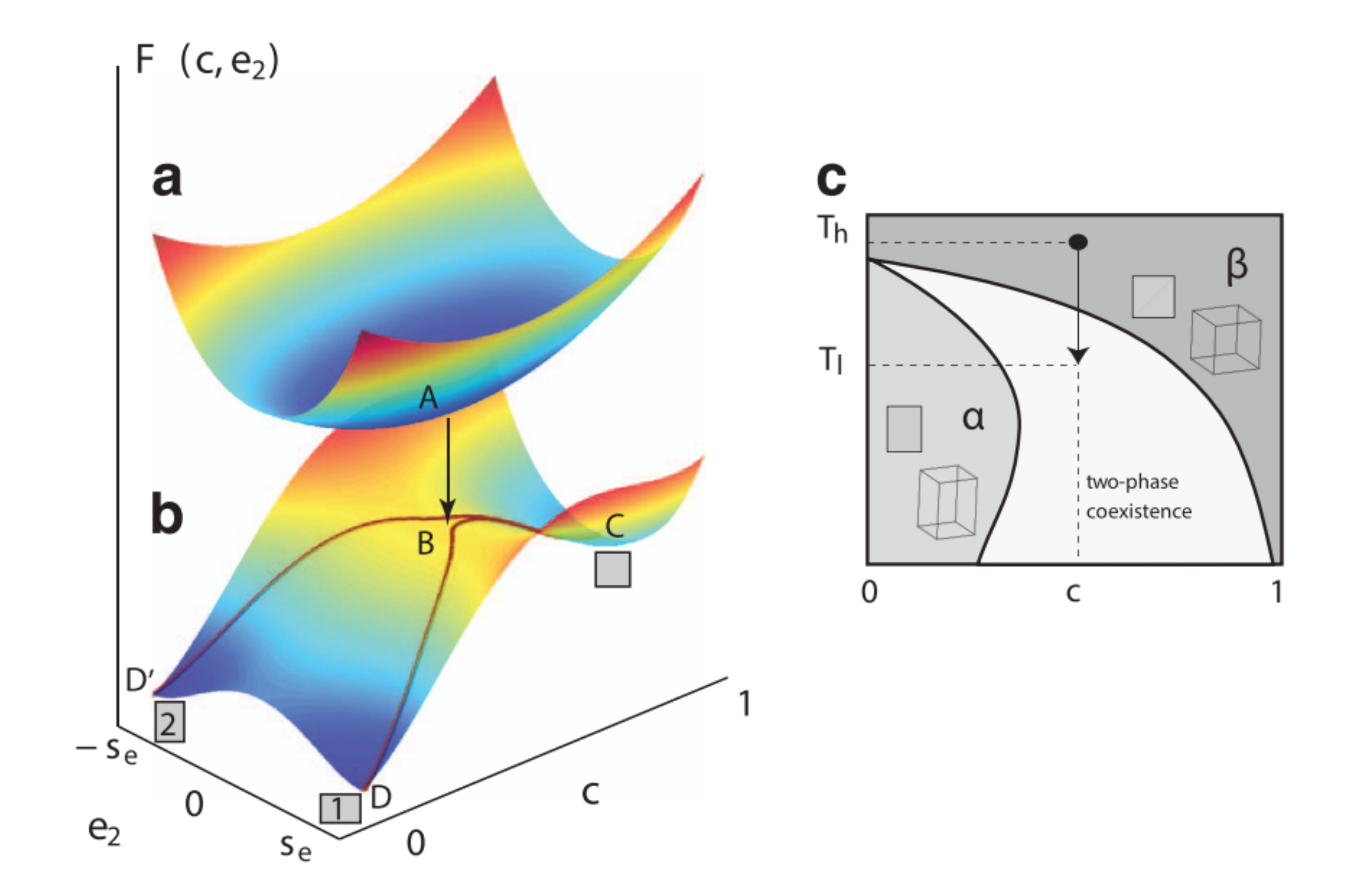}
	\caption{Free energy density landscape over the composition $c$ and strain component $e_{2}$ order parameters, showing changes in curvature and resulting phases of square and rectangular lattice symmetries (1 and 2), \textbf{a} being at high temperature and \textbf{b} at a quenched, low temperature state, from Rudraraju \etal. \cite{Rudraraju2016} On the right, the effect of temperature and composition on the phase diagram of stable regions of purely square or rectangular phases, or a mixture of phases is shown.}
	\label{fig:energyphase}
\end{figure}

To describe these phase transitions, we follow Zhang \& Garikipati. \cite{Zhang2020} At low temperatures, these systems consist of stable phases related to the local crystalline symmetry, specifically cubic, or tetragonal symmetry. The system considered is restricted to two spatial dimensions, and so the material will consist of phases with either square or rectangular structure, or symmetry, as described by linear combinations of the strain components $\vec{e}(\vec{E})$. Chemically the system is modeled as a binary mixture, with scalar composition $0\leq c \leq 1$ of one of the two constituents. With the composition $c$ and symmetry-coding strain parameter $\vec{e}$ and their gradients, the non-uniform free energy density is written: $\psi_{\textrm{grad}}(\gradient[][]{c}, \gradient[][]{\vec{e}})$. The non-convex  dependence of $\psi_{\textrm{hom}}$ on $c$ and $\vec{e}$ is shown in the reduced two-dimensional case by \cref{fig:energyphase}, and is regularized by, typically, quadratic dependence on the gradients $\gradient[][]{c}$ and $\gradient[][]{\vec{e}}$ in $\psi_{\textrm{grad}}$. The gradient contributions model interface energies. Together, $\psi_{\textrm{hom}}$ and $\psi_{\textrm{grad}}$ control a complex mechano-chemical evolution of microstructure \cite{Rudraraju2016} 

In the two-dimensional case considered here, $(c,e_1,e_{2},e_6)$ are order parameters. Of these, $c$ and $e_2$ define the microstructure, leading to phase volume fractions $\varphi_{\alpha}$, where $\alpha = \{\square,\hrectangle,\vrectangle\}$ correspond, respectively, to the square and two rectangular variants illustrated as 1 and 2 in \cref{fig:energyphase}. Furthermore, the gradient regularization that represents the energy of interfaces also manifests in a number $N_{\alpha}$ of each variant, and total interfacial lengths $l_{\alpha}$.

A variational formalism can then be applied to yield the relevant generalized chemical potentials. For this particular work, the homogeneous free energy density is a smooth polynomial function of the order parameters,
\begin{equation} 
	\psi_{\textrm{hom}}(c,\vec{e}) \equiv 16 \alpha_{c} c^4 - 32\alpha_{c} c^3 + \alpha_{c} c^2 + 2\frac{\alpha_{e}}{{\beta_{e}}^2}({e_{1}}^2 + {e_{6}}^2) + \frac{\alpha_{e}}{{\beta_{e}}^2}{e_{2}}^4 + 2\frac{\alpha_{e}}{{\beta_{e}}^2}(2c-1){e_{2}}^2,
	\label{eq:psi_hom}
\end{equation}
where $\{\alpha_{c},\alpha_{e},\beta_{e}\}$ are constant coefficients.

The gradient contributions will take the form of strictly quadratic forms in the order parameter gradients,
\begin{equation}
	\psi_{\textrm{grad}}(\gradient[][]{c},\gradient[][]{\vec{e}}) = \frac{1}{2} \gradient[][]{c} \cdot\vec{\kappa} \gradient[][]{c} + \frac{1}{2} \gradient[][]{\vec{e}} \cdot \boldsymbol{\gamma} \gradient[][]{\vec{e}}, \label{eq:psi_grad}
\end{equation}
 and only constant isotropic tensors $\vec{\kappa} = \kappa \textrm{I}$ and $\boldsymbol{\gamma} = \gamma \vec{I}$ are considered here for the gradient terms. 

When stationarity of the free energy is imposed in a variational framework, the Euler-Lagrange equations yield the weak form of the balance of linear momentum, and further variational arguments lead to the strong form. Using this variational approach, based on this form of the free energy, and given its dependence on the strain and composition fields, the governing PDE's that emerge are the coupled Cahn-Hilliard time-dependent parabolic equation for the composition \cite{Cahn1958}, and Toupin's non-linear gradient elasticity equations at equilibrium \cite{Toupin1962,Rudraraju2016}

\subsubsection{Cahn-Hilliard dynamics} \label{sec:physicalsystems_microstructures_CahnHilliard}
Cahn-Hilliard dynamics are first order in time for the composition and take the form of a transport equation
\begin{align}
	\derivative[1]{c}{t} + \divergence[][{\vec{J}}] ~=&~ 0, \label{eq:pde_cahnhilliard} 
\end{align}
where
\begin{align}
	\vec{J} ~=&~ -\vec{L}\gradient[][][\mu], \label{eq:cflux}
\end{align}
and $\vec{L} = \vec{L}(c,\vec{e})$ is a mobility transport tensor, and is assumed to be isotropic $\vec{L} = \textrm{L}\vec{I}$. The variational treatment yields the chemical potential $\mu = \partial\psi_\textrm{hom} - \kappa\nabla^2 c$, whose gradient yields the flux \cref{eq:cflux} in a form that  satisfies the thermodynamic dissipation inequality \cite{degrootmazur1984}, and when substituted in \cref{eq:pde_cahnhilliard} leads to  mechano-chemically coupled equations involving fourth-order spatial derivatives in composition. Initial conditions, and Dirichlet/Neumann boundary conditions of appropriate order are discussed in Rudraju et al.\cite{Rudraraju2016}

\subsubsection{Toupin model of elasticity} \label{sec:physicalsystems_microstructures_Toupin}
Toupin's \cite{Toupin1962} strain gradient elasticity formulation arises as the Euler-Lagrange equations obtained from seeking extremization of the total free energy with respect to displacement (elastic equilibrium). It leads to equations for the first Piola-Kirchhoff $\vec{P}$ and higher order $\vec{B}$ stress tensors,
\begin{align}
	\vec{P} ~=&~~ \derivative[1]{\psi}{\vec{F}} \label{eq:piola1}\\
	\vec{B} ~=&~~ \derivative[1]{\psi}{\gradient[][]{\vec{F}}} \label{eq:piola3}
\end{align}
which are derivatives of the elastic free energy density and are conjugate to the deformation gradient $\vec{F}$. Here, the Green-Lagrange strains are $\vec{E} = {1}/{2}\left(\vec{F}^{T}\vec{F} - \vec{I}\right)$ and gradients $\gradient[][]$ are derivatives with respect to the reference configuration. The elasticity governing equations, shown here in the strong form, are
\begin{align}
	\divergence[]{\vec{P}} - \divergence[]{(\divergence[]{\vec{B}})} = 0, \label{eq:pde_elasticity}
\end{align}
plus additional Dirichlet and higher order traction boundary conditions.

\subsubsection{Direct numerical simulations} \label{sec:physicalsystems_microstructures_DNS}

\begin{figure}[ht]
  \centering
  \subfloat[$c$ at step 51]{\includegraphics[trim={15cm 3.5cm 15cm 3.5cm}, clip, width=0.24\linewidth]{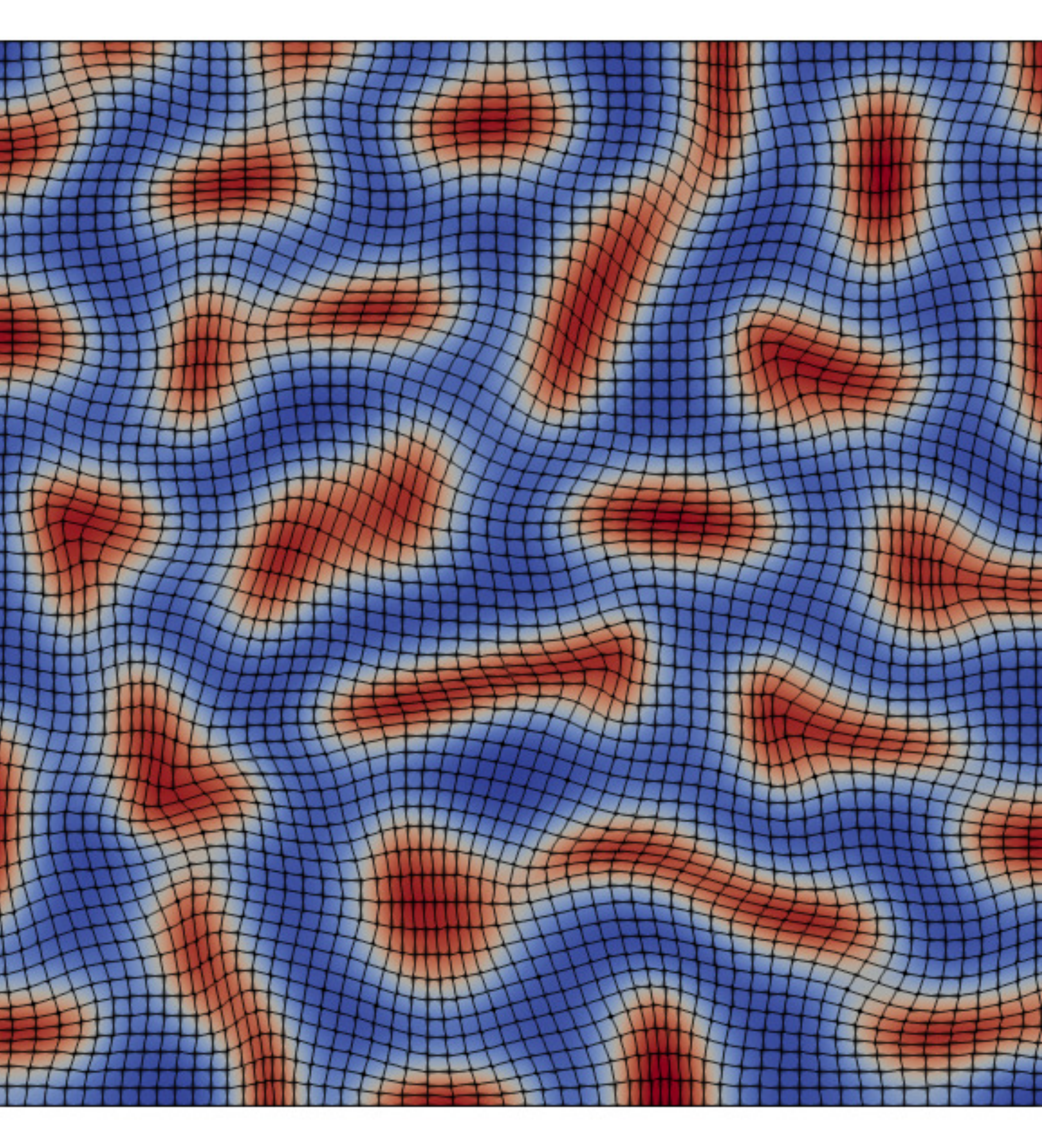}} \hfill
  \subfloat[$c$ at step 150]{\includegraphics[trim={15cm 3.5cm 15cm 3.5cm}, clip, width=0.24\linewidth]{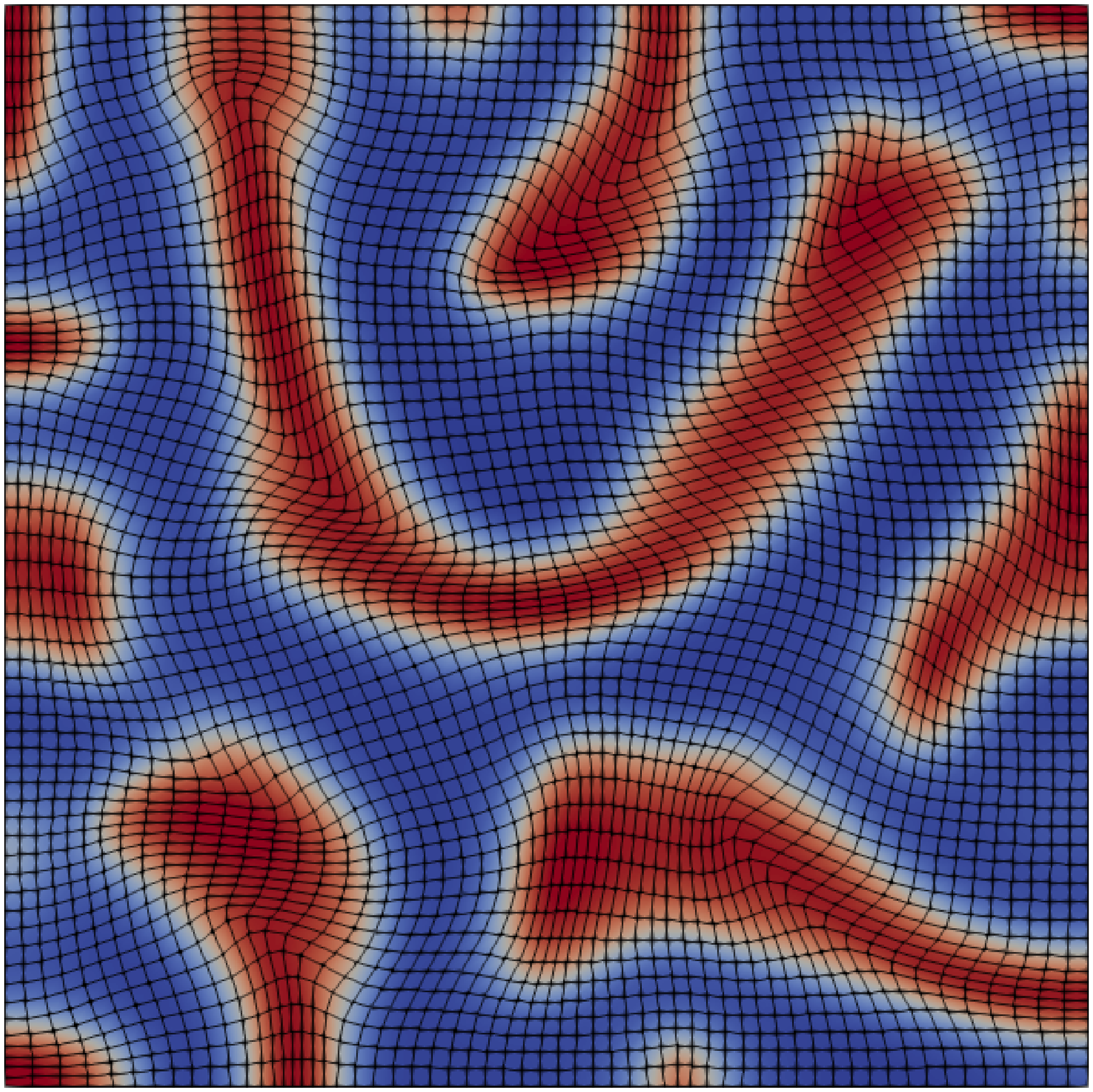}} \hfill
  \subfloat[$c$ at step 400]{\includegraphics[trim={15cm 3.5cm 15cm 3.5cm}, clip, width=0.24\linewidth]{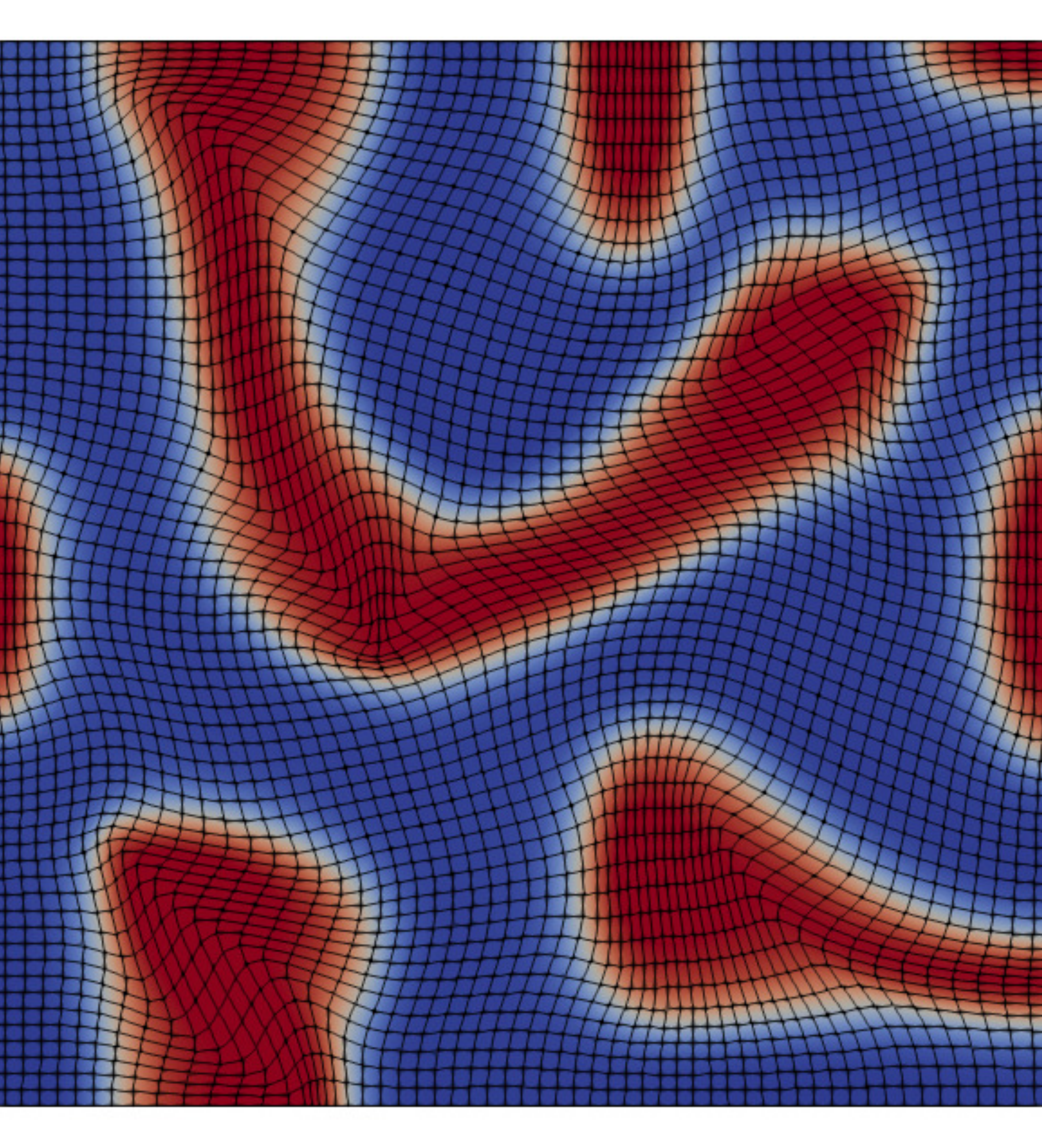}} \hfill
  \subfloat[$c$ at step 900]{\includegraphics[trim={15cm 3.5cm 15cm 3.5cm}, clip, width=0.24\linewidth]{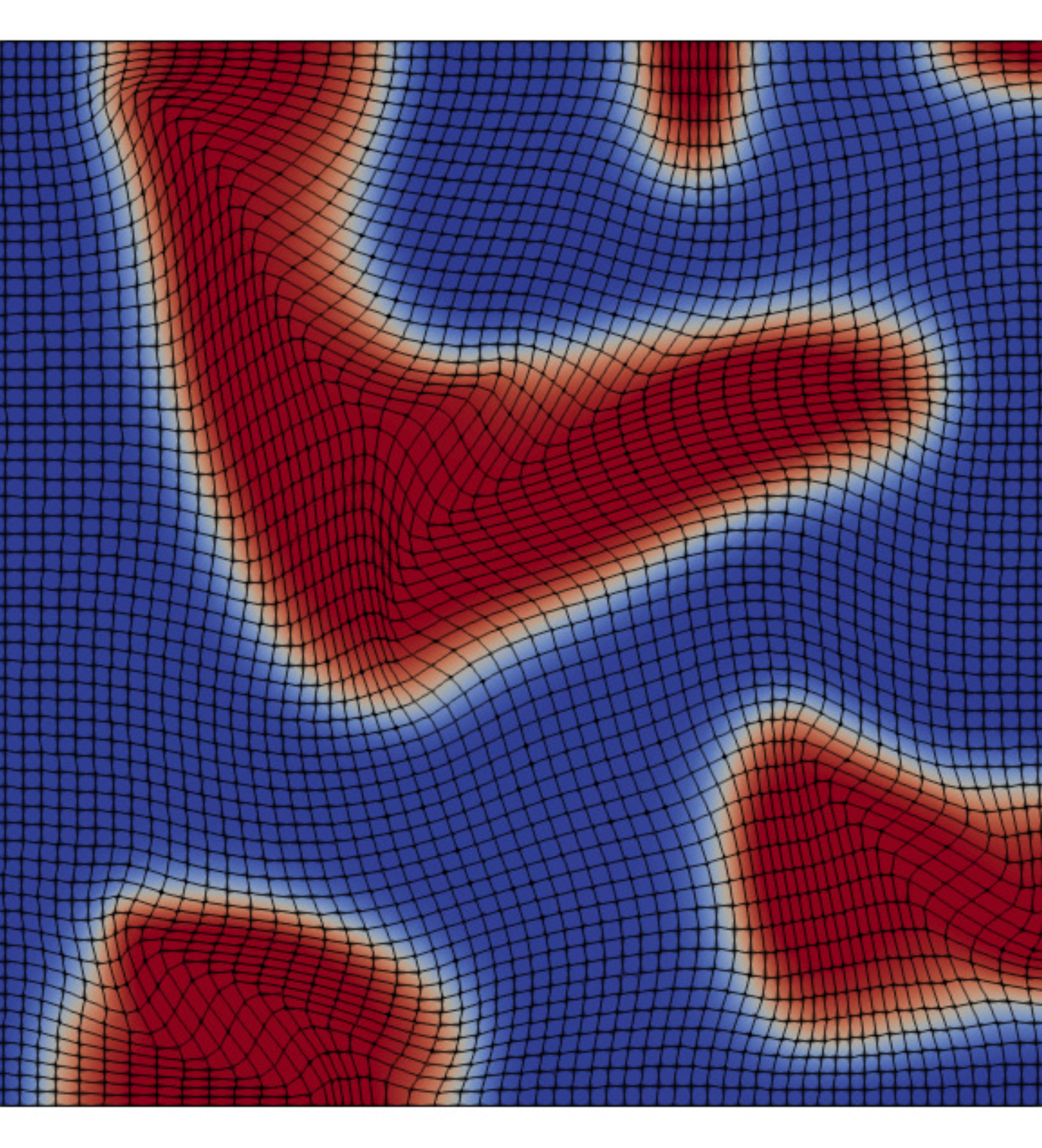}} \\
  \subfloat[$e_2$ at step 51]{\includegraphics[trim={15cm 3.5cm 15cm 3.5cm}, clip, width=0.24\linewidth]{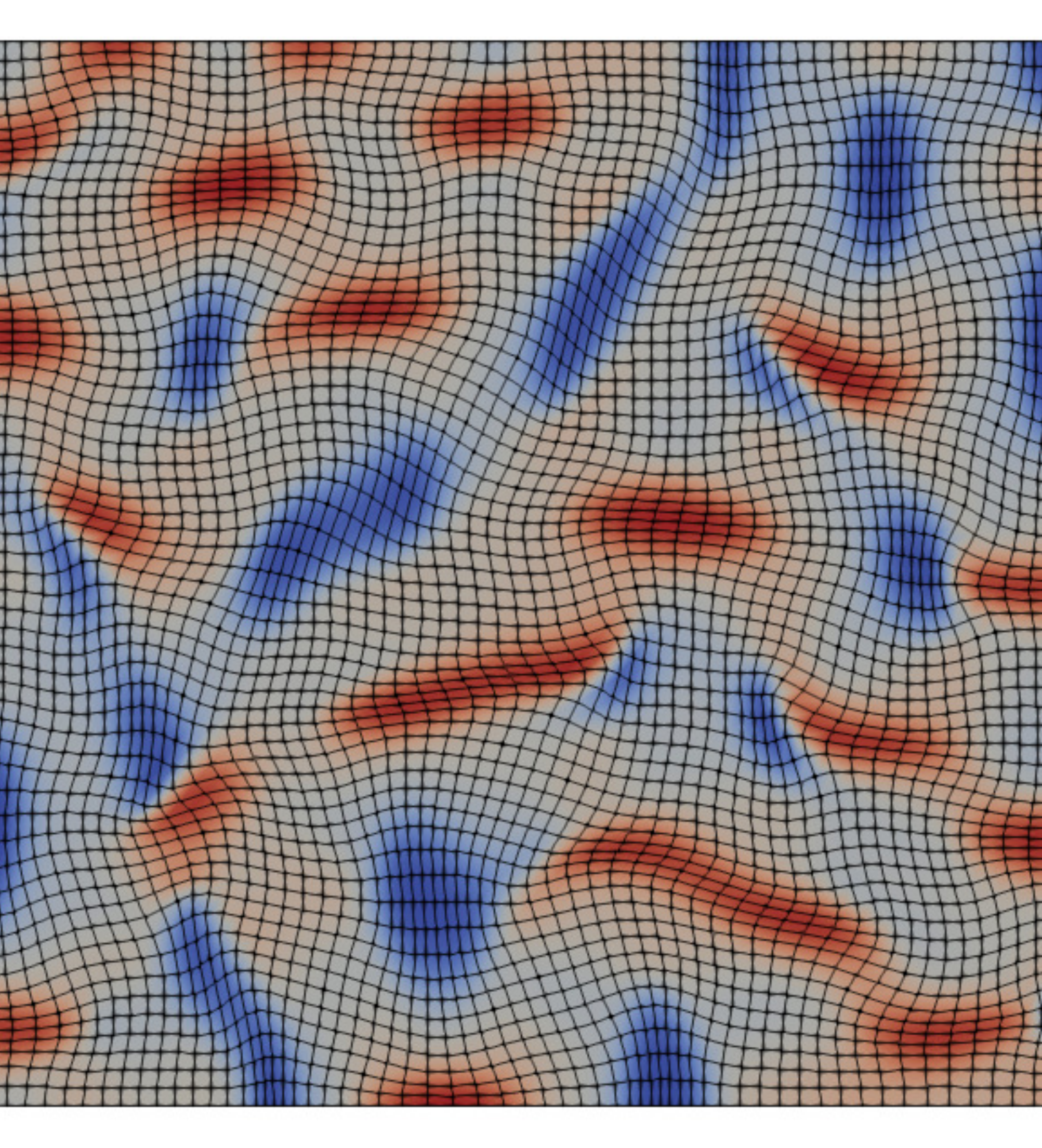}} \hfill
  \subfloat[$e_2$ at step 150]{\includegraphics[trim={15cm 3.5cm 15cm 3.5cm}, clip, width=0.24\linewidth]{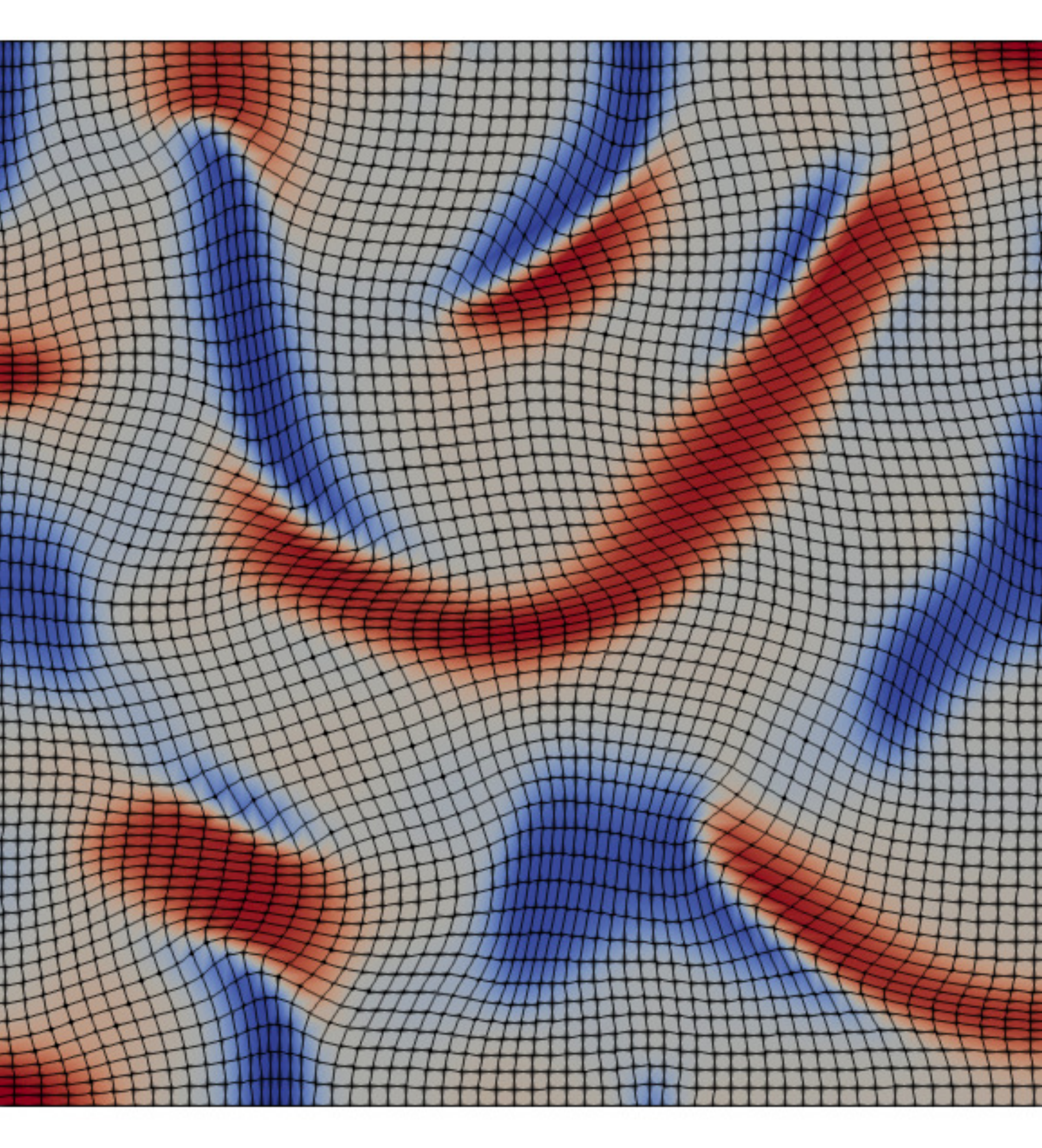}} \hfill
  \subfloat[$e_2$ at step 400]{\includegraphics[trim={15cm 3.5cm 15cm 3.5cm}, clip, width=0.24\linewidth]{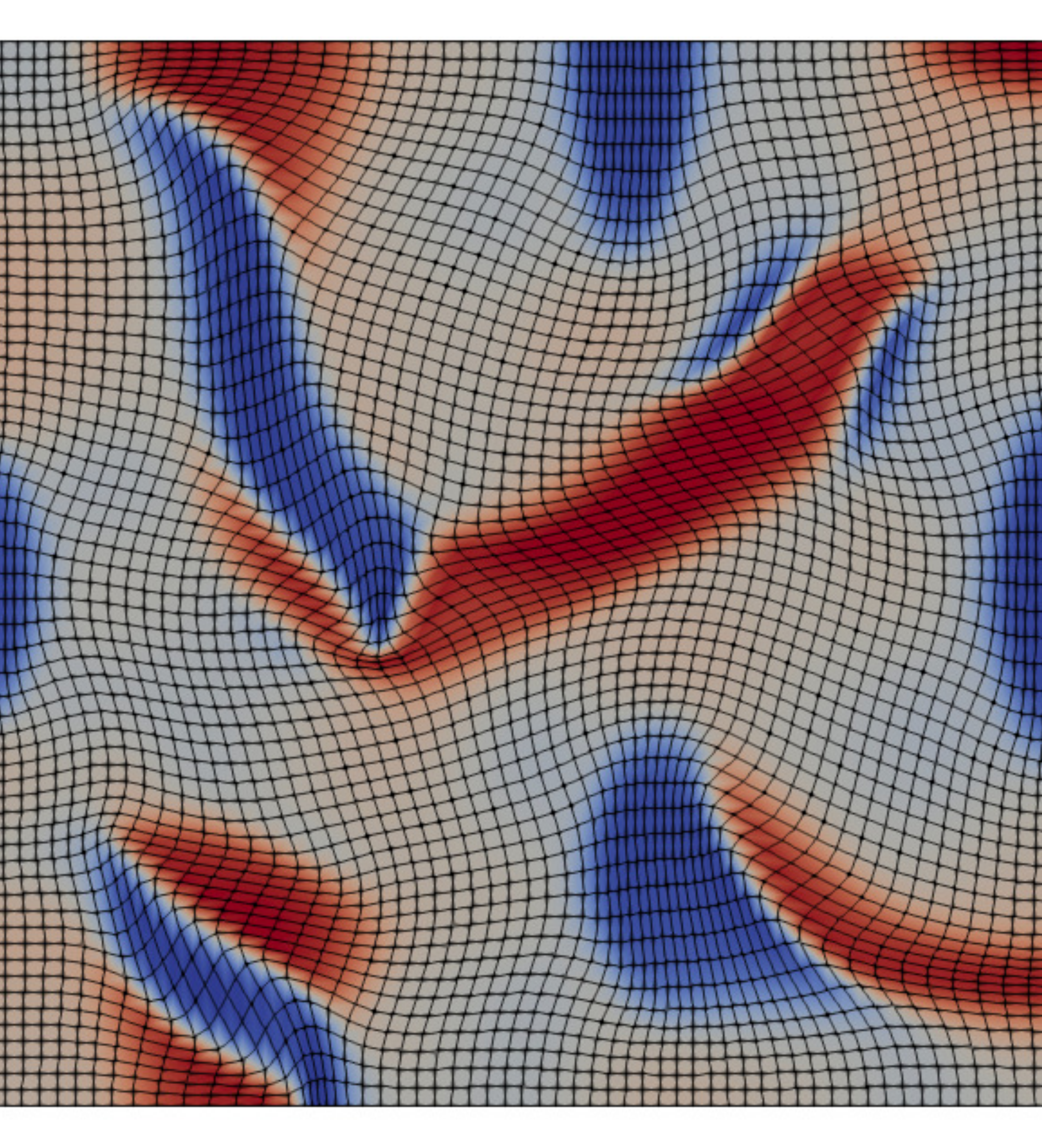}} \hfill
  \subfloat[$e_2$ at step 900]{\includegraphics[trim={15cm 3.5cm 15cm 3.5cm}, clip, width=0.24\linewidth]{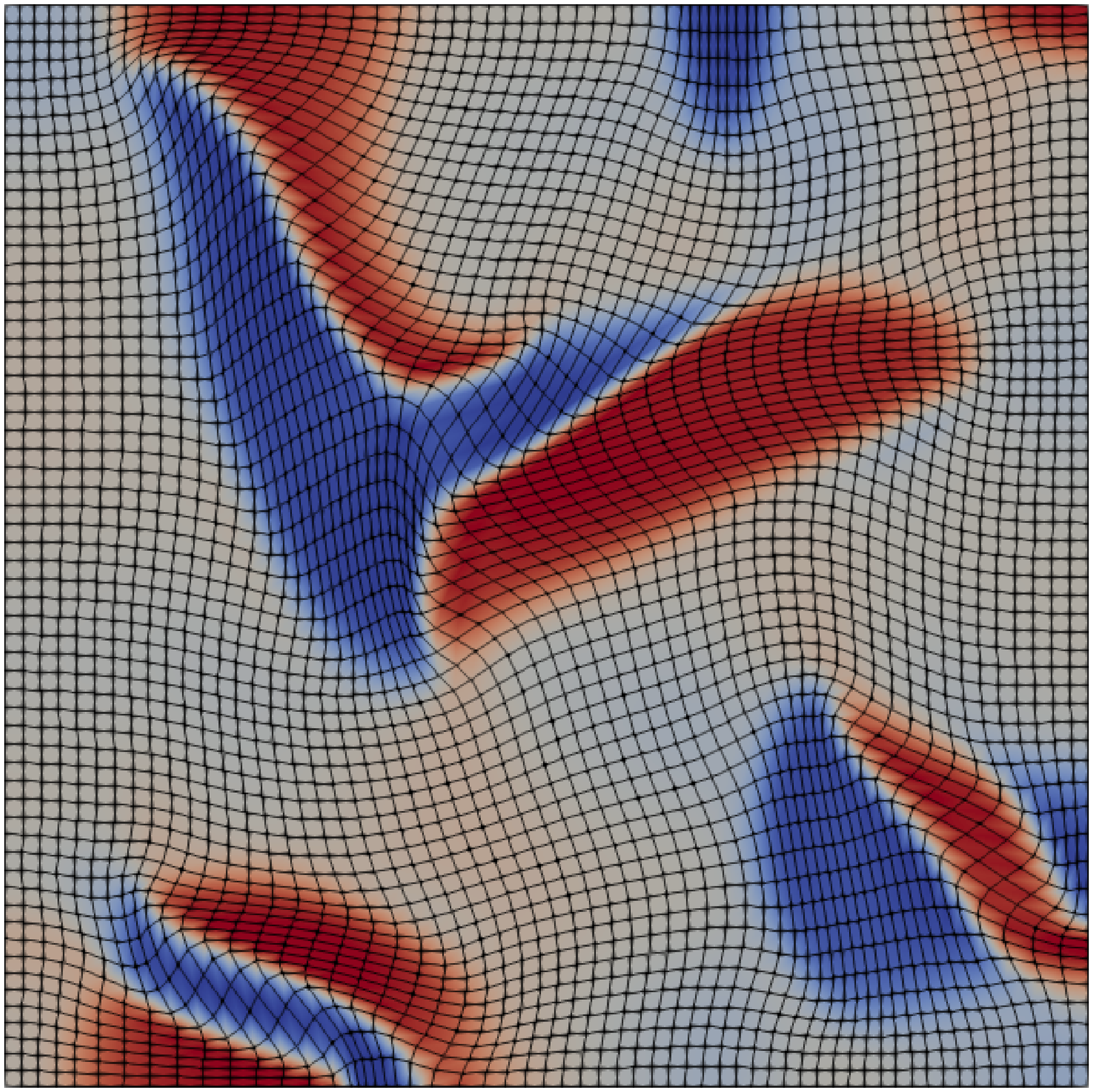}}
  \caption{Microstructures evolving under the coupled Cahn-Hilliard, gradient elasticity equations.
  (a-d) Composition field: red for $c=1$ and blue for $c=0$.
  (e-h) $e_2$ field: red for $e_2=0.1$ (``positive'' rectangle phase) and blue for $e_2=-0.1$ (``negative'' rectangle phase). 
In the region where $c=0$, $e_2$ has a value of $0.0$, corresponding to the square phase.}
  \label{fig:dns-results}
\end{figure}
\begin{figure}[hpt]
\centering
\includegraphics[width=0.6\textwidth]{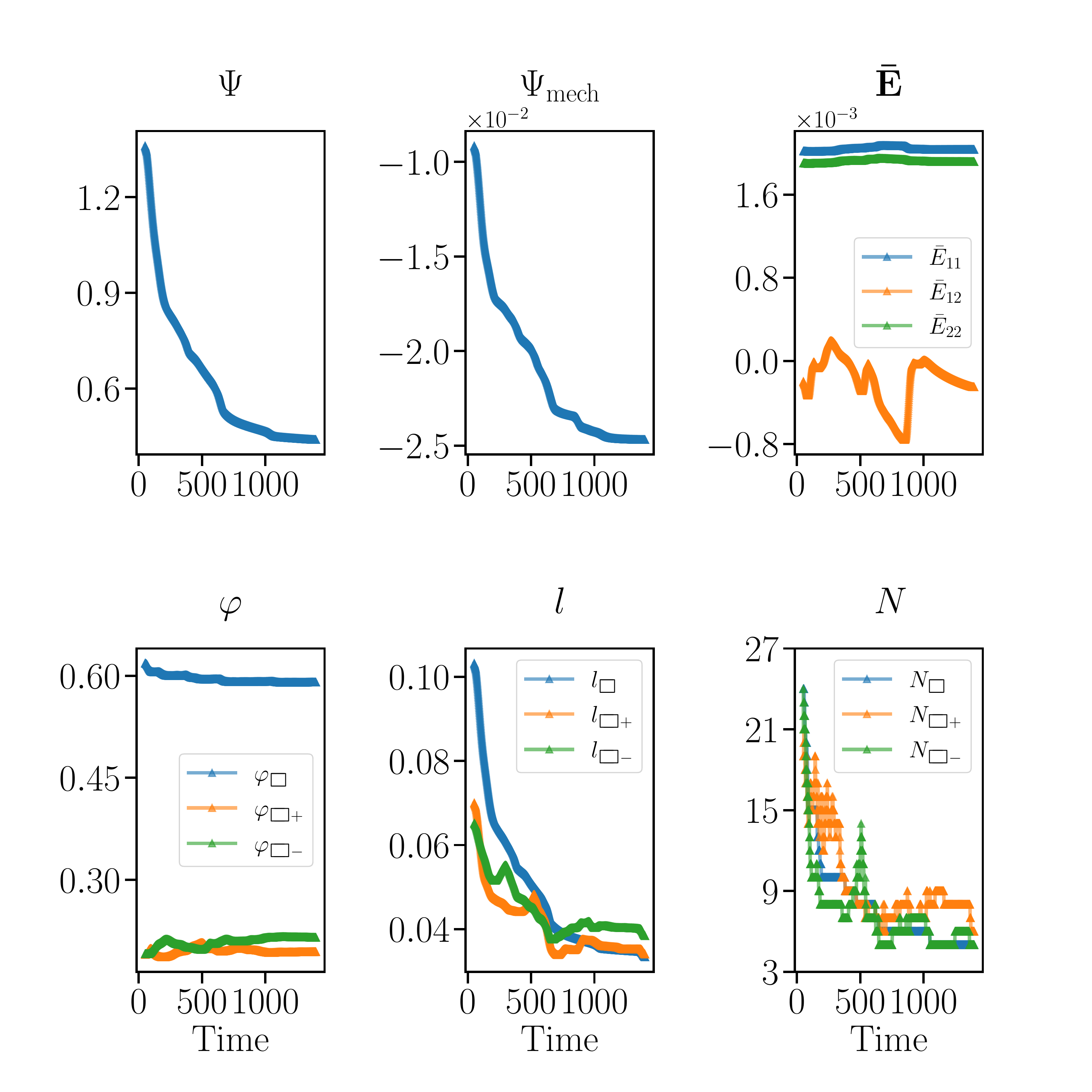}
\caption{DNS calculated state vector components. See the text for definitions.}
\label{fig:DNSobservables}
\end{figure}

These governing PDEs may now be solved following Zhang \& Garikipati, \cite{Zhang2020} (see \cref{fig:dns-results}) yielding time series data that are functionals of the high dimensional solutions. These quantities, denoted as barred variables unless otherwise specified, are effective volume averaged quantities $u \to \bar{u}$ over the material domain, with example direct numerical simulation observables shown in \cref{fig:DNSobservables} In this work, we consider a single initial and boundary value problem of a solid subjected to plane strain conditions in a two-dimensional domain $\Omega = (0, 0.01)\times (0, 0.01)$ with a mesh size of $60 \times 60$. The solid is loaded by a steady biaxial Dirichlet boundary conditions. The solid has a randomly fluctuating initial composition in the range of $c = 0.46 \pm 0.05$, with a uniform initial $e_2 = 0$ field, which corresponds to a single square phase that exists at high temperature. Zero chemical flux boundary conditions are applied to the solid. We study one trajectory of the evolution of the microstructure and aim to model its quantities of interest. In future communications we will consider the effectiveness of the graph theoretic approach at producing general models for a family of microstructures with different boundary and initial conditions. 

\noindent The resulting state vector for the proposed models is
\begin{equation}
	x = \{\Psi,\Psi_{\textrm{mech}} ,\vecbar{E},\varphi_{\alpha},l_{\alpha},N_{\alpha}\}, \label{eq:observables}
\end{equation}
where $\Psi = \int_\Omega\psi~\mathrm{d}V$, $\Psi_\textrm{mech}$ is the elastic component of $\Psi$ extracted from \cref{eq:psi_hom} and \cref{eq:psi_grad}. Vertices in this graph correspond to the time index. Due to the dissipative nature of phase field dynamics, the microstructure states at different times are all related, and therefore the graph is considered to be fully connected when considering edges and edge weights.

\subsubsection{Reduced-order modelling} \label{sec:physicalsystems_reducedordermodelling}
As with the Allen-Cahn example of \cref{sec:physicalsystems_diffusion}, we aim to obtain a reduced-order model for the dynamics in terms of the state vector \cref{eq:observables}. As the dynamics progresses and the microstructure evolves, see \cref{fig:dns-results}, data is generated in the form of fields $c,\vec{e}$ from which the state vector \cref{eq:observables} can be extracted. An analysis of \crefrange{eq:psi_hom}{eq:piola3} suggests which algebraic and differential terms should be proposed for the reduced-order model. The free energy contains no explicit spatial dependence, and contains monomials with powers up to $c^4$ and ${{e_{\alpha}}}^4$, plus coupled terms that are up to $c {e_{2}}^2$ order, plus quadratic gradient terms. Therefore the Cahn-Hilliard first order dynamics for the composition are of the order
\begin{align}
	\derivative[1]{c}{t} &\sim \mathcal{O}(\uniderivative[3]{\psi}{c} c^2) + \mathcal{O}(\derivative[3]{\psi}{c,c,{e_{\alpha}}} c {e_{\alpha}})
	+ \mathcal{O}(\derivative[3]{\psi}{c,{{e_{\alpha}}},{{e_{\beta}}}}~ {{e_{\alpha}}}{e_{\beta}}) \label{eq:cahnhilliard_approx}	\\
	&~ + \mathcal{O}(\uniderivative[2]{\psi}{c} c) + \mathcal{O}(\derivative[2]{\psi}{c,{e_{\alpha}}} {e_{\alpha}}) + \mathcal{O}(c) \nonumber.
\end{align}
A physics informed basis for the phase volume fraction dynamics and free energy can then be selected consisting of polynomials of terms such as those in \cref{eq:cahnhilliard_approx}, where all derivative terms are represented using the non-local calculus definitions in \cref{eq:diff_p} with a desired order of accuracy.

\noindent Upon selection of a physically relevant basis for the models, system inference methods, specifically linear stepwise regression\cite{Wang2019,Kochunas2020}, allow the model to be made as parsimonious as possible, while remaining an accurate representation of the underlying physics. Here, we choose to perform regression to compute linear coefficients $\gamma$ for each term in the basis using Ordinary Least Squares (OLS), as well as Ridge regression \cite{Lv2009,Tibshirani1996} schemes which constrain the norm of the coefficients. An important aspect of the fitting procedure is the choice of loss function to characterize the error between the known data and regression predictions, particularly for highly oscillatory data, such as this coupled mechano-chemical system, which displays rapidly evolving local phases during its transient response. We define the total loss function to be minimized over the stepwise regression as a weighted sum of loss functions $\ell = \sum_{\chi \in \{1,2,\infty,\cdots\}} w_{\chi}\cdot \ell_{\chi}$, where $\chi$ denotes the type of norm $L_{\chi}$.


\subsubsection{Phase volume fraction dynamics}\label{sec:physicalsystems_vol}
To model the phase volume fraction dynamics, we are guided by the Cahn-Hilliard equation in \cref{eq:pde_cahnhilliard}, and the order of terms present in \cref{eq:cahnhilliard_approx}, and propose a model for first order dynamics of the square phase $\bar{\varphi} = \varphi_{\square}$ that is dependent on polynomials of products of the free energy derivatives up to $\difference[3]{\Psi}{\bar{\varphi},\bar{E}_{\xi \kappa}}$, and phase dependent algebraic terms of  $\{\bar{\varphi},\bar{E}_{\xi \kappa},\bar{l}_{\alpha}, \bar{N}_{\alpha},\vecbar{E}\}$. This model takes the form
\begin{equation}
	\derivative[1]{\bar{\varphi}}{t} = f(\nderivative[q]{\Psi}{\bar{\varphi}}{\bar{E}_{\alpha \beta}},\bar{\varphi},\bar{E}_{\alpha \beta},\bar{l}_{\alpha}, \bar{N}_{\alpha}), \label{eq:model_vol}
\end{equation}
where $\nderivative[q]{\Psi}{\bar{\varphi}}{\bar{E}_{\alpha \beta}} = \{\derivative[2]{\Psi}{\bar{\varphi},\bar{E}_{\alpha \beta}},\derivative[2]{\Psi}{\bar{\varphi},\bar{\varphi}},\derivative[3]{\Psi}{\bar{\varphi},\bar{E}_{\alpha \beta},\bar{E}_{\xi \kappa}},\derivative[3]{\Psi}{\bar{\varphi},\bar{\varphi},\bar{E}_{\alpha \beta}},\uniderivative[3]{\Psi}{{\bar{\varphi}}}\}$, and $f$ is a polynomial function of the basis.

\noindent Stepwise regression is performed for a  third-order polynomial function, $f$, using a mixture of the $\ell_{2}$ and $\ell_{1}$ loss functions, for OLS and Ridge regression, where the optimal ridge parameter is found to be small $\lambda = 10^{-17}$. Each choice of loss weighting yields similar trends in the loss curves increasing with increased model parsimony, offset by how much the $\ell_{1}$ contributes to the total loss. We note that the high frequency oscillations in the phase volume fractions are fit well by the $\ell_{1}$ loss as seen in \cref{fig:DNSobservables}. We understand this in terms of the more stringent penalization that the $L_1$ norm imposes on the error, in comparison with, e.g. $L_2$.  When looking at the exact ordering of operators that the stepwise regression presents as most important to the model, the different loss weights do not result in identical orderings, even at low numbers of operators where the loss curves are more similar. However, the general trends of which operators are present at various stages of the stepwise regression, denoted by the various plateaus in the loss curves, are very similar across choices of loss weights. We therefore choose to focus on analyzing the fits where solely the $\ell_{1}$ loss is used.

\noindent With loss weights $w_2 = 0,~ w_1 = 1$, the first 10 terms from the stepwise regression for the first order dynamics of $\bar{\varphi}$ are found using OLS and Ridge regression to be
\input{./figures/model_vol_OLS_Ridge17_local.tex}
\noindent It can be seen that both models consist of solely algebraic terms. The OLS and Ridge approaches also generate models with similar terms, with the Ridge regression model having 9 terms compared to 8 terms for OLS regression that contain predominantly strain dependencies. 

\begin{figure}[hpt]
\centering
\begin{subfigure}[t]{0.49\textwidth}
	\centering
	\includegraphics[width=\textwidth]{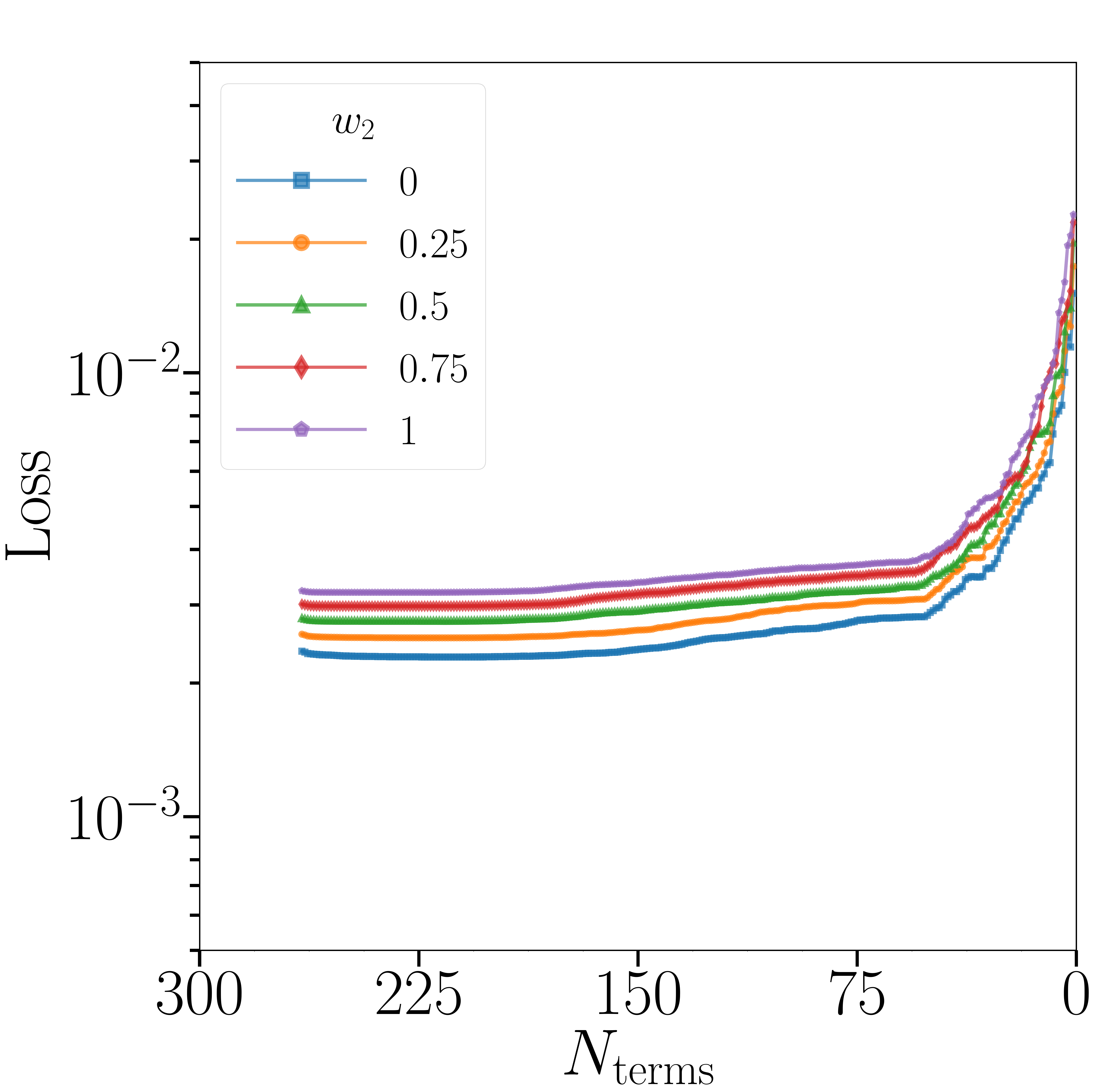}
	\subcaption{OLS regression.}
	\label{fig:vol_loss_OLS_local}
\end{subfigure}
\hfill
\begin{subfigure}[t]{0.49\textwidth}
	\centering
	\includegraphics[width=\textwidth]{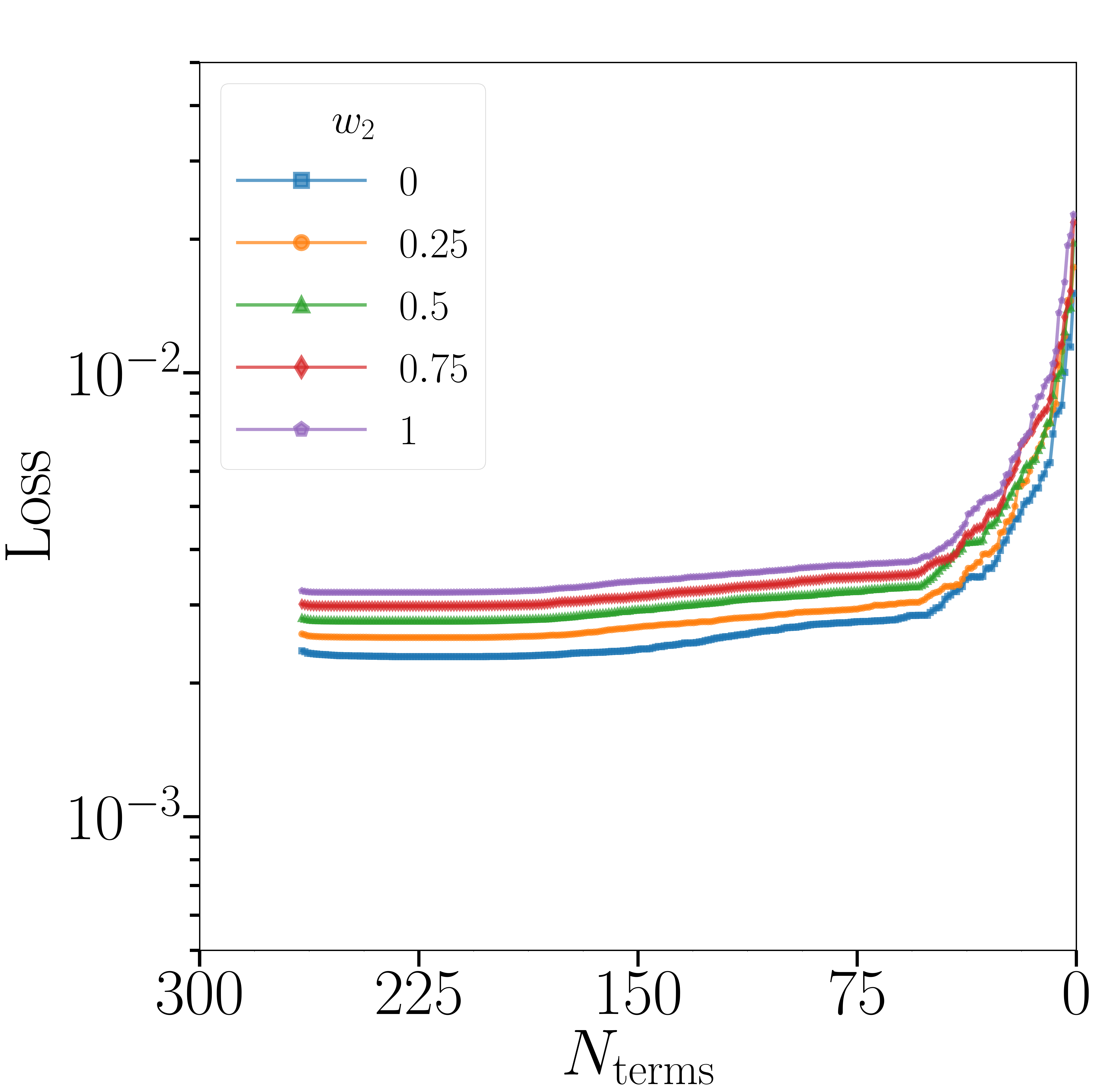}
	\subcaption{Ridge regression with $\lambda = 10^{-17}$.}
	\label{fig:vol_loss_Ridge17}
\end{subfigure}
\caption{Stepwise regression loss curves for the phase volume fractions first order dynamics, labeled by their weighting of the $\ell_{2}$ loss in the weighted residual loss function.}
\label{fig:vol_loss_OLS_Ridge17_local}
\end{figure}

\begin{figure}[h]
\centering
\begin{subfigure}[t]{0.49\textwidth}
	\centering
	\includegraphics[width=\textwidth]{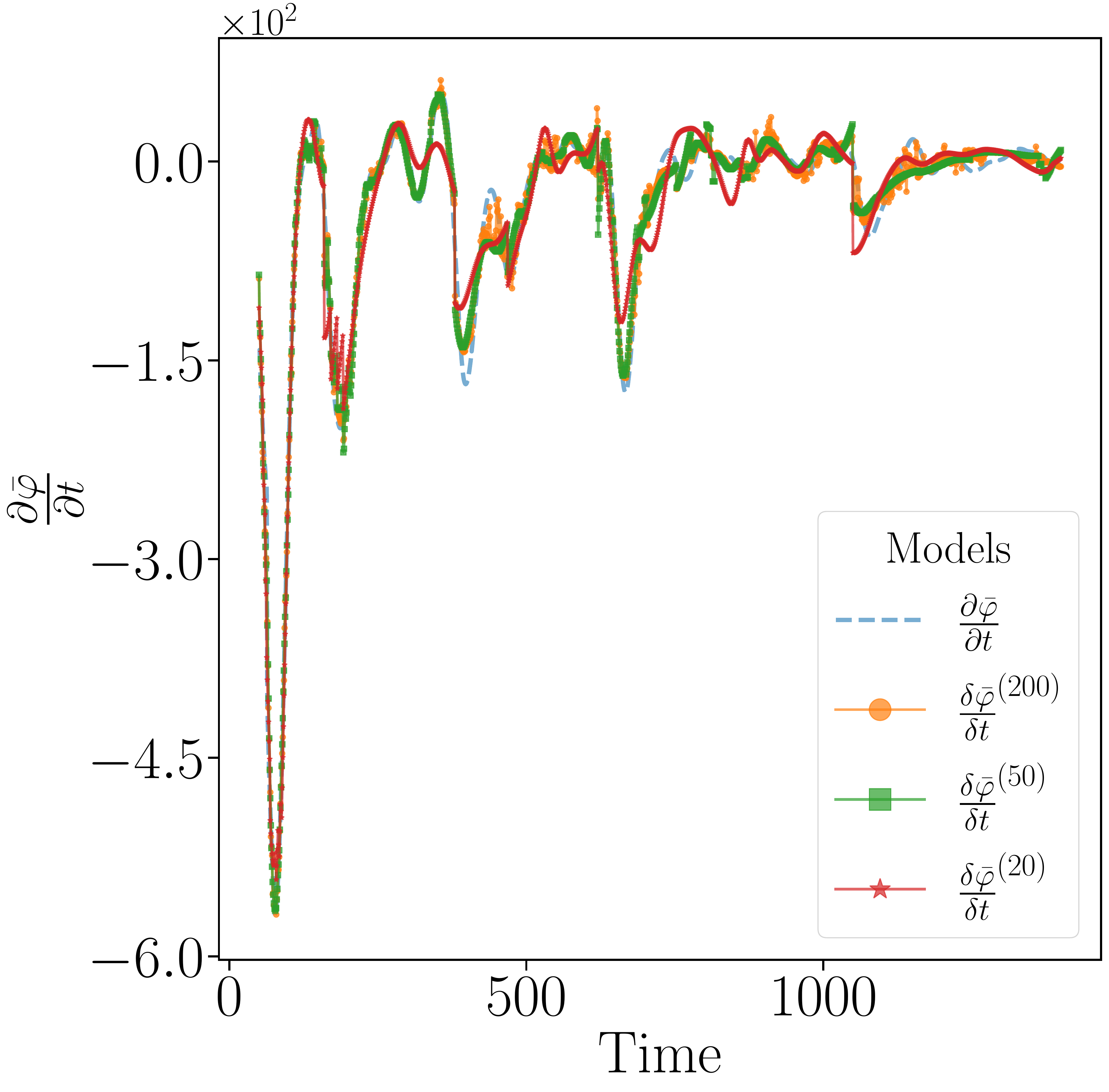}
	\subcaption{OLS regression.}
	\label{fig:vol_fit_OLS_local}
\end{subfigure}
\hfill
\begin{subfigure}[t]{0.49\textwidth}
	\centering
	\includegraphics[width=\textwidth]{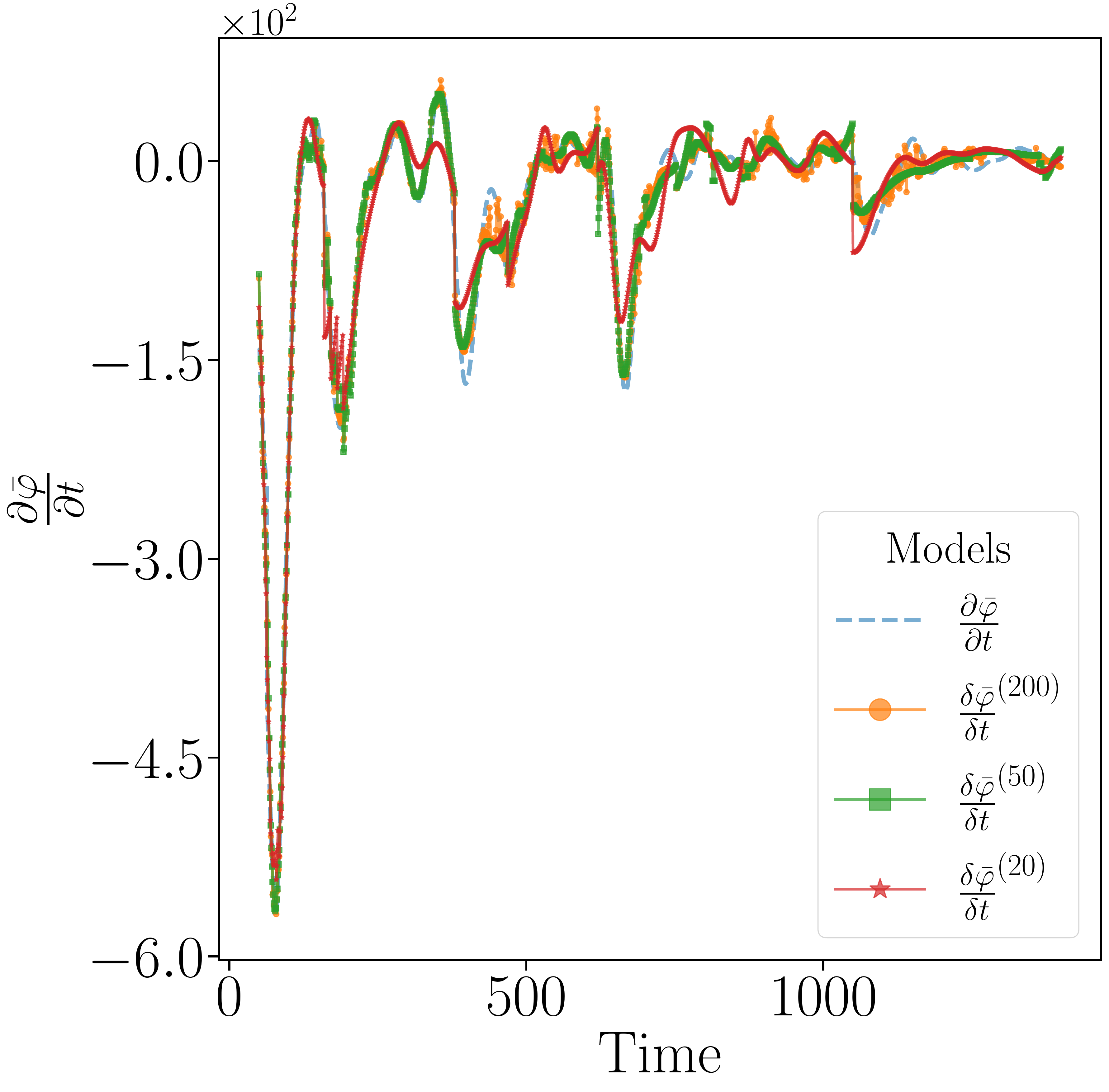}
	\subcaption{Ridge regression with $\lambda = 10^{-17}$.}
	\label{fig:vol_fit_Ridge17_local}
\end{subfigure}
\caption{First order dynamics of phase volume fractions for $20$, $50$, and $200$ terms in the forcing function, and use of the $\ell_{1}$ loss in stepwise regression. Backward Euler first time derivative data is shown with the dashed blue curve.}
\label{fig:vol_fit_OLS_Ridge17_local}
\end{figure}

\noindent Looking at the plateaus and sharp increases in the weighted loss curves, particularly for the OLS model in \cref{fig:vol_loss_OLS_local} over the stepwise iterations, the plateau extends from $265$ down to approximately past $70$ operators. The first order dynamics, and fits for models with $20$, $50$, and $200$ terms is shown in \cref{fig:vol_fit_OLS_Ridge17_local}, with both OLS and Ridge regression yielding comparable results. For the practically full model with $200$ terms, the fits correspond to the data very well at early times, but have difficulty matching all of the small oscillations at late times. As the model becomes more parsimonious to $50$ terms where the sharp increase in the loss occurs, the general trend of the dynamics is captured, however this model misses capturing the full amplitude of the major peaks compared to the more complex models. Finally, with $20$ terms, which corresponds to the number of terms present in the original Cahn-Hilliard equation in $d=2$ dimensions, the large peaks, and small oscillations are not entirely captured, however the general trends are still shown. The algebraic polynomial terms are potentially over-fitting, however manage to retain the general behavior of the phase volume fractions changing over time. 

\noindent In evaluating the usefulness of the reduced-order models, the number of terms should be compared with the high-dimensional PDE solutions that had $\mathcal{O}(10^4)$ degrees of freedom. \cite{Zhang2020} The dominance of the solely algebraic operators in the reduced-order models with 10 terms is clearly a case of overfitting with  models that are inadequate representations of the physics. The models with between 20 and 50 terms contain many more of the derivative operators that would be expected from \cref{eq:cahnhilliard_approx}. These models also do resolve the oscillations well, if not the peaks and troughs in \crefrange{fig:vol_fit_OLS_local}{fig:vol_fit_Ridge17_local} The reduced-order models with operator numbers in this range perhaps should be understood as incorporating some of the physics with the oscillatory derivative operators, while relying on the smoother algebraic terms to follow the mean trends. These trade-offs warrant further study in the context of the graph theoretic approach presented here. Importantly, however, this treatment admits these examinations of interpretability via explicit representation of non-local derivative terms.



\section{Conclusions}
\label{sec:concl}
In this work, we develop a graph theoretic approach for reduced-order modelling of physical systems. By defining a non-local calculus with small local neighborhoods of vertices, we show that derivatives of any order can have any desired order of accuracy to their differential counterparts, without any assumptions or symmetry required for the underlying graph. These well behaved derivatives allow physics informed bases to be selected, and reduced-order models found for any system of interest.

Applying this calculus to studying volume averaged quantities from high dimensional direct numerical simulations by representing these quantities on a graph, shows the intuitive, general, and effective nature of this approach. With example physical system of phase evolution of a chemical (the Allen-Cahn example) and a mechano-chemical multi-crystalline solid (the Cahn Hilliard-non-convex gradient elasticity example), we see that using a physics informed basis of operators can resolve functional representations well enough to be realized in reduced-order, time-dependent ODE models. We note that the smoother evolution of the Allen-Cahn PDE in \cref{sec:physicalsystems_diffusion} was represented with as few as three terms in the reduced-order form, while the more oscillatory and complex coupled Cahn Hilliard-non-convex gradient elasticity example \cref{sec:physicalsystems_microstructures} needed 20-50 terms.  However,  the original Cahn-Hilliard PDE itself has nearly 20 terms when fully expanded \cref{eq:cahnhilliard_approx}, and this number could be regarded as a lower bound for the complexity of the reduced-order models. As we also observed in \cref{sec:physicalsystems_microstructures}, these reduced-order models, while having tens of terms should be compared with the $\mathcal{O}(10^4)$ degrees of freedom in the generating PDEs. The explicit derivatives in these models, using the non-local calculus on graphs, confers interpretability upon them. Finally, the losses, on the order of $\mathcal{O}(10^{-4})$, attained in the model identification should be compared with recent studies of the same microstructure system using neural networks \cite{Zhang2020}, which achieve losses on the order of $\mathcal{O}(10^{-6})$.

\section*{Supplementary information}
\noindent Supplementary information is available as a separate document. Data and code used to produce this work can be found at \url{https://github.com/sidsriva/graph_calculus}.

\section*{Acknowledgments}
\noindent The authors gratefully acknowledge the support of the National Science Foundation, United States grant \textbf{\#1729166}. "DMREF/GOALI: Integrated Framework for Design of Alloy-Oxide Structures".

\pagestyle{bibliography}
\bibliography{main}

 
\newpage
\begin{appendices}  \label{sec:appendix}
\appendix
\renewcommand\thesection{\Alph{section}}
\renewcommand\theequation{\getcurrentref{section}.\arabic{equation}}
\renewcommand\thefigure{\getcurrentref{section}.\arabic{figure}}
\renewcommand\thetable{\getcurrentref{section}.\arabic{table}}

\section{Error analysis of local weights model}\label{app:error_analysis}
\setcounter{equation}{0}
\setcounter{figure}{0}
\setcounter{table}{0}
\subsection{Modified Taylor series model} \label{app:error_analysis_model}

An error analysis will be conducted for a function $u(x)$ of a $p$ dimensional variable $x = \{x^{\mu}\}$, represented by a modified $k$ order Taylor series $u_k(x)$ with non-local derivatives. 

The Taylor series functional representations of the $K$ order function $u(x)$ and $k\leq K$ order model $u_k(x)$, based at a point $\widetilde{x}$, are in one dimension
\begin{align}
	u(x|\widetilde{x}) =&~ u(\widetilde{x}) + \derivative[1]{u(\widetilde{x})}{x}(x-\widetilde{x}) + \frac{1}{2!}\uniderivative[2]{u(\widetilde{x})}{x}(x-\widetilde{x})^2 + \cdots + \frac{1}{K!}\uninderivative[K]{u(\widetilde{x})}{x}(x-\widetilde{x})^{K}, \label{eq:differentialtaylor} \\
	u_k(x|\widetilde{x}) =&~ u(\widetilde{x}) + \gamma_{1}(\widetilde{x})\unidifference[1]{u(\widetilde{x})}{x}(x-\widetilde{x}) + \frac{\gamma_{2}(\widetilde{x})}{2!}\unidifference[2]{u(\widetilde{x})}{x}(x-\widetilde{x})^2 + \cdots + \frac{\gamma_{k}(\widetilde{x})}{k!}\unindifference[k]{u(\widetilde{x})}{x}(x-\widetilde{x})^k \label{eq:modeltaylor_analysis_p}, 	
\end{align}
and in higher dimensions
\begin{align}
	u(x|\widetilde{x}) =&~ u(\widetilde{x}) + \sum_{\mu}\derivative[1]{u(\widetilde{x})}{x^{\mu}}(x-\widetilde{x})^{\mu} + \sum_{\mu\nu}\frac{1}{2!}\derivative[2]{u(\widetilde{x})}{x^{\mu},x^{\nu}}(x-\widetilde{x})^{\mu\nu} + \cdots \label{eq:differentialtaylor_p}\\
	&~~~~~~~+~ \sum_{\mu_{0}\cdots\mu_{K-1}}\frac{1}{K!}\nderivative[K]{u(\widetilde{x})}{x^{\mu_{0}}}{x^{\mu_{K-1}}}(x-\widetilde{x})^{\mu_{0}\cdots\mu_{l-1}}, \nonumber\\
	u_k(x|\widetilde{x}) =& u(\widetilde{x}) + \sum_{\mu}\gamma_{1}^{\mu}(\widetilde{x})\difference[1]{u(\widetilde{x})}{x^{\mu}}(x-\widetilde{x})^{\mu} + \sum_{\mu\nu}\frac{\gamma_{2}^{\mu\nu}(\widetilde{x})}{2!}\difference[2]{u(\widetilde{x})}{x^{\mu},x^{\nu}}(x-\widetilde{x})^{\mu\nu} + \cdots \label{eq:modeltaylor_p} \\
	&~~~~~~+~ \sum_{\mu_{0}\cdots\mu_{k-1}}\frac{\gamma_{k}^{\mu_{0}\cdots\mu_{k-1}}(\widetilde{x})}{k!}\ndifference[k]{u(\widetilde{x})}{x^{\mu_{0}}}{x^{\mu_{k-1}}}(x-\widetilde{x})^{\mu_{0}\cdots\mu_{k-1}} \nonumber, 	
\end{align}
where the coefficients $\gamma_{l}(\widetilde{x})$ are the linear coefficients fit to the non-local derivative model based around $\widetilde{x}$. For this analysis, we denote products of $l$ powers of vector elements as
\begin{align}
	\xi^{\mu_{0}^{t_{0}}\cdots\mu_{l-1}^{t_{l-1}}} =&~ \prod_{s=0}^{l-1}\xi^{\mu_{s}^{t_{s}}}
\end{align}
for any vector $\xi$, indices $\{\mu_{q}\}$ and powers $\{t_{q}\}$.

\noindent The model is comprised of non-local derivatives, defined as the differential derivatives plus an error term
\begin{align}
	\unindifference[l]{u}{x} =&~ \uninderivative[l]{u}{x} + \varepsilon_l,
\end{align}
and the linear coefficients $\gamma$ that are fit given data. The behavior of the derivatives and coefficients in the model will depend on the non-local nature of the derivatives, represented by weights in the graph, and the choice of data used to fit the model. 

\noindent The true Taylor series of the function has the linear coefficients being identically $1$ and the derivatives being the true differential derivatives. Therefore how the derivative and coefficient differ from their values in a true Taylor series, and how these errors scale with the spacing of the data will determine the scaling of the local and global error.

\subsection{Data mesh} \label{app:error_analysis_mesh_p}

\begin{figure}[hpt]
\centering
\includegraphics[width=0.4\textwidth]{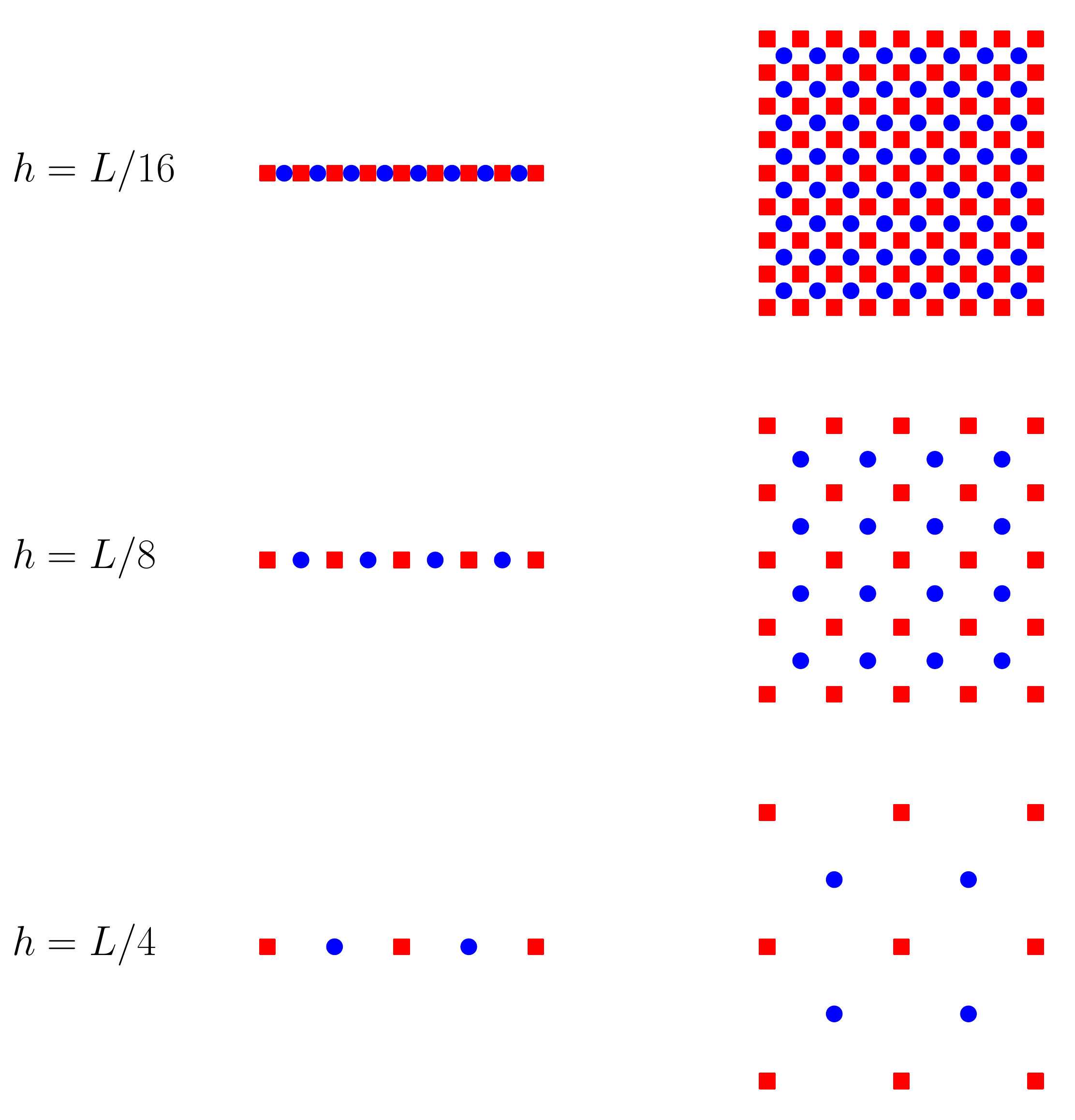}
\caption{Interlacing of red square training points and blue circle testing points for $p=1$ (left) and $p=2$ (right) dimensional datasets. Here $n$ is a power of 2, the domain length $L = 2nh$ is fixed, and $n+1$ training points are placed along each dimension to ensure the boundaries of the domain are always included for all spacings $h$.}
\label{fig:meshes_analysis_p}
\end{figure}

The graph vertices are partitioned into a training data $\widetilde{V}$ with $(n+1)^p$ uniformly spaced points, interlaced by testing data $V$ with $n^p$ uniformly spaced points, as per \cref{fig:meshes_analysis_p}. The data points are over an interval $L$ along each of the $p$ dimensions, with uniform spacing 
\begin{align}
	h =&~ \frac{L}{2n} \label{eq:hspacing_p}, 
\end{align}
giving $N = (n+1)^p + n^p$ total points. Each point will be assigned the vertex index vector $j$ such that
\begin{align}
	x_j = 2h(j^0,\cdots,j^{p-1}).
\end{align}
For simplicity of showing explicit scaling, we denote $x_j = \frac{j}{n}L = \varchi_j L$, where $\varchi_j \in [0,1]$ is fixed for any $n$.

We will also use the slight abuse of notation $x \in V$ to denote the point $x$ at a vertex in the set of vertices $V$, and denote sums of functions $f(x)$ over a set of vertices as
\begin{align}
	\sum_{V} f =&~ \sum_{j \in V} f(x_j).
\end{align}
We may define around each vertex $\widetilde{x} \in \widetilde{V}$ a neighborhood of adjacent vertices
\begin{align}
	\mathcal{N}(\widetilde{x}) \subseteq&~ \widetilde{V},
	\intertext{where}
	d = d(\widetilde{x})&~ = \abs{\mathcal{N}(\widetilde{x})}
\end{align}
is the size of the neighborhood, and also can be thought of as the total degree of that vertex. 

\subsection{Local error} \label{app:error_analysis_local}
We will now investigate the local pointwise scaling of the error in spacing $h$ of the modified Taylor series model and we denote the components of the error in the follow definitions. The objective of this analysis is to determine 
\begin{align}
	r_{\textrm{global}},~r_{\textrm{local}},~r_{\textrm{der}_{l}},~\textrm{and}~ r_{\textrm{coef}_l} \nonumber
\end{align}
which are the leading order scaling $\mathcal{O}(h^{r})$ of the global model error, the local model error, the local derivative error, and local coefficient error, for a given $l$ order of derivative. They each depend on the dimension $p$, the order of the model $k$, the order of the function $K$, and the data $\widetilde{V}$ used in fitting the model.

\noindent We now define the exact forms of the errors in these components of the model, which we write as sums over powers in $h$, with minimum scalings $r$. 

\noindent The local error for a base point $\widetilde{x}$ in the model is defined as
\begin{align}
	e(x|\widetilde{x}) = u_k(x|\widetilde{x}) -&~ u(x|\widetilde{x}) = \sum_{q \geq r_{\textrm{local}}} C_q(x|\widetilde{x}) h^{q}. \\
	\intertext{The local error in the non-local derivatives is defined in one dimension as}
	\varepsilon_l(\widetilde{x}) = \unindifference[l]{u(\widetilde{x})}{x} -&~ \uninderivative[l]{u(\widetilde{x})}{x} = \sum_{q \geq r_{\textrm{der}_{l}}} E_{l_q}(\widetilde{x})h^{q} \\
	\intertext{and in higher dimensions as}
	\varepsilon_l^{\mu_{0}\cdots\mu_{l-1}}(\widetilde{x}) = \ndifference[l]{u(\widetilde{x})}{x^{\mu_{0}}}{x^{\mu_{l-1}}} -&~ \nderivative[l]{u(\widetilde{x})}{x^{\mu_{0}}}{x^{\mu_{l-1}}} = \sum_{q \geq r_{\textrm{der}_{l^{\mu_{0}\cdots\mu_{l-1}}}}} E_{l_q}^{\mu_{0}\cdots \mu_{l-1}} (\widetilde{x})h^{q}. \\
	\intertext{The error in the linear coefficients in one dimension is}
	\gamma_l(\widetilde{x}) - 1 =&~ \sum_{q \geq r_{\textrm{coef}_l}} G_{l_q}(\widetilde{x}) h^{q} \\
	\intertext{and in higher dimensions as}
	\gamma_l^{\mu_{0}\cdots\mu_{l-1}}(\widetilde{x}) - 1 =&~ \sum_{q \geq r_{\textrm{coef}_{l^{\mu_{0}\cdots\mu_{l-1}}}}} G_{l_q}^{\mu_{0}\cdots\mu_{l-1}}(\widetilde{x}) h^{q}.
\end{align}
Here the $\abs{C_l(x|\widetilde{x})} \leq C_l,~\abs{E_{k_l}(\widetilde{x})}\leq E_{k_l},\textrm{ and } \abs{G_{k_l}(\widetilde{x})} \leq G_{k_l}$ are constant with respect to $h$ and are assumed to have upper bounds over the data domain. 

\noindent The local error at $x \in V$, for a given $\widetilde{x} \in \widetilde{V}$ to base the Taylor series in $p=1$ dimensions is
\begin{align}
	e(x | \widetilde{x}) =&~ \sum_{l=0}^{k}\frac{\gamma_l(\widetilde{x}) - 1}{l!} \uninderivative[l]{u(\widetilde{x})}{x} (x-\widetilde{x})^l \\
	~+&~ \sum_{l=0}^{k} \frac{\gamma_l(\widetilde{x})}{l!} \varepsilon_l(\widetilde{x}) (x-\widetilde{x})^l \nonumber \\
	~-&~ \sum_{l=k+1}^{K}\frac{1}{l!} \uninderivative[l]{u(\widetilde{x})}{x} (x-\widetilde{x})^{l} \nonumber \\
	\sim&~ \sum_{l=1}^{k}\mathcal{O}(h^{r_{\textrm{coef}_l}+l}) + \sum_{l=1}^{k}(1 + \mathcal{O}(h^{r_{\textrm{coef}_l}}))\mathcal{O}(h^{r_{\textrm{der}_{l}}+l}) + \mathcal{O}(h^{k+1}) \\
	=&~ \mathcal{O}(h^{r_{\textrm{local}}})
\end{align}
and in higher dimensions is
\begin{align}
	e(x | \widetilde{x}) =&~ \sum_{l=0}^{k}\sum_{\mu_{0}\cdots\mu_{l-1}}\frac{\gamma_l^{\mu_{0}\cdots\mu_{l-1}}(\widetilde{x}) - 1}{l!} \nderivative[l]{u(\widetilde{x})}{x^{\mu_{0}}}{x^{\mu_{l-1}}} (x-\widetilde{x})^{\mu_{0}\cdots\mu_{l-1}} \\
	~+&~ \sum_{l=0}^{k}\sum_{\mu_{0}\cdots\mu_{l-1}} \frac{\gamma_l^{\mu_{0}\cdots\mu_{l-1}}(\widetilde{x})}{l!} \varepsilon_l^{\mu_{0}\cdots\mu_{l-1}}(\widetilde{x}) (x-\widetilde{x})^{\mu_{0}\cdots\mu_{l-1}} \nonumber \\
	~-&~ \sum_{l=k+1}^{K}\sum_{\mu_{0}\cdots\mu_{l-1}}\frac{1}{l!} \nderivative[l]{u(\widetilde{x})}{x^{\mu_{0}}}{x^{\mu_{l-1}}} (x-\widetilde{x})^{\mu_{0}\cdots\mu_{l-1}} \nonumber \\
	\sim&~ \sum_{l=1}^{k}\sum_{\mu_{0}\cdots\mu_{l-1}} \mathcal{O}(h^{r_{\textrm{coef}_{l^{\mu_{0}\cdots\mu_{l-1}}}}+l}) \\
	+&~ \sum_{l=1}^{k}\sum_{\mu_{0}\cdots\mu_{l-1}} \left(1 + \mathcal{O}(h^{r_{\textrm{coef}_{l^{\mu_{0}\cdots\mu_{l-1}}}}})\right) \mathcal{O}(h^{r_{\textrm{der}_{l^{\mu_{0}\cdots\mu_{l-1}}}}+l}) \\
	+&~ \mathcal{O}(h^{k+1}) \\
	=&~ \mathcal{O}(h^{r_{\textrm{local}}})
\end{align}

\noindent We see that the local error in the method is due to the non-local derivatives having an error $\varepsilon_l \neq 0$ from the differential derivatives, due to the fit linear coefficient $\gamma_l(\widetilde{x}) \neq 1$ being not identically equal to $1$ like in a differential calculus Taylor series, and due to the presence of higher order terms in the function when $K>k$. To complete this scaling analysis, we must determine how the fit linear coefficients are dependent on the training data.

\subsection{Number of terms in polynomial basis in higher dimensions}\label{app:nonlocalerror_number_terms}
We will briefly review how to compute the combinatorics of the number of basis terms $q(p,r)$ for an order $r$ polynomial basis in $p$ dimensions. For a given order $k \in \{0,1,\dots,r\}$, a term in this basis can be thought of as a product of $p$ monomials $x^{\mu}$, each with power $a_{\mu}$:
\begin{align}
	x^{a} \equiv x^{{0}^{a_{0}}} \cdots x^{{p-1}^{a_{p-1}}} ~:~ a_{0} + \cdots + a_{p-1} = k ~,~ 0 \leq a_{\mu} \leq k.
\end{align}
The polynomial basis representation of a function $f_{p,r}(x)$, with coefficients $\alpha$, therefore has the form
\begin{align}
	f_{p,r}(x;\alpha) =&~ \sum_{k=0}^{r} \sum_{a_{0} + \cdots + a_{p-1} = k} \binom{k}{a_0 \cdots a_{p-1}} \alpha_{k}^{a_0\cdots a_{p-1}}~ x^{{0}^{a_{0}}} \cdots x^{{p-1}^{a_{p-1}}}
\end{align}
where the multinomial coefficients $\binom{k}{a_0 \cdots a_{p-1}} = \frac{k!}{a_{0}! \cdots a_{p-1}!}$ represent the number of ways of distributing a total power of $k$ across the $p$ dimensions, given the fixed numbers $\{a_{0},\cdots,a_{p-1}\}$. The total number of terms in the basis is then
\begin{align}
	q(p,r) =&~ \sum_{k=0}^{r} q_k(p).
\end{align}
We will consider either the case of treating each permutation of the order of the $\{a_{0},\cdots,a_{p-1}\}$ in each term as distinct, yielding many non-unique terms to count, or the case of counting each permutation as identical and only counting unique terms. For example when $k=2$ and $\mu \neq \nu$, one term may be $x^{\mu}x^{\nu}$, and whether it is counted as being distinct from $x^{\nu}x^{\mu}$ affects how we count the number of terms.

\noindent For each $k \in \{0,1,\dots,r\}$, if we consider the first case of counting all permutations of the ordering of the $\{a_{0},\cdots,a_{p-1}\}$ as distinct, then the total number of terms with order $k$ is
\begin{align}
	q_k(p) =&~ \sum_{a_{0} + \cdots + a_{p-1} = k} \binom{k}{a_0 \cdots a_{p-1}} = p^{k} ~\textrm{(non-unique)}.
\end{align}
If consider only the unique terms, the stars and bars method can be used to show the total number of terms with order $k$ is
\begin{align}
	q_k(p) =&~ \binom{p+k-1}{k} ~\textrm{(unique)}.
\end{align}
The combinatorics of the total non-unique number of terms can be calculated by summing the geometric series of $q_k$ for $k$ from $0$ to $r$. The combinatorics of the total unique number of terms can be done using the stars and bars method, or by summing the above terms of $q_k$ and proving through induction.

\begin{figure}[hpt]
	\centering
	\includegraphics[width=0.8\textwidth]{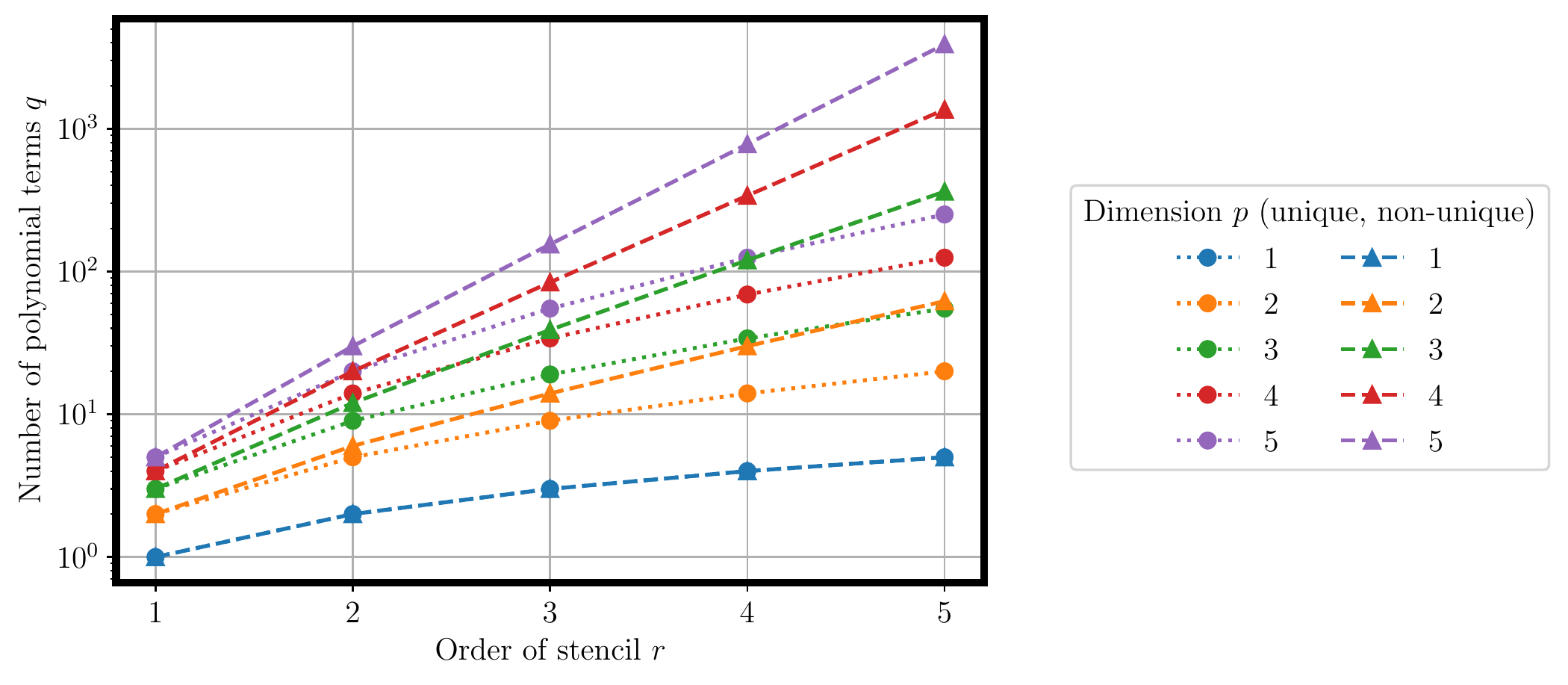}
	\caption{Number of polynomial terms when considering unique and non-unique terms due to commutativity of multiplication of monomials in each term.}
\label{fig:npolyterms}
\end{figure}

\noindent There are therefore
\begin{align}
	q(p,r) =&~ \frac{p^{r+1}-1}{p-1} ~\textrm{(non-unique)}~ \to \frac{(p+r)!}{p!~r!} ~\textrm{(unique)} \label{eq:npolyterms}
\end{align}
terms in a polynomial basis in high dimensions. Please refer to \cref{fig:npolyterms} for a comparison of the number of terms for a certain order $r$, for each possible dimension $p$.

\noindent In general, we must therefore choose our neighborhoods to have $d\geq q$ points for our systems of constraints on our weights to be well posed. For example, if the symmetric set of up to $l$ nearest neighbors on a uniform Cartesian mesh are considered at vertices far from the boundaries of the graph, then there will be $d = d(l) \geq 2pl$ neighbors for some integer $l \geq 1$ such that $d \geq q$. For most $p$ and $r$, $d > q$, and for completely structured data, or in particular cases where there happen to be equidistant neighboring vertices, the choice must be made either to solve this under-determined system, or choose a set of $d=q$ neighbors that does not have the symmetry of the data. Decisions on the choice of neighborhood for a given $q$ and $p$ must also be made near the boundaries of the graph, where there is generally less symmetry in the data. It is also an open question, related to Gauss' circle problem \cite{Lowry-Duda2011,Ahmed2018}, of the total number of nearest neighbors that are at most $l$ nearest-neighbors away.

\begin{figure}[H]
	\centering
	\begin{subfigure}[t]{0.33\textwidth}
	\centering
	\includegraphics[width=1\textwidth]{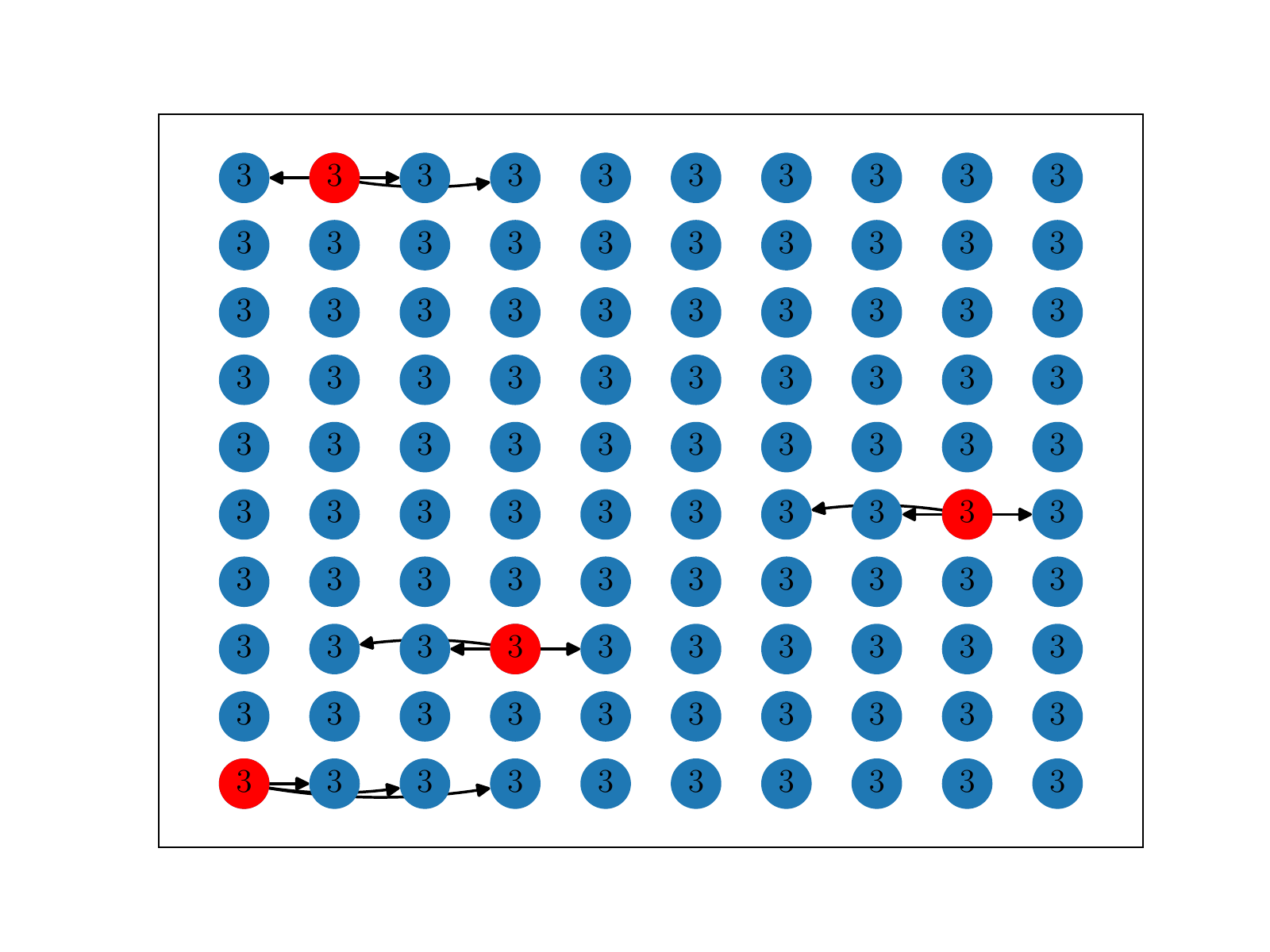}
	\subcaption{Uniform data with constant integer spacing.}
	\label{fig:stencil_n10_p2_uniform}
	\end{subfigure}
	\hspace{1cm}
	\begin{subfigure}[t]{0.33\textwidth}
	\centering
	\includegraphics[width=1\textwidth]{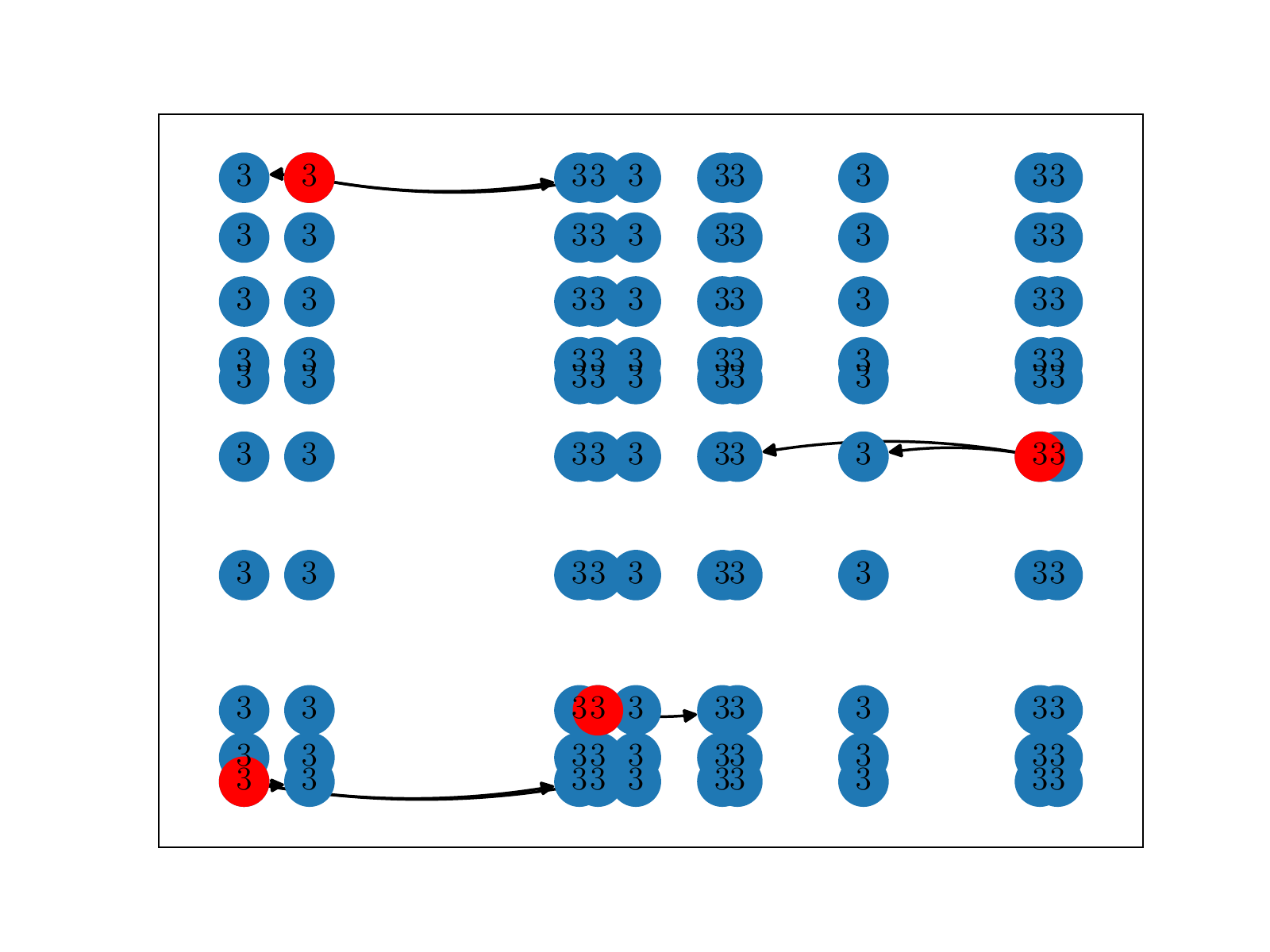}
	\subcaption{Uniform data with non-constant integer spacing.}
	\label{fig:stencil_n10_p2_randint}
	\end{subfigure}
	\hspace{1cm}	
	\begin{subfigure}[t]{0.33\textwidth}
	\centering
	\includegraphics[width=1\textwidth]{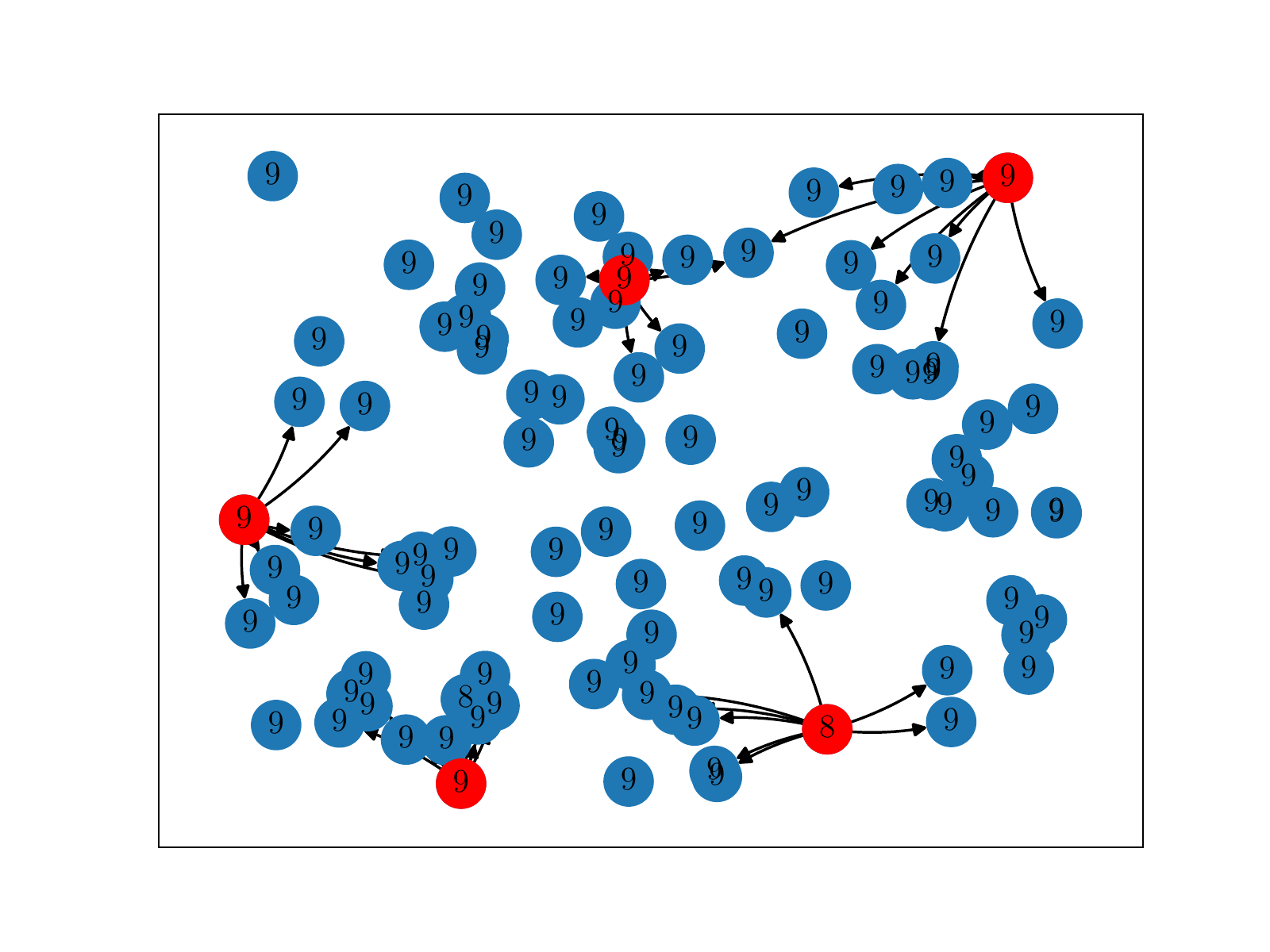}
	\subcaption{Unstructured data with random spacing.}
	\label{fig:stencil_n10_p2_rand}
	\end{subfigure}
	\caption{Example stencils found for $n = 10$ points per dimension in $p = 2$ dimensions. $r = 3$ order accurate stencils for derivatives along the horizontal direction are found, and thus $q = 9$ constraints must be satisfied. Example points are shown in red with arrows to non-zero weights. The number of non-zero weights for the stencil at each point is shown for each vertex.}
	\label{fig:stencil_n10_p2}
\end{figure}

\noindent We also may observe which points in a stencil have non-zero weights, in the case of symmetric and unstructured data, as shown in \cref{fig:stencil_n10_p2}. Here, we see the stencils are found to be the simplest and most intuitive solution when the data is structured with adequate number of points along the dimension of the derivative. Here, only data points along that dimension are assigned non-zero weights, whereas for unstructured data when data along that dimension do not exist, data points along other directions must have non-zero weights.

\subsection{Local weight definitions}\label{app:nonlocalerror}
We will now perform an error analysis of the non-local derivatives with local weight definitions

\subsubsection{Local weight definitions for first order derivatives in one dimension}\label{app:nonlocalerror_localderivative}
We will start with first order derivatives in one dimension.

\noindent Given a point $\widetilde{x}$ and a neighborhood $\mathcal{N}(\widetilde{x}) \subseteq \widetilde{V}$ that is a subset of the whole graph of vertices $\widetilde{V} \subseteq V$ in $p=1$ dimensions, we may evaluate non-local derivatives of the form
\begin{align}
	\difference[1]{u(\widetilde{x})}{x} =&~ \frac{1}{\abs{\mathcal{N}(\widetilde{x})}} \sum_{x \in \mathcal{N}(\widetilde{x})} \left( u(x) - u(\widetilde{x})\right)(x - \widetilde{x})w(\widetilde{x},x), \label{eq:localderivative_app}
\end{align}
with weight functions having support over the neighborhood $w(\widetilde{x},x)$ and size of neighborhood $\abs{\mathcal{N}(\widetilde{x})}$.

\noindent Guided by stencils that are generated by discretization methods, we define weights to be of the form
\begin{align}
	w(\widetilde{x},x) =&~ \frac{\abs{\mathcal{N}(\widetilde{x})}}{(x - \widetilde{x})^2}a(x-\widetilde{x}) \label{eq:localweightfunction_app}
\end{align}
where $a = a(x - \widetilde{x})$ are (dimensionless) reduced weight functions that are independent of vertex spacing. The derivatives therefore take the form: 
\begin{align}
	\difference[1]{u(\widetilde{x})}{x} =&~ \sum_{x \in \mathcal{N}(\widetilde{x})} \frac{u(x) - u(\widetilde{x})}{(x - \widetilde{x})}a(x-\widetilde{x}). \label{eq:localderivative_a_app}
\end{align}
The reduced weights are strictly functions of the vectors $x - \widetilde{x} \in \reals[p]$, and represent the weight of each point $x$ in the stencil of points around $\widetilde{x}$. For example, a symmetric difference two-point stencil in $p=1$ dimensions would have $a(x-\widetilde{x}) = \frac{1}{2}$. 

\noindent For a vertex $\widetilde{x}$ and its neighborhood $\mathcal{N}(\widetilde{x})$, we seek a certain order of accuracy between the non-local and differential calculus derivatives, and find the $a(x-\widetilde{x})$ that yields this order of accuracy. Depending on the local neighborhoods, as per \cref{fig:localweights_app}, the reduced weights are not necessarily symmetric: $a(x-y) \neq a(y-x)$, $x,y \in \widetilde{V}$. 

\noindent In previous works such as \cite{Gilboa2008,Elmoataz2008}, continuous weights $w(\widetilde{x},x)$ are chosen, that have support over the entire graph. These weights generally are chosen to decay away from $\widetilde{x}$, such as a Gaussian weight $w \sim e^{-\abs{x-\widetilde{x}}^2}$, with the conjecture that faster than polynomial decay will lead to the non-local derivatives converging to the differential derivatives. We have conducted several initial analyses and numerical studies with these continuous weights, and have shown that in fact weights with non-local support $\abs{\mathcal{N}(\widetilde{x})} \approx n$ (the dimension of $G$), lead to a constant error between the derivative definitions, for any form of the weights. Local neighborhoods of weights with strictly finite support, $\abs{\mathcal{N}(\widetilde{x})} \ll n$ are therefore required for rigorous convergence for any distribution of vertices on the graph.

\noindent \textbf{Remark}: As an aside, a modified $k$-order Taylor series, as discussed in \cref{app:error_analysis_model} of the \appsecname, may be constructed using a training dataset as a functional representation for $u(x) \approx u_k(x|\widetilde{x})$, where a different Taylor series expansion is developed around each possible point $\widetilde{x}$. Please refer to \cref{app:error_analysis_model} for a definition of this modified Taylor series, \cref{app:error_analysis_mesh_p} for details concerning the graph of data used. A complete error analysis for this modified Taylor series approach can be referred to in \cref{app:error_analysis_model}, which reveals the expected order of accuracy: a $\ith[k]$ order Taylor series has error scaling as $h^{k+1}$.

\noindent The graph of vertices can therefore be thought of as being potentially directed and multi-edged, where pairs of vertices $x$ and $y$ may have an edge $a(y-x)$ from $x$ to $y$, and an edge $a(x-y)$ from $y$ to $x$, depending on the state vectors. Depending on the state vectors of the graph, such as a uniformly spaced mesh of data, or various boundary conditions on the data, there may be additional symmetries that yield a symmetric, undirected graph where $a(y-x) = a(x-y)$. 

\noindent The analysis will not impose these symmetries and each $a(y-x) ~\forall y \in \mathcal{N}(x)$ will be found to ensure each derivative in the graph has a specified order of accuracy. We will impose $q_l$ constraints on the weights for the $\ith[l]$ derivative and we will denote the resulting scaling of the error in each non-local derivative from its differential counterpart as $r_{\textrm{der}_{l}}$ such that
\begin{align}
	\unindifference[l]{u(\widetilde{x})}{x} =&~ \uninderivative[l]{u(\widetilde{x})}{x} + \mathcal{O}(h^{r_{\textrm{der}_{l}}}).
\end{align}
Depending on the definitions of the derivatives via $\delta{\widetilde{x}}$ and the behavior of the error, there will emerge a relationship found between $l$, $k$, $q_l$ and the final scaling $r_{\textrm{der}_{l}}$.

\begin{figure}[hpt]
\centering
\tikzfig{0.9}{figures/localweights}
\caption{Local neighborhoods of adjacent vertices $x$ and $y$, where the absence of symmetry in the weight functions for the edge weight between $x$ and $y$ is shown.}
\label{fig:localweights_app}
\end{figure}

\noindent We may expand the differences of functions in a Taylor series about $\widetilde{x}$:
\begin{align}
	\difference[1]{u(\widetilde{x})}{x} =&~ \sum_{s=0}^{\infty}\frac{1}{(s+1)!}\uninderivative[s+1]{u(\widetilde{x})}{x} \sum_{\mathcal{N}(\widetilde{x})} (x-\widetilde{x})^{s}a(x - \widetilde{x}).
\end{align}

\noindent We now define the monomials for the $d$ data points (vertices) included in the neighborhood $\mathcal{N}(\widetilde{x})$
\begin{align}
	z = z(\widetilde{x}) =&~ x - \widetilde{x} \in \reals[d \times p].
\end{align}
The non-local derivative Taylor series expansion can be written as a function of $z$ 
\begin{align}
	\difference[1]{u(\widetilde{x})}{x} =&~ \sum_{s=0}^{\infty}\frac{1}{(s+1)!}\uninderivative[s+1]{u(\widetilde{x})}{x} \sum_{\mathcal{N}(\widetilde{x})} z(\widetilde{x})^{s}a(z(\widetilde{x})),
\end{align}
which can be written succinctly as a product of two infinite dimensional matrices
\begin{align}
	\difference[1]{u(\widetilde{x})}{x} =&~ d_{\infty}^T V_{\infty}^T a .
\end{align}
Here, $V_{\infty} = V_{\infty}(z) \in \reals[d \times \infty]$ is the Vandermonde-like matrix of powers of $z$ with elements
\begin{align}
	V_{\infty_l}^T =&~ z^{l-1} \in \reals[d],
\end{align}
$d_{\infty} = d_{\infty}(\widetilde{x}) \in \reals[\infty]$ is the vector of derivatives with elements
\begin{align}
	d_{\infty_l} =&~ \frac{1}{l!}\uninderivative[l]{u(\widetilde{x})}{x} \in \reals[],
\end{align}
and $a = a(z) \in \reals[d]$ is the vector of weights with elements
\begin{align}
	a_{l} = a(z_l).
\end{align}
Here we have used base $1$ indexing to be consistent with previous definitions of these matrices.

\noindent If we desire that the non-local derivatives are $r_{\textrm{der}_1} = r \leq d$ order accurate: 
\begin{align}
	\difference[1]{u(\widetilde{x})}{x} =&~ \derivative[1]{u(\widetilde{x})}{x} + \sum_{s=r}^{\infty}\frac{1}{(s+1)!}\uninderivative[s+1]{u(\widetilde{x})}{x} \sum_{\mathcal{N}(\widetilde{x})} z(\widetilde{x})^{s}a(z(\widetilde{x})), \\
	=&~ \derivative[1]{u(\widetilde{x})}{x} + \mathcal{O}(z^{r}),
\end{align}
then the weights can be found from solving the linear system of equations of the first $r$ moments:
\begin{align}
	\sum_{\mathcal{N}(\widetilde{x})} z(\widetilde{x})^{s}a(z(\widetilde{x})) = \delta_{0s} \quad s = \{0,\dots,r-1\} \label{eq:localweightconstraint_moments_1}.
\end{align}
This can be written as a linear problem, given we partition the matrices as
\begin{align}
	V_{\infty} =&~ \left[\begin{array}{ccc} V & \bar{V} \end{array} \right],
	\intertext{where $V = \left[1~Z\right] \in \reals[d \times r],Z \in \reals[d \times r-1], \bar{V} \in \reals[d \times \infty]$ are Vandermonde matrices,}
	V =&~	\left[ 
		\begin{array}{ccccc} 
		1 & z & z^2 & \cdots & z^{r-1}
		\end{array} 
		\right],\\
	Z =&~	\left[ 
		\begin{array}{cccc} 
		z & z^2 & \cdots & z^{r-1}
		\end{array} 
		\right],\\
	\intertext{and}
	d_{\infty} =&~ \left[\begin{array}{c} d \\ \bar{d} \end{array} \right],
\end{align}
where $d = \left[d_1 ~ \hat{d}^T \right]^T \in \reals[r]$, $\hat{d} \in \reals[r-1]$, and $\bar{d} \in \reals[\infty]$.
The linear problem to be solved is
\begin{align}
	d^TV^Ta =&~ d_1 \\
	\intertext{which can be written as a linear combination of the derivatives}
	d_1(1^Ta - 1) +& \hat{d}^TZ^Ta = 0.
\end{align}
The weights must yield consistent derivatives for all functions $u$, and so the vector of derivatives $d$ is linearly independent and the problem can be simplified to
\begin{align}
	V^Ta =&~ e_1 \label{eq:linear_localweights_app}
\end{align}
where $e_1 \in \reals[r]$ has elements $e_{1_l} = \delta_{1l}$.

\noindent For $d > r$, the problem is under-determined, and the minimum norm of $a$ general inverse solution is
\begin{align}
	a = V(V^TV)^{-1}e_1
\end{align}
and when $d = r$, the ordinary least squares pseudo-inverse solution is
\begin{align}
	a = (VV^T)^{-1}Ve_1.
\end{align}

\noindent We will generally choose $d$ points in the stencil to ensure $r=d$ order convergence of the non-local derivative to the differential derivative. Therefore when $r=d$, $V$ is a square and full rank Vandermonde matrix, with a known inverse. \cite{Turner1966,Pugliese2000}

\noindent Given this linear problem for $a(x-\widetilde{x})$ in \cref{eq:linear_localweights_app}, the weight functions for all expansion points $\widetilde{x}$ for the modified Taylor series models can be found using \cref{eq:localweightfunction_app}. The weights are then used in the original non-local calculus definitions in \cref{eq:localderivative_app}, yielding non-local first derivatives that are $r = r_{\textrm{der}_{1}}$ order accurate. This procedure can be repeated for higher order derivatives, and as will be shown in \cref{app:nonlocalerror_localderivative_higher}, distinct sets of edge weights will be found to ensure that each derivative at each expansion point has the desired order of accuracy.


\subsubsection{Local weight definitions for higher order derivatives in one dimension} \label{app:nonlocalerror_localderivative_higher}
Given the definitions of the first derivatives in \cref{eq:localderivative_app}, we may take several approaches to the form of the non-local $l>1$ higher derivatives. One approach is to define the derivatives recursively, and use the same edge weights $w(\widetilde{x},x)$ and same local neighborhoods $\mathcal{N}(\widetilde{x})$ for every order of derivative. This approach is potentially more intuitive, and maintains one set of weights in the graph. However we will initially not impose this recursive constraint, and observe how the different weights, neighborhoods and orders of accuracy are related for different orders of derivatives. 

In order for a given $\ith[l]$ order derivative to have a desired order of accuracy $r_{\textrm{der}_{l}}$, and requiring that each derivative potentially have an independent order of accuracy from other derivatives, we choose to take the approach of defining different sets of weights 
\begin{align}
	w^{(l)}(\widetilde{x},x) =&~ \frac{\abs{\mathcal{N}^{(l)}(\widetilde{x})}}{(x - \widetilde{x})^2}a^{(l)}(x-\widetilde{x}) \label{eq:localweightfunction_higher}
\end{align}
with different local neighborhoods $\mathcal{N}^{(l)}(\widetilde{x})$ for each $\ith[l]$ order of derivative computed on the graph, and impose that the weights $a^{(l)}$ satisfy the $q_l$ linear constraints
\begin{align}
	\sum_{\mathcal{N}^{(l)}(\widetilde{x})} z(\widetilde{x})^{s}a(z(\widetilde{x})) = \delta_{0s}, \quad s = \{0,\dots,r_{\textrm{der}_{l}}-1\}, ~\forall \widetilde{x}.
	\label{eq:localweightfunction_higher_constraints}
\end{align}
In $p=1$ dimensions, the number of constraints equals the desired order of accuracy and $q_l = r_{\textrm{der}_{l}}$.

\begin{figure}[hpt]
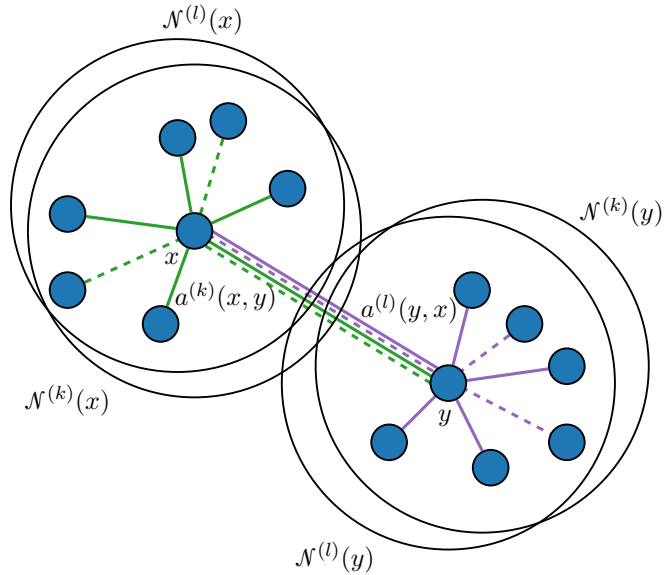

\centering
\tikzfig{0.9}{figures/localweights_higher}
\caption{Local neighborhoods of adjacent vertices $x$ and $y$, for each order neighborhood and edge weights induced by the desired derivative operator. Solid and dashed edges denote different induced graphs for the $k$ and $l$ order derivatives respectively.}
\label{fig:localweights_higher_app}
\end{figure}

\noindent These different weights $w^{(l)}(\widetilde{x},x)$, as well as possibly different neighborhoods $\mathcal{N}^{(l)}(\widetilde{x})$ used for different derivatives can be thought of as being different graph representations of the system, induced by each non-local calculus operator of interest. Illustrations of the induced graphs for the $k$ and $l$ order derivatives, with their distinct neighborhoods and edges, are shown in \cref{fig:localweights_higher_app}.

\noindent The analysis for the higher order derivatives will follow the same procedure as for the first derivatives, where we expand about the vertex of interest the higher order derivatives in Taylor series of the lower order derivatives. The higher-order non-local derivatives are estimated using the following recursive relation obtained by estimating the non-local partial derivative of the function, $\unindifference[l-1]{u(x)}{x}$:
\begin{align}
	\unindifference[l]{u(\widetilde{x})}{x} =&~ \sum_{\mathcal{N}^{(l)}(\widetilde{x})} \left( \unindifference[l-1]{u(x)}{x} - \unindifference[l-1]{u(\widetilde{x})}{x}\right)\frac{a^{(l)}(z(\widetilde{x}))}{z(\widetilde{x})}. \label{eq:localderivative_higher_l_app}
\end{align}
Here, we will denote the monomials $z(x) = y - x$ for $ y \in \mathcal{N}(x)$ with the argument of the vertex $x$, about which the expansion is carried out, and may use primed notation $\mathcal{N}^{\prime}(x) = \{ y^{\prime}~:~ y^{\prime} \in \mathcal{N}(x)\}$ to differentiate between sums over identical neighborhoods. 	

\noindent We will examine the behavior of the $l=2$ non-local derivatives, where given the expansion for the lower non-local first derivative about the base vertex $\widetilde{x}$
\begin{align}
	\difference[1]{u(\widetilde{x})}{x} =&~ \derivative[1]{u(\widetilde{x})}{x} + \sum_{s=r_1}^{\infty}\frac{1}{(s+1)!}\uninderivative[s+1]{u(\widetilde{x})}{x} \sum_{\mathcal{N}^{(1)}(\widetilde{x})} {z(\widetilde{x})}^{s}a^{(1)}(z(\widetilde{x})),
\end{align}
Noting that the above equation for $\difference[1]{u}{x}$ is valid for any vertex, the non-local derivatives at other vertices $x \in \mathcal{N}^{(2)}(\widetilde{x})$ can be expanded about $\widetilde{x}$ as a Taylor series. In the higher order terms, the summations over the neighborhoods are left in terms of $x$:
\begin{align}
	\unidifference[1]{u(x)}{x} =&~ \uniderivative[1]{u(\widetilde{x})}{x} + \uninderivative[2]{u(\widetilde{x})}{x} z(\widetilde{x}) \\
	&+~ \sum_{s=1}^{\infty}\frac{1}{(s+1)!}\uninderivative[s+2]{u(\widetilde{x})}{x} {z(\widetilde{x})}^{s+1} \nonumber \\
	&+~ \sum_{s=r_1}^{\infty}\frac{1}{(s+1)!}\uninderivative[s+1]{u(\widetilde{x})}{x} \sum_{\mathcal{N}^{(1)}(x)} {z(x)}^{s}a^{(1)}(z(x)) \nonumber \\
	&+~ \sum_{s=r_1}^{\infty}\sum_{s^{\prime}=0}^{\infty}\frac{1}{(s+1)!(s^{\prime}+1)!}\uninderivative[s+s^{\prime}+2]{u(\widetilde{x})}{x} {z(\widetilde{x})}^{s^{\prime}+1} \sum_{\mathcal{N}^{(1)}(x)} {z(x)}^{s}a^{(1)}(z(x)). \nonumber
\end{align}
Therefore the non-local second derivative has the expansion
\begin{align}
	\unindifference[2]{u(x)}{x} =&~ \uninderivative[2]{u(\widetilde{x})}{x} \label{eq:localderivative_2} \\
	&+~ \sum_{s=r_2}^{\infty}\frac{1}{(s+1)!}\uninderivative[s+2]{u(\widetilde{x})}{x} \left[\sum_{\mathcal{N}^{(2)}(\widetilde{x})} {z(\widetilde{x})}^{s}a^{(2)}(z(\widetilde{x}))\right] \nonumber \\
	&+~ \sum_{s=r_1}^{\infty}\frac{1}{(s+1)!}\uninderivative[s+1]{u(\widetilde{x})}{x} \left[\sum_{\mathcal{N}^{(2)}(\widetilde{x})} \frac{a^{(2)}(z(\widetilde{x}))}{{z(\widetilde{x})}^{}}\left[\sum_{\mathcal{N}^{(1)^{\prime}}(x)} {z^{\prime}(x)}^{s}a^{(1)}(z^{\prime}(x)) \nonumber \right.\right.\\
	&\quad\quad\quad\quad\quad\quad\quad\quad\quad\quad\quad\quad\quad\quad\quad\quad\quad~~ - \left.\left. \sum_{\mathcal{N}^{^{(1)^{\prime}}} (\widetilde{x})} {z^{\prime}(\widetilde{x})}^{s}a^{(1)}(z^{\prime}(\widetilde{x})) \right]\right] \nonumber \\
	&+~ \sum_{s=r_1}^{\infty}\sum_{s^{\prime}=0}^{\infty}\frac{1}{(s+1)!(s^{\prime}+1)!}\uninderivative[s+s^{\prime}+2]{u(\widetilde{x})}{x} \left[\sum_{\mathcal{N}^{(2)}(\widetilde{x})} {z(\widetilde{x})}^{s^{\prime}}a^{(2)}(z(\widetilde{x})) \left[\sum_{\mathcal{N}^{^{(1)^{\prime}}}(x)} {z^{\prime}(x)}^{s}a^{(1)}(z^{\prime}(x)) \right]\right], \nonumber
\end{align}
where we have substituted the assumed constraints and order of accuracy from the neighborhood for second derivatives in \cref{eq:localweightfunction_higher_constraints}, such that the second summation term has lower bound $r_2$.

\noindent An important observation about the convergence of the higher order $\ith[l]$ derivatives, particularly in looking at the third line in \cref{eq:localderivative_2} as a difference in sums over the vertex neighborhoods $\mathcal{N}(x)$ and $\mathcal{N}(\widetilde{x})$, is that they depend on the homogeneity of the moments of the weights at each vertex. Due to the presence of $z'(x)^{r_{l-1}} / z(\widetilde{x})$ term in the expansion, given an $r_{l-1}$ convergence of the lower derivatives, unless the moments of the weights resulting from constraints of the neighborhood about each vertex are identically equal, there will be a decrease in scaling for each higher derivative and $r_{\textrm{der}_{l}} = r_{\textrm{der}_{l-1}} - 1$ for $l>1$. Therefore we must choose the stencils of the lower order derivatives to have adequate order of accuracy to ensure the highest derivatives still have the desired accuracy $r_{\textrm{der}_{l}}$. One choice is to calculate all weights recursively, and impose that the initial first order derivatives have accuracy $r_{\textrm{der}_{l}} + l$.

\noindent Given this approach to the weights with $q_l = q_1$ number of constraints, the higher order derivatives will scale with decreasing accuracy as
\begin{align}
	\unindifference[l]{u(\widetilde{x})}{x} =&~ \uninderivative[l]{u(\widetilde{x})}{x} + \mathcal{O}(z^{r+1-l}) \\
	\intertext{and so}
	r_{\textrm{der}_{l}} =&~ r + 1 - l \label{eq:derivative_scaling_local_app}.
\end{align}

\noindent Fixing $r \geq k$ constraints on all weights and specifying $d \ll n$ local neighborhoods ensures that all non-local derivatives have $r_{\textrm{der}_{l}}>0$ and converge to their differential derivatives in the limit of an infinite amount of data; i.e., as the number of vertices, $n \to \infty$. Concurrently, in the model error analysis in \cref{eq:localmodelerror}, each $\ith[l] $ derivative error term is a product with $h^{l}$, meaning the error of that term scales as a constant $r_{\textrm{der}_{l}}+l = r + 1$ for all $l$. The decrease in scaling with $l$ is fortuitously cancelled by the scaling with $h$ in the Taylor series model. Having a fixed $r = \bar{r}_{\textrm{der}_{}}$ also means that in the simplest case, there is a single set of $r$ constraints on the moments that all weights $a^{(l)}$ must satisfy in \cref{eq:localweightconstraint_moments_1}, meaning each neighborhood $\mathcal{N}^{(l)}(\widetilde{x}) = \mathcal{N}(\widetilde{x})$ and weights $w^{(l)} = w$ are identical for each order $l$. Higher order derivatives in $p=1$ dimensions therefore have the form
\begin{align}
	\unindifference[l]{u(\widetilde{x})}{x} =&~ \sum_{\mathcal{N}(\widetilde{x})} \left( \unindifference[l-1]{u(x)}{x} - \unindifference[l-1]{u(\widetilde{x})}{x}\right)\frac{a(z(\widetilde{x}))}{z(\widetilde{x})}. \label{eq:localderivative_higher_app}
\end{align}


\subsubsection{Local weight definitions for first order derivatives in higher dimensions}\label{app:nonlocalerror_localderivative_p}

\noindent Guided by stencils underlying discretization methods, we will define the weights along each dimension to be of the symmetric, radial form $w^{\mu}(\widetilde{x},x) = w^{\mu}(x - \widetilde{x})$.
Weights that have a tensor product form built upon one-dimensional weights, $w(\widetilde{x},x) = \prod_{\mu} w^{\mu}(\widetilde{x}^{\mu},x^{\mu})$, lead to weight constraints, particularly for arbitrary vertex neighborhoods, being non-linear and difficult to impose. We therefore define the weights separately for each order of derivative and along each dimension as
\begin{align}
	{w^{\mu}}(\widetilde{x},x) = \frac{\abs{\mathcal{N}^{\mu}(\widetilde{x})}}{(x^{\mu} - \widetilde{x}^{\mu})^{2}}{a^{\mu}}(x-\widetilde{x}), \label{eq:localweightfunction_p_higher}
\end{align}
where the ${a^{\mu}}$ are functions of the whole state vector, and the neighborhoods $\mathcal{N}^{{\mu}}(\widetilde{x})$ are defined for each dimension of derivative.

\noindent The non-local first partial derivatives now take the form 
\begin{align}
	\difference[1]{u(\widetilde{x})}{x^{\mu}} =&~ \sum_{\mathcal{N}^{{\mu}}(\widetilde{x})} \frac{u(x) - u(\widetilde{x})}{x^{\mu} - \widetilde{x}^{\mu}} {a^{\mu}}(x-\widetilde{x}) \label{eq:localderivative_p_app}.
\end{align}
We may expand the differences of functions about $\widetilde{x}$ using a Taylor series as a function of $z$ and so
\begin{align}
	\difference[1]{u(\widetilde{x})}{x^{\mu}} =&~ \sum_{s=0}^{\infty}\sum_{\mu_0\cdots\mu_s}\frac{1}{(s+1)!}\nderivative[s+1]{u(\widetilde{x})}{x^{\mu_0}}{x^{\mu_{s}}} \sum_{\mathcal{N}^{{\mu}}(\widetilde{x})} \frac{z^{\mu_0\cdots\mu_{s}}}{z^{\mu}}{a^{\mu}}(z).
\end{align}
If $r_{\textrm{der}_1^{\mu}} = r_1^{\mu} = r$ order scaling is imposed on the weights, the first derivatives will scale as: 
\begin{align}
	\difference[1]{u(\widetilde{x})}{x^{\mu}} =&~ \derivative[1]{u(\widetilde{x})}{x^{\mu}} + \sum_{s=r_1^{\mu}}^{\infty}\sum_{\mu_0\cdots\mu_s}\frac{1}{(s+1)!}\nderivative[s+1]{u(\widetilde{x})}{x^{\mu_0}}{x^{\mu_{s}}} \sum_{\mathcal{N}^{{\mu}}(\widetilde{x})} \frac{z^{\mu_0\cdots\mu_{s}}}{z^{\mu}}{a^{\mu}}(z) \label{eq:localderivative_p_1} \\
	=&~ \derivative[1]{u(\widetilde{x})}{x^{\mu}} + \mathcal{O}(z^r).
\end{align}
This analysis confirms that for fixed $r$-order accurate stencils, all higher $\ith[l]$-order derivatives in any number of dimensions have error that scales as
\begin{align}
	r_{\textrm{der}_{l}} =&~ r+1-l
\end{align}
and therefore choosing $r \geq k$ will ensure all derivatives have error that decreases with spacing $z$.
Given this non-local derivative error scaling, the error analysis also confirms that a $k$-order Taylor series model comprised of non-local derivatives will have global model of $e = \mathcal{O}(z^{k+1})$, corresponding with Sobolev error analysis for a $k$-order model.

\noindent To find how to define the constraints leading to $r$-order accurate weights, we observe that expanding out our definitions of the non-local derivatives yield constraints on the first $q$ multidimensional moments of the weight distribution:
\begin{align}
	\sum_{\mathcal{N}^{{\mu}}(\widetilde{x})} \frac{z^{\mu_0\cdots\mu_{s}}}{z^{\mu}}{a^{\mu}}(z) =&~ \delta_{0s}\delta^{\mu_0\mu} \cdots \delta^{\mu_{s}\mu}
	. \label{eq:localweightconstraint_moments_p}
\end{align}
The resulting linear problem to be solved is
\begin{align}
	V_{\mu}^T{a^{\mu}}^{} =&~ e_{\mu} \label{eq:linear_localweights_p}
\end{align}
where $a^{\mu} \in \reals[d]$ are the unknown weights, $e_{\mu} \in \reals[q]$ has elements $e_{{\mu}_{l}} = \delta_{\mu l}$, and $V_{\mu} \in \reals[{d \times q}]$ is a Veronese map, represented by a Vandermonde-like block matrix:
\begin{align}
	V_{\mu} =&~	\left[ 
		\begin{array}{cccc} 
		[\frac{z^{\nu}}{z^{\mu}}] & [\frac{z^{\nu\eta}}{z^{\mu}}] & \cdots & [\frac{z^{\mu_{0}\cdots\mu_{r-1}}}{z^{\mu}}]
		\end{array} 
		\right].
\end{align}	
For the $\ith{k}$ block, when considering the non-unique case, has $p^{k}$ elements, and when considering the unique case, has $\binom{p+k-1}{k}$ elements.

\noindent For this Vandermonde-like block matrix, we will set the size of the neighborhoods to be $d = q$ to ensure the linear problem is of full rank and the matrix is potentially invertible. 

\noindent For example, in $p=2$ dimensions and $r=3$, we have coordinate indices $\mu=1,\nu=2$ and the Vandermonde-like block matrix looks like, $z = x - \widetilde{x} \in \reals[d]$:
\begin{align}
	V_{1} =&~	\left[ 
		\begin{array}{ccccccccc} 
		1 & \frac{z^{2}}{z^{1}} & z^{1} & z^{2} & \frac{z^{2^2}}{z^{1}} & z^{1^2} & z^{1}z^{2} & z^{2^2} & \frac{z^{2^3}}{z^{1}}
		\end{array} 
		\right].
\end{align}	

\noindent If $V_{\mu}$ is exactly a Vandermonde matrix, with no denominators, or if the neighborhoods are constrained such that there are no neighboring vertices that are aligned with the base point along the dimension $\mu$ such that $z^{\mu} = 0$, then this matrix has a known inverse. The elements of this matrix have denominators of $z^{\mu}$, and numerators with non-negative powers of $z^{\mu}$, and when the numerator power of $z^{\mu}$ is greater than 0, the divergence as $z^{\mu} \to 0$ will be cancelled. 

\noindent Given the $q(p,r) = \bar{q}(p,r)-1$ unique constraints, where $\bar{q}(p,r) = {(p+r)!}/{(p!~r!)}$, there are $q^{\prime}(p,r) = \bar{q}(p,r-1)$ constraints where $V^{\prime}_{\mu}$ = $\left[1~[z^{\nu}]~\cdots~[z^{\mu_{0}}\cdots z^{\mu_{r-2}}] \right] \in \reals[q^{\prime} \times d]$ has elements with numerators with powers of $z^{\mu}$ greater than $0$. There are similarly $q^{\prime\prime}(p,r) = \bar{q}(p-1,r) - 1$ constraints where $V^{\prime\prime}_{\mu} \in \reals[q^{\prime\prime} \times d]$ has elements with numerators with powers of $z^{\mu}$ equal to $0$. Therefore, the problem can be posed as the sum of regular and possibly divergent contributions, where the total number of constraints is 
\begin{align}
	q =&~ q^{\prime} + q^{\prime\prime},\\
	\intertext{and the linear problem is}
	V_{\mu}^Ta^{\mu} =&~ V_{\mu}^{\prime T}a^{\mu} + V_{\mu}^{\prime\prime T}a^{\mu} = e_{\mu}.
\end{align}

\noindent Care must be taken to ensure that $a$ is identically $0$ at the elements in $V_{\mu}^{\prime\prime}$ that diverge for all $\{\mu_s\}$ when $z^{\mu} \to 0$. This means the generalized invertability of $V_{\mu}$ depends on the the local neighborhoods and the global distribution of data. All $q$ constraints are required, so this can be resolved numerically by choosing different neighborhoods 
\begin{align}
	\mathcal{N}^{{\mu}}(\widetilde{x}) =&~ \{ x \in \mathcal{N}(\widetilde{x}) ~:~ x^{\mu}-\widetilde{x}^{\mu} \neq 0,~ \textrm{rank}(V_{\mu}) = q\}
\end{align}
for each dimension $\mu$ and base point $\widetilde{x}$ such that only neighboring points along the dimension are chosen and there are no divergent constraints. This will also mean additional neighbors along each dimension will have to be found, to ensure that $d^{{\mu}}(\widetilde{x}) = \abs{\mathcal{N}^{{\mu}}(\widetilde{x})} \geq q$ and the linear problem of weight constraints is still well posed. 

The points must also be selected in such a way that the polynomial basis, represented by the matrix $V$ is full rank and non-singular. The neighborhood will be selected by first sorting points using some defined metric, for example, minimum euclidean distance $\norm{z}$, and conditionally adding points to the neighborhood if they do not affect the invertability of $V$. Points are added until $d^{{\mu}}(\widetilde{x}) = q$ and the constrained problem can be solved.


\subsection{Local weight definitions for higher order derivatives in higher dimensions} \label{app:nonlocalerror_localderivative_p_higher}
Given the definitions of the first derivatives in \cref{eq:localderivative_p_app}, and the recursive definition of higher order derivatives
\begin{align}
	\ndifference[l+1]{u(\widetilde{x})}{x^{\mu_0}}{x^{\mu_{l-1}}\delta x^{\mu}} =&~ \sum_{\mathcal{N}^{{\mu}}(\widetilde{x})} \frac{\ndifference[{l}]{u(x)}{x^{\mu_0}}{x^{\mu_{l-1}}} - \ndifference[{l}]{u(\widetilde{x})}{x^{\mu_0}}{x^{\mu_{l-1}}}}{x^{\mu} - \widetilde{x}^{\mu}} {a^{\mu}}(x-\widetilde{x}), \label{eq:localderivative_p_higher},
\end{align}
and the same constraints in \cref{eq:linear_localweights_p} may be used to solve for $a^{{\mu}}$.

\noindent We will examine the behavior of the $l=2,\mu,\nu$ non-local derivatives as an example. Given the expansion for the lower non-local first derivative about the base vertex $\widetilde{x}$
\begin{align}
	\difference[1]{u(\widetilde{x})}{x^{\mu}} =&~ \derivative[1]{u(\widetilde{x})}{x^{\mu}} + \sum_{s=r_{1}^{\mu}}^{\infty}\sum_{\mu_0\cdots\mu_s}\frac{1}{(s+1)!}\nderivative[s+1]{u(\widetilde{x})}{x^{\mu_0}}{x^{\mu_{s}}} \sum_{\mathcal{N}^{{\mu}}(\widetilde{x})} \frac{z(\widetilde{x})^{\mu_0\cdots\mu_{s}}}{z(\widetilde{x})^{\mu}}{a^{\mu}}(z(\widetilde{x}))
\end{align}
the non-local derivatives at other vertices $x \in \mathcal{N}(\widetilde{x})$ can be expanded about $\widetilde{x}$ as
\begin{align}
	\difference[1]{u(x)}{x^{\mu}} =&~ \uniderivative[1]{u(\widetilde{x})}{x^{\mu}} + \derivative[2]{u(\widetilde{x})}{x^{\mu},x^{\nu}} z(\widetilde{x})^{\nu} \\
	&+~ \sum_{s=0}^{\infty}\sum_{\substack{\mu_{0}\cdots\mu_{s}\\\mu_{0}\cdots\mu_{s} \neq \nu ~\textrm{for}~ s = 0}}\frac{1}{(s+1)!}\nderivative[s+2]{u(\widetilde{x})}{x^{\mu}\partial x^{\mu_{0}}}{x^{^{\mu_{s}}}} {z(\widetilde{x})}^{\mu_{0}\cdots \mu_{s}} \nonumber \\
	&+~ \sum_{s=r_1^{\mu}}^{\infty}\sum_{\mu_{0}\cdots\mu_{s}}\frac{1}{(s+1)!}\nderivative[s+1]{u(\widetilde{x})}{x^{\mu_{0}}}{x^{\mu_{s}}} \sum_{\mathcal{N}^{{\mu}}(x)} \frac{z(x)^{\mu_0\cdots\mu_{s}}}{z(x)^{\mu}}{a^{\mu}}(z(x)) \nonumber \\
	&+~ \sum_{s=r_1^{\mu}}^{\infty}\sum_{s^{\prime}=0}^{\infty}\sum_{\mu_{0}\cdots\mu_{s}}\sum_{\mu_{0}^{\prime}\cdots\mu_{s^{\prime}}^{\prime}}\frac{1}{(s+1)!(s^{\prime}+1)!}\nderivative[s+s^{\prime}+2]{u(\widetilde{x})}{x^{\mu_{0}}}{x^{\mu_{s^{\prime}}^{\prime}}} {z(\widetilde{x})}^{\mu_{0}^{\prime}\cdots\mu_{s^{\prime}}^{\prime}} \sum_{\mathcal{N}^{{\mu}}(x)} \frac{z(x)^{\mu_0\cdots\mu_{s}}}{z(x)^{\mu}}{a^{\mu}}(z(x)). \nonumber
\end{align}
Therefore the non-local second derivative has the expansion
\begin{align}
	\difference[2]{u(x)}{x^{\mu},x^{\nu}} =&~ \derivative[2]{u(\widetilde{x})}{x^{\mu},x^{\nu}} \label{eq:localderivative_p_2} \\
	&+~ \sum_{s=r_2^{\nu}}^{\infty}\sum_{\mu_{0}\cdots\mu_{s}}\frac{1}{(s+1)!}\nderivative[s+2]{u(\widetilde{x})}{x^{\mu}\partial x^{\mu_{0}}}{x^{^{\mu_{s}}}}\left[\sum_{\mathcal{N}^{{\nu}}(\widetilde{x})} \frac{z(x)^{\mu_0\cdots\mu_{s}}}{z(x)^{\nu}}{a^{\nu}}(z(x))\right] \nonumber \\
	&+~ \sum_{s=r_1^{\mu}}^{\infty}\sum_{\mu_{0}\cdots\mu_{s}}\frac{1}{(s+1)!}\nderivative[s+1]{u(\widetilde{x})}{x^{\mu_{0}}}{x^{\mu_{s}}} \left[\sum_{\mathcal{N}^{{\nu}}(\widetilde{x})} \frac{{a^{\nu}}(z(\widetilde{x}))}{{z(\widetilde{x})}^{\nu}}\left[ 	\sum_{\mathcal{N}^{{\mu}}(x)} \frac{z^{\prime}(x)^{\mu_0\cdots\mu_{s}}}{z^{\prime}(x)^{\mu}}{a^{\mu}}(z^{\prime}(x))\nonumber \right.\right.\\
	&\quad\quad\quad\quad\quad\quad\quad\quad\quad\quad\quad\quad\quad\quad\quad\quad\quad\quad\quad\quad\quad\quad~ - \left.\left. \sum_{\mathcal{N}^{{\mu}}(\widetilde{x})} \frac{z^{\prime}(\widetilde{x})^{\mu_0\cdots\mu_{s}}}{z^{\prime}(\widetilde{x})^{\mu}}{a^{\mu}}(z^{\prime}(\widetilde{x})) \right]\right] \nonumber \\
	&+~ \sum_{s=r_1^{\mu}}^{\infty}\sum_{s^{\prime}=0}^{\infty}\sum_{\mu_{0}\cdots\mu_{s}}\sum_{\mu_{0}^{\prime}\cdots\mu_{s^{\prime}}^{\prime}}\frac{1}{(s+1)!(s^{\prime}+1)!}\nderivative[s+s^{\prime}+2]{u(\widetilde{x})}{x^{\mu_{0}}}{x^{\mu_{s^{\prime}}^{\prime}}} \Bigg[ \nonumber \\
	&\quad\quad\quad\quad\quad\quad\quad\quad\quad\quad\quad\quad \left. \sum_{\mathcal{N}^{{\nu}}(\widetilde{x})} \frac{{z(\widetilde{x})}^{\mu_{0}^{\prime}\cdots\mu_{s^{\prime}}^{\prime}}}{z(\widetilde{x})^{\nu}} {a^{\nu}}(z(\widetilde{x})) \left[\sum_{\mathcal{N}^{{\mu}}(x)} \frac{z(x)^{\mu_0\cdots\mu_{s}}}{z(x)^{\mu}}{a^{\mu}}(z(x)) \right]\right], \nonumber
\end{align}
where we have substituted the assumed constraints such that the terms from the second derivatives have $r_2^{\nu}$ order scaling.

\noindent Identically to in the one dimensional case, the scaling of the higher order derivatives depends on the homogeneity of the moments of the weights at each vertex. Therefore the same decrease in scaling when no particular symmetries are present for higher order derivatives is $r_{\textrm{der}_{l^{\mu_{0}\cdots\mu_{l}}}} = r_{\textrm{der}_{l^{\mu_{0}\cdots\mu_{l-1}}}} - 1$ for $l>1$. Given this approach to the weights with fixed $r_{l}^{\mu_{0}\cdots\mu_{l-1}} = r = \bar{r}_{\textrm{der}_{}}$ for all orders of derivatives, the higher order derivatives will scale with decreasing accuracy as
\begin{align}
	\ndifference[l]{u(\widetilde{x})}{x^{\mu_{0}}}{x^{\mu_{l-1}}} =&~ \nderivative[l]{u(\widetilde{x})}{x^{\mu_{0}}}{x^{\mu_{l-1}}} + \mathcal{O}(z^{r+1-l}) \\
	\intertext{and so}
	r_{\textrm{der}_{l^{\mu_{0}\cdots\mu_{l-1}}}} =&~ r + 1 - l.
\end{align}
Fixing $r \geq k$ constraints on all weights and specifying $d \ll n$ local neighborhoods ensures that all non-local derivatives have $r_{\textrm{der}_{l}}>0$ and converge to their differential derivatives in the limit of an infinite amount of data.


\subsection{Linear coefficients in one dimension} \label{app:error_analysis_coef}
We can investigate the form of the fit linear coefficients when using the modified $k$ order Taylor series basis in $p=1$ dimensions. We conduct the analysis at a base point $\widetilde{x} \in \widetilde{V}$ and we fix $\gamma_0(\widetilde{x}) = 1$ for model consistency.

\noindent The linear problem to be solved is 
\begin{align}
	y =&~ X\gamma.
\end{align}

\noindent The Taylor series model is dependent on the training data as follows, where the basis of derivatives $\unindifference[l]{u(\widetilde{x})}{x}$ is denoted as $X = X(x|\widetilde{x}) \in \reals[{\widetilde{d} \times k}]$, the linear coefficients are $\gamma = \gamma(\widetilde{x}) \in \reals[k]$, and the known function values are $y = y(x|\widetilde{x}) = u(x) - u(\widetilde{x}) \in \reals[\widetilde{d}]$. Here we are denoting the $\widetilde{d}$ training points $x \subseteq \widetilde{V} \in \reals[\widetilde{d}]$, and the base point $\widetilde{x} \in \reals[]$. Therefore functions of this data are $u(x) \in \reals[\widetilde{d}]$, and $u(\widetilde{x}) \in \reals[]$.

\noindent We may take several approaches to solving this linear system, and choose for this analysis the ordinary least squares solution. The linear coefficients therefore take the form
\begin{align}
	\gamma = (X^TX)^{-1}X^Ty \label{eq:coef_gamma_1d_pseudo}.
\end{align}
Given the Taylor series basis is a polynomial basis, we define monomials in terms of the distances between neighboring vertices to a vertex $\widetilde{x}$
\begin{align}
	z = z(\widetilde{x}) =&~ x-\widetilde{x} \in \reals[\widetilde{d}]
\end{align}
and we will denote generalized dot products, or sums over powers of the components of $z$ as
\begin{align}
	\sum_{\widetilde{\mathcal{N}}(\widetilde{x})} z^l \equiv&~ \sum_{x \in \widetilde{\mathcal{N}}(\widetilde{x})} z(\widetilde{x})^l.
\end{align}

\noindent We can now decompose the basis as
\begin{align}
	X =&~ Z(D+ \Varepsilon).
\end{align}
This decomposition includes a factor of a rectangular Vandermonde-like matrix of $z$ polynomials
\begin{align}
	Z =&~	\left[ 
		\begin{array}{cccc} 
		z & z^2 & \cdots & z^{k}
		\end{array} 
		\right] \in \reals[{\widetilde{d} \times k}],
\end{align}	
and a factor of a diagonal matrix of derivatives that is written as the sum of the differential derivatives, plus error terms
\begin{align}
	D + \Varepsilon =&~ \textrm{diag}(d) + \textrm{diag}(\varepsilon) \in \reals[{k \times k}],
\end{align}
where
\begin{align}
	d =&~ \left[ 
		\begin{array}{cccccccc} 
		\uniderivative[1]{u(\widetilde{x})}{x} & \frac{1}{2!}\uninderivative[2]{u(\widetilde{x})}{x} & \cdots & \frac{1}{k!}\uninderivative[k]{u(\widetilde{x})}{x}
		\end{array} 
		\right] \in \reals[{k}]
	\intertext{and}
	\varepsilon =&~ \left[ 
		\begin{array}{cccccccc} 
		\varepsilon_{1}(\widetilde{x}) & \frac{1}{2!}\varepsilon_{2}(\widetilde{x}) & \cdots & \frac{1}{k!}\varepsilon_{k}(\widetilde{x})
		\end{array} 
		\right] \in \reals[{k}].
\end{align}
The linear coefficients given the ordinary least squares solution are therefore
\begin{align}
	\gamma = (D+\Varepsilon)^{-1}(Z^TZ)^{-1}Z^Ty \label{eq:coef_gamma_1d_korder}.
\end{align}
The diagonal derivative matrix is trivially invertible and so the analytic form of $\gamma$ then can now be derived based on the form of the pseudo-inverse of the Vandermonde-like $Z$. 

For the following analysis, unless specified, we will use base $1$ indexing for the vector and matrix element indices. The objective is to keep the indices consistent with previous analysis for the $\gamma$ coefficients, given we are fixing $\gamma_0$, and $\gamma_{1,\cdots k}$ are being analyzed.

\noindent The Gram matrix of $Z$ has the elements that are sums over the data
\begin{align}
	(Z^TZ)_{ls} =&~ \sum_{\widetilde{\mathcal{N}}(\widetilde{x})} z^{l+s}
\end{align}
and the diagonal derivative matrix, given the errors in the derivatives, has an inverse whose diagonal elements can be written as
\begin{align}
	(D + \Varepsilon)^{-1}_{ll} =&~ \left[\frac{l!}{\uninderivative[l]{u(\widetilde{x})}{x}}\right]\left[1 + \sum_{q=1}^{\infty}\left(-\frac{\varepsilon_{l}(\widetilde{x})}{\uninderivative[l]{u(\widetilde{x})}{x}}\right)^{q}\right].
\end{align}

\begin{figure}[hpt]
	\centering
	\includegraphics[width=0.6\textwidth]{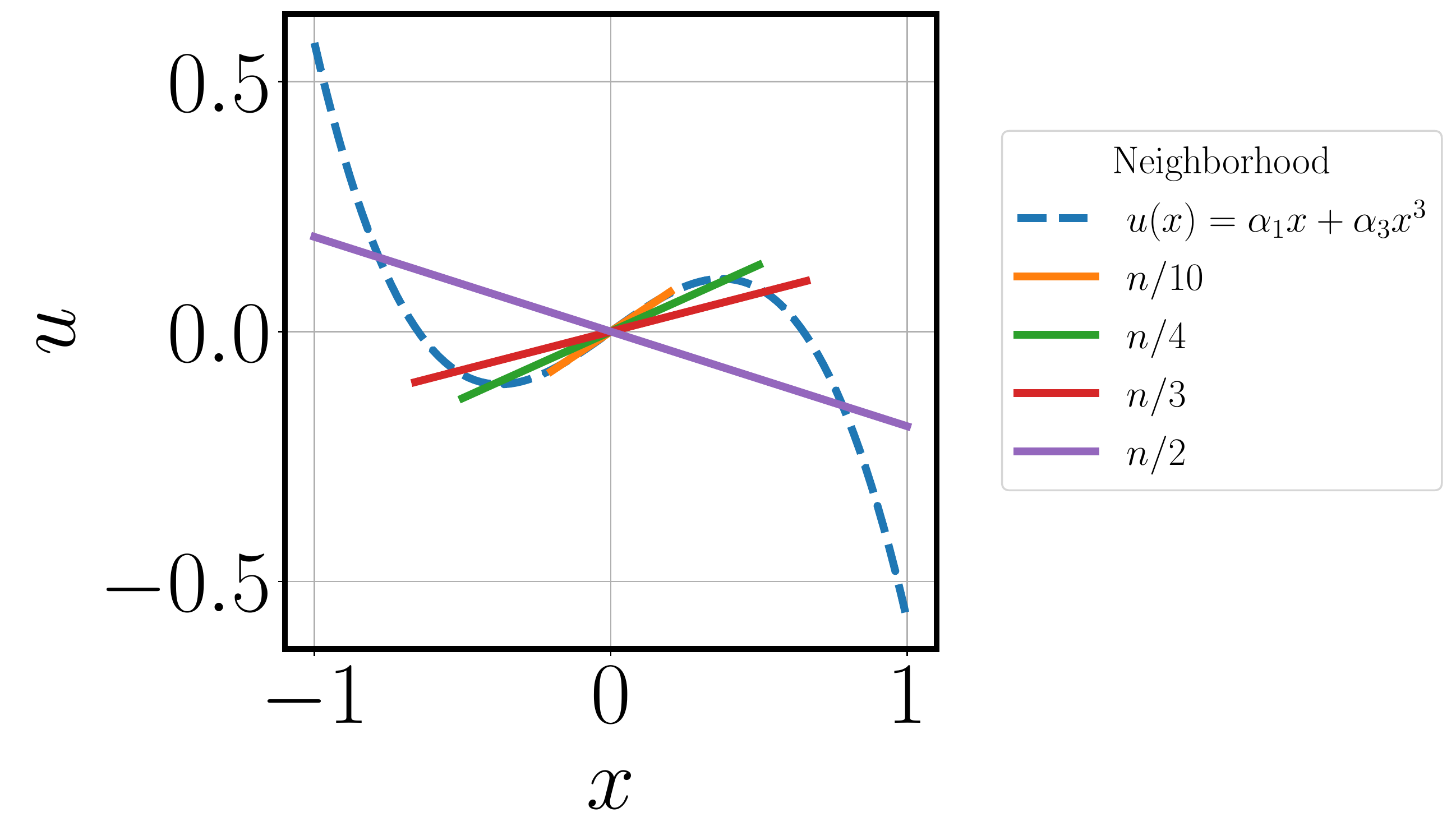}
	\caption{Effect of using neighborhoods with $n/10,n/4,n/3,n/2$ nearest neighbors on the resulting fits with $n$ total data points. Here, a cubic function $u(x|\widetilde{x}) = \alpha_1(x - \widetilde{x}) + \alpha_3(x - \widetilde{x})^3$ is fit with a linear function $u_1(x|\widetilde{x}) = u(\widetilde{x}) + \gamma_1(x-\widetilde{x})$ with ordinary least squares regression at $\widetilde{x} = 0$.}
\label{fig:localfitting}
\end{figure}

\noindent We will now determine the leading behavior of the $\gamma$ coefficients in the limit of an infinite amount of data, for different choices of neighborhoods around the base points $\widetilde{x}$ for the fitting. As can be seen in \cref{fig:localfitting}, if the function to be modeled has adequate changes in curvature and trends far from the base point that cannot be fit by the model, such as the case of a linear model for a cubic function, the size of the neighborhood around the point of fitting greatly affects the resulting linear coefficients and the trends of the fits. In this example, as fitting data becomes more local, the model more adequately fits the local trends of the function.

\noindent The dot products between variables such as in $Z^TZ$ and $Z^Ty$ will be computed as sums over these nearest neighbors, which in this case is evaluated using Faulhaber's formula and we denote as
\begin{align}
	\varphi_{l}[\widetilde{\mathcal{N}}(\widetilde{x})] =&~ \sum_{j \in \widetilde{\mathcal{N}}(\widetilde{x})} j^{l} = \frac{1}{l+1}\widetilde{d}^{l+1} + \frac{1}{2}\widetilde{d}^{l} + O(\widetilde{d}^{l-1}) 
\end{align}
which may scale with $h$ if $\widetilde{d} \sim n$. If the neighborhood is the complete set of integers $\widetilde{\mathcal{N}}(\widetilde{x}) = \{1,\cdots,\widetilde{d}\}$, we denote the argument of $\varphi_{l}[\cdot]$ with the integer $\widetilde{d}$.

\noindent The Gram matrix is in the form of a Hankel-like matrix, which has a closed form for the inverse, however we simply desire the scaling of the $\gamma$ coefficients with $h$, where $\widetilde{\varchi} = \widetilde{j}/n$:
\begin{align}
	(Z^TZ)_{ls} =&~ \sum_{j \in \widetilde{\mathcal{N}}(\widetilde{x})} (j-\widetilde{j})^{l+s} \left(\frac{L}{n}\right)^{l+s} \\
	=&~ \sum_{j \in \widetilde{\mathcal{N}}(\widetilde{x})} \sum_{q=0}^{l+s}\binom{l+s}{q}j^{q}\widetilde{\varchi}^{l+s-q}n^{-q} L^{l+s} \\
	=&~ \sum_{q=0}^{l+s}\binom{l+s}{q}\widetilde{\varchi}^{l+s-q}\varphi_{q}[\widetilde{\mathcal{N}}(\widetilde{x})]n^{-q} L^{l+s}.
\end{align}
We can see that for $\widetilde{d} \sim n$, the sum over the training data in $n^{-q}\varphi_{l}[\widetilde{\mathcal{N}}(\widetilde{x})] \sim \frac{1}{q+1}n$ cancels most of the $h$ scaling of the distance between data points, otherwise for $\widetilde{d} \ll n$, $n^{-q}\varphi_{l}[\widetilde{\mathcal{N}}(\widetilde{x})] \sim n^{-q}$ and the scaling remains. For intermediate $\widetilde{d}$ between $k$ and $n$, the scaling must be investigated numerically.

\noindent We now define the Gram matrix as
\begin{align}
	Z^TZ =&~ W \in \reals[k \times k]
\end{align}
and so the inverse is simply written as
\begin{align}
	(Z^TZ)^{-1} =&~ W^{-1},
\end{align}
which as discussed, may scale more or less with $h^{-1}$ depending on the size of the neighborhood.

\noindent The other factor in the pseudo-inverse calculation of $\gamma$ is $Z^Ty$, which has elements
\begin{align}
	(Z^Ty)_{s} =&~ \sum_{\widetilde{\mathcal{N}}(\widetilde{x})} (u(x)-u(\widetilde{x}))z^{s} \\
	\intertext{which can be expanded as a Taylor series}
	(Z^Ty)_{s} =&~ \sum_{q=1}^{\infty} \frac{1}{q!}\uninderivative[{q}]{u(\widetilde{x})}{x} \sum_{\widetilde{\mathcal{N}}(\widetilde{x})} z^{s+q}.
\end{align}
This matrix is almost identical to the form of $W$ in the Gram matrix, and this term can now be thought of as the product of two infinite dimensional matrices
\begin{align}
	Z^Ty =&~ W_{\infty}d_{\infty}.
\end{align}
The infinite dimensional matrices can be written as partitioned into a $k$ length block, plus the remaining block:
\begin{align}
	W_{\infty} =&~ \left[ \begin{array}{cc} W & \bar{W} \end{array} \right] \in \reals[{k \times \infty}], \\
	\intertext{and}
	d_{\infty} =&~ \left[ \begin{array}{c} d \\ \bar{d} \end{array} \right] \in \reals[{\infty}], \\
	\intertext{where the elements of the matrices are}
	W_{\infty_{sq}} =&~ \sum_{t=0}^{s+q}\binom{s+q}{t} \widetilde{\varchi}^{s+q-t}\varphi_{t}[\widetilde{\mathcal{N}}(\widetilde{x})]n^{-t} L^{s+q}, \\	
	\intertext{and}
	d_{\infty_{q}} =&~ \frac{1}{q!}\uninderivative[{q}]{u(\widetilde{x})}{x}
\end{align}
such that $D^{-1}d = 1 \in \reals[k]$ is the vector of ones.

\noindent Therefore using the form of the Gram matrix, the product of $Z^Ty$, and the inverse of the diagonal derivative matrix, the coefficients take the form
\begin{align}
	\gamma =&~ (I + D^{-1}\Varepsilon)^{-1}D^{-1}W^{-1}(Wd + \bar{W}\bar{d}),
\end{align}
and so the coefficients can be written to leading order, once the resolvant $(I + D^{-1}\Varepsilon)^{-1}$ is expanded, as being unity plus an error term:
\begin{align}
	\gamma =&~ 1 + D^{-1}W^{-1}\bar{W}\bar{d} - D^{-1}\Varepsilon1.
	\label{eq:gamma_in_terms_of_W}
\end{align}

\noindent The scaling of the coefficients we denote as
\begin{align}
	\gamma - 1 =&~ \mathcal{O}(h^{r_{\textrm{coef}}})
\end{align}
which is dependent on the scaling of the coefficient separate from the derivative error, denoted as
\begin{align}
	W^{-1}\bar{W} \sim&~ \mathcal{O}(h^{\bar{r}_{\textrm{coef}}})
	\intertext{and the derivative error, denoted as}
	\Varepsilon \sim&~ \mathcal{O}(h^{r_{\textrm{der}}})
	\intertext{such that}
	r_{\textrm{coef}} =&~ \min{(\bar{r}_{\textrm{coef}},r_{\textrm{der}})}.
\end{align}
Regardless of the error in the derivatives $\gamma - 1$ is proportional to $W^{-1}\bar{W}$ and so the scaling of $W$ with $h$ will determine whether $\lim_{n\to \infty} \gamma \to 1$. 


\subsubsection{Model trained with entire training dataset} \label{app:error_analysis_coef_fulldataset}
If the entire training dataset is used to fit the model at each base point, $\widetilde{\mathcal{N}}(\widetilde{x}) = \widetilde{V}/\{\widetilde{x}\}$ and $\widetilde{d}(\widetilde{x}) = n$. The dot products between variables in the terms of the coefficients will be computed as sums over the entire graph.
For this case, the Gram matrix will tend to the form
\begin{align}
	\lim_{n\to \infty} (Z^TZ) =&~ n\hat{W} + \mathcal{O}(1),
\end{align}
where we denote the matrix $\hat{W}$ with elements
\begin{align}
	\hat{W}_{ls} =&~ \sum_{q=0}^{l+s}\binom{l+s}{q} \frac{1}{q + 1}\widetilde{\varchi}^{l+s-q} L^{l+s} \sim \mathcal{O}(1)
\end{align}
that do not scale with $h$ and so
\begin{align}
	\lim_{n\to \infty} W^{-1} =&~ \frac{1}{n} \hat{W}^{-1} + \mathcal{O}(1).
\end{align}
The product $X^Ty$ will tend to the form
\begin{align}
	\lim_{n\to \infty} (Z^Ty) =&~ n\hat{W}_{\infty}D_{\infty} + \mathcal{O}(1),
\end{align}
where the infinite dimensional $\hat{W}_{\infty} \sim \mathcal{O}(1)$ has the same elements as $\hat{W}$. 

\noindent It can be seen using the entire training dataset to fit the coefficients cancels the $h$ scaling in the leading term of the coefficient. Therefore 
\begin{align}
	\hat{W}^{-1}\bar{\hat{W}} \sim \mathcal{O}(1)
	\intertext{and so}
	r_{\textrm{coef}_l} = 0.
\end{align}
The resulting $\gamma$ for $\widetilde{d} = n$ do not approach $1$ as $n \to \infty$ and a different choice of neighborhood must be made.


\subsubsection{Model trained with nearest neighbors training dataset} \label{app:error_analysis_coef_nearestneighbors}
We also may train the models at each base point using only local data within a $k \leq \widetilde{d} \ll n$ sized neighborhood around the base point such that the linear system is not under-determined. The $W$ matrix elements therefore have the form
\begin{align}
	W_{ls} =&~ \varphi_{l+s}[\widetilde{\mathcal{N}}(\widetilde{x})] \cdot (2h)^{l+s} \label{eq:Gram_nearestneighbor}	
\end{align}
where the sum over the powers of index integers of the data are contained in the $\varphi[\widetilde{\mathcal{N}}(\widetilde{x})]$.

\noindent Given this Hankel-like matrix, inverse elements of a matrix are proportional to the corresponding minors of the matrix, the inverse matrix will also have a similar form and the same reciprocal scaling with $h$. Therefore
\begin{align}
	W_{ls}^{-1} =&~ \varphi^{-1}_{l+s}[\widetilde{\mathcal{N}}(\widetilde{x})] \cdot (2h)^{-(l+s)}.
\end{align}
Here, the $\varphi_{l+s}^{-1}[\widetilde{d}/2]$ matrix, not to be confused with the reciprocal of $\varphi_{l+s}[\widetilde{d}/2]$, do not scale with $h$ for $\widetilde{d} \ll n$.

\noindent Given this form for $W$, and using \cref{eq:gamma_in_terms_of_W}, the coefficients have elements that scale with $h$ like
\begin{align}
	\gamma_{l} =&~ 1 + \sum_{q=1}^{\infty}\frac{l!}{(k+q)!}\frac{\uninderivative[{k+q}]{u(\widetilde{x})}{x}}{\uninderivative[{l}]{u(\widetilde{x})}{x}} \left[\sum_{s=1}^{k}\varphi_{l+s}^{-1}[\widetilde{\mathcal{N}}(\widetilde{x})] \cdot\varphi_{s+k+q}[\widetilde{\mathcal{N}}(\widetilde{x})] \right] (2h)^{k+q-l} \label{eq:coef_scaling_local_neighborhood} \\
	~&~~+~ \frac{\varepsilon_{l}(\widetilde{x})}{\uninderivative[{l}]{u(\widetilde{x})}{x}} \nonumber \\
	=&~ 1 + \mathcal{O}(h^{r_{\textrm{coef}_l}}),
\end{align}
where the product of the matrix and its inverse is contained in the sum over the $\varphi[\widetilde{\mathcal{N}}(\widetilde{x})]$ elements, which are constant with respect to $h$ for this local neighborhood. Therefore by restricting the fitting to being within a $\widetilde{d} \ll n$ sized neighborhood of $\widetilde{x}$, we obtain linear coefficients that approach unity, thus recovering the true differential Taylor series in the limit of infinite data. 

\noindent The minimum scaling of the $\ith[l]$ coefficient with $h$ can now be seen from \cref{eq:coef_scaling_local_neighborhood} to be
\begin{align}
	r_{\textrm{coef}_l} =&~ \min{(k+1-l,r_{\textrm{der}_{l}})} \label{eq:r_coef_scaling}
\end{align}
due to the Taylor series expansion in the $Z^Ty$ product of the pseudoinverse having at least second derivative terms and the scaling of the coefficients independent of the derivative error is $\bar{r}_{\textrm{coef}} = k + q - l$. This $k-l$ dependent means higher order derivatives will have lower order scaling. 

If the data used in fitting also contains symmetry, then there may be additional constraints on the sums over the data in the matrix elements in \cref{eq:Gram_nearestneighbor}, which will in turn affect the scaling bounds in \cref{eq:r_coef_scaling}. For example if there is even (odd) angular symmetry about the $\widetilde{x}$, then $W$ and $W^{-1}$ will have a checkerboard pattern of non-zero elements only at even (odd) indexed elements. This will mean there will only be even powers of $h$ in the scaling of $\gamma$, and there will be a relationship between the order of the model $k$ and the index of the coefficient on its scaling $l$, and $\chi = \mathcal{O}(1)$ for even (odd) symmetry. The coefficients will therefore scale as
\begin{align}
	r_{\textrm{coef}_{l_{\textrm{symmetric}}}} =&~ \min{(k+2 + (k + \chi) ~\textrm{mod}~2 -(l+\chi)~\textrm{mod}~2-l,r_{\textrm{der}_{l}})} \label{eq:r_coef_scaling_symmetric}
\end{align}
meaning models with symmetric neighborhoods have $ k+1-l \leq \bar{r}_{\textrm{coef}_l} \leq k+3-l$, and depending on the choice of $r_{\textrm{der}_{l}}$ some linear coefficients may scale better than for general unstructured neighborhoods. 

\noindent The size of an adequately local neighborhood does not influence the scaling of the coefficient error, however it remains to be seen the relationship between the choice of local neighborhood for the derivatives, and the choice of local neighborhood for the fitting, and how this indirectly influences $\bar{r}_{\textrm{coef}}$ and directly influences $r_{\textrm{der}}$.


\subsection{Linear coefficients in higher dimensions} \label{app:error_analysis_coef_p}
We can investigate the form of the fit linear coefficients when using the modified $k$ order Taylor series basis in $p$ dimensions. This model has $q = (p^{k+1}-1)/(p-1) - 1$ terms if the derivatives do not commute, or $q = (p+k)!/(p!~k!)$ terms if the derivatives do commute, or a full rank polynomial basis is desired. The dataset has a maximum $N = (n+1)^p$ points. We conduct the analysis at a base point $\widetilde{x} \in \widetilde{V}$ and we fix $\gamma_0(\widetilde{x}) = 1$ for model consistency. 

\noindent For the following analysis, unless specified, we will use base $1$ indexing for the vector and matrix element indices. The objective is to keep the indices consistent with previous analysis for the $\gamma$ coefficients, given we are fixing $\gamma_0$, and $\gamma_{1,\cdots k}$ are being analyzed.

\noindent The linear problem to be solved is 
\begin{align}
	y =&~ X\gamma.
\end{align}

\noindent The Taylor series model is dependent on the training data as follows, where the basis of derivatives $\ndifference[l]{u(\widetilde{x})}{x^{\mu_0}}{x^{\mu_{l-1}}}$ is denoted as $X = X(x|\widetilde{x}) \in \reals[{\widetilde{d} \times q}]$, the linear coefficients are $\gamma = \gamma(\widetilde{x}) \in \reals[q]$, and the known function values are $y = y(x|\widetilde{x}) = u(x) - u(\widetilde{x}) \in \reals[\widetilde{d}]$. Here we are denoting the $\widetilde{d}$ training points $x \subseteq \widetilde{V} \in \reals[\widetilde{d} \times p]$, and the base point $\widetilde{x} \in \reals[p]$. Therefore functions of this data are $u(x) \in \reals[\widetilde{d}]$, and $u(\widetilde{x}) \in \reals[]$.

\noindent We may take several approaches to solving this linear system, and choose for this analysis the ordinary least squares solution. The linear coefficients therefore take the form
\begin{align}
	\gamma = (X^TX)^{-1}X^Ty \label{eq:coef_gamma_pd_pseudo}.
\end{align}
Given the Taylor series basis is a polynomial basis, we define the monomials
\begin{align}
	z = z(x|\widetilde{x}) =&~ x-\widetilde{x} \in \reals[\widetilde{d} \times p].
\end{align}

\noindent We also denote block matrices with $k$ blocks
\begin{align}
	A = A(\xi) = \left[ 
		\begin{array}{ccc} 
		A_{1}(\xi) & \cdots & A_{k}(\xi)
		\end{array} 
		\right] \in \reals[m \times q] \label{eq:block_matrix}
\end{align}
where $\xi \in \reals[m \times p]$, and $m$ and $q = \sum_{l=1}^{k}q_l$ are the total dimensions of the matrix. Each block $A_l$ is a rank $k_l$ matrix, where the sub-indices range between $1$ and $p_l$ such that $m$ and $q_l = p_l^{k_l}$ or $q_l = (p_l+k_l-1)!/(p_l-1 !~ k_l!)$ are the dimensions of the block, depending how unique terms are counted. We may also denote these blocks by their components
\begin{align}
	A_l = [a_{\mu_{1}\cdots \mu_{k_l}}] = [a_{11~\cdots ~1}~a_{11~\cdots~2}~\cdots~a_{p_l~p_l~\cdots~p_l}] \in \reals[m \times q_l]
\end{align}
to indicate the form of the elements in each block. The indices of the matrix $A$ are denoted by the block index, with sub-indices for the indices within the block such that
\begin{align}
	A_{l_{\mu_{1}\cdots \mu_{k_l}}} =&~ a_{\mu_1\cdots\mu_{k_l}}.
\end{align}

\noindent We can now decompose the basis as
\begin{align}
	X =&~ Z(D+ \Varepsilon).
\end{align}
This decomposition includes a factor of a rectangular Vandermonde-like matrix of $z$ polynomials
\begin{align}
	Z =&~	\left[ 
		\begin{array}{cccc} 
		[z^{\mu}] & [z^{\mu\nu}] & \cdots & [z^{\mu_{1}\cdots\mu_{k}}]
		\end{array} 
		\right] \in \reals[{\widetilde{d} \times q}],
\end{align}	
and a factor of a diagonal matrix of derivatives that is written as the sum of the differential derivatives, plus error terms
\begin{align}
	D + \Varepsilon =&~ \textrm{diag}(d) + \textrm{diag}(\varepsilon) \in \reals[{k \times k}],
\end{align}
where
\begin{align}
	d =&~ \left[ 
		\begin{array}{cccccccc} 
		[\uniderivative[1]{u(\widetilde{x})}{x^{\mu}}] & [\frac{1}{2!}\derivative[2]{u(\widetilde{x})}{x^{\mu},x^{\nu}}] & \cdots & [\frac{1}{k!}\nderivative[k]{u(\widetilde{x})}{x^{\mu_{1}}}{x^{\mu_{k}}}]
		\end{array} 
		\right] \in \reals[{q}]
	\intertext{and}
	\varepsilon =&~ \left[ 
		\begin{array}{cccccccc} 
		[\varepsilon_{1}^{\mu}(\widetilde{x})] & [\frac{1}{2!}\varepsilon_{2}^{\mu\nu}(\widetilde{x})] & \cdots & [\frac{1}{k!}\varepsilon_{k}^{\mu_{1}\cdots\mu_{k}}]
		\end{array} 
		\right] \in \reals[{q}].
\end{align}
For reference, using the block notation in \cref{eq:block_matrix}, the $Z$ matrix has $k$ blocks $Z_l = [z^{\mu_{1}\cdots\mu_{l}}]$ with rank $k_l = l$ and sub-indices ranging between $1$ and $p_l=p$.

\noindent The linear coefficients given the ordinary least squares solution are therefore
\begin{align}
	\gamma = (D+\Varepsilon)^{-1}(Z^TZ)^{-1}Z^Ty \label{eq:coef_gamma_pd_korder}.
\end{align}
The diagonal derivative matrix is trivially invertible and so the analytic form of $\gamma$ then can now be derived based on the form of the pseudo-inverse of the Vandermonde-like $Z$.

\noindent We now define the Gram matrix as
\begin{align}
	Z^TZ =&~ W \in \reals[q\times q]
\end{align}
and so the inverse is simply
\begin{align}
	Z^TZ^{-1} =&~ W^{-1},
\end{align}
which as discussed, may scale more or less with $h^{-1}$ depending on the size of the neighborhood.

\noindent The other factor in the pseudo-inverse calculation of $\gamma$ is $Z^Ty$, which can now be thought of as the product of two infinite dimensional matrices
\begin{align}
	Z^Ty =&~ W_{\infty}d_{\infty}.
\end{align}
The infinite dimensional matrices can be written as
\begin{align}
	W_{\infty} =&~ \left[ \begin{array}{cc} W & \bar{W} \end{array} \right] \in \reals[{q \times \infty}], \\
	\intertext{and}
	d_{\infty} =&~ \left[ \begin{array}{c} d \\ \bar{d} \end{array} \right] \in \reals[{\infty}].
\end{align}
such that $D^{-1}d = 1 \in \reals[q]$ is the vector of ones.

\noindent Therefore using the form of the Gram matrix, the product of $Z^Ty$, and the inverse of the diagonal derivative matrix, the coefficients take the form
\begin{align}
	\gamma =&~ (I + D^{-1}\Varepsilon)^{-1}D^{-1}W^{-1}(Wd + \bar{W}\bar{d}),
\end{align}
and so the coefficients can be written to leading order, once the resolvant $(I + D^{-1}\Varepsilon)^{-1}$ is expanded, as being unity plus an error term:
\begin{align}
	\gamma =&~ 1 + D^{-1}W^{-1}\bar{W}\bar{d} - D^{-1}\Varepsilon1.
\end{align}

\noindent The scaling of the coefficients we denote as
\begin{align}
	\gamma - 1 =&~ \mathcal{O}(r_{\textrm{coef}})
\end{align}
which is dependent on the scaling of the coefficient separate from the derivative error
\begin{align}
	W^{-1}\bar{W} \sim&~ \mathcal{O}(h^{\bar{r}_{\textrm{coef}}})
	\intertext{and the derivative error, denoted as}
	\Varepsilon \sim&~ \mathcal{O}(h^{r_{\textrm{der}}})
	\intertext{such that}
	r_{\textrm{coef}} =&~ \min{(\bar{r}_{\textrm{coef}},r_{\textrm{der}})}.
\end{align}
Regardless of the error in the derivatives $\gamma - 1$ is proportional to $W^{-1}\bar{W}$ and so the scaling of $W$ with $h$ will determine whether $\lim_{n\to \infty} \gamma \to 1$.

\noindent The dot products between variables such as in $Z^TZ$ and $Z^Ty$ will be computed as sums over these nearest neighbors, which we denote as
\begin{align}
	\varphi_{l_{\mu_{1}^{t_{1}}\cdots\mu_{l}^{t_{l}}}}[\widetilde{\mathcal{N}}(\widetilde{x})] =&~ \sum_{j \in \widetilde{\mathcal{N}}(\widetilde{x})} j^{\mu_{1}^{t_{1}}\cdots\mu_{l}^{t_{l}}}
\end{align}
for particular indices $\{\mu_{q}\}$ and powers $\{t_q\}$. This sum may scale with $h$ if $\widetilde{d} \sim n$.

\noindent The Gram matrix is in the form of a Hankel-like matrix, which has a closed form for the inverse \cite{Pugliese2000}, however we simply desire the scaling of the $\gamma$ coefficients with $h$:
\begin{align}
	W_{l_{\mu_1\cdots\mu_{l}}s_{\nu_1\cdots\nu_{s}}} =&~ \sum_{j \in \widetilde{\mathcal{N}}(\widetilde{x})} \left[\prod_{t=1}^{l}(j^{\mu_{t}}-\widetilde{j}^{\mu_{t}})\right]\left[\prod_{t=1}^{s}(j^{\nu_{t}}-\widetilde{j}^{\nu_{t}})\right] \left(\frac{L}{n}\right)^{l+s} \\
	=&~ \sum_{j \in \widetilde{\mathcal{N}}(\widetilde{x})} \sum_{\substack{\substack{t_{l1}\cdots t_{ll}\\t_{s1}\cdots t_{ss}}=0\\t = \sum_{q=1}^{l} t_{l_q} + \sum_{q=1}^{s} t_{s_q}}}^{1}j^{\mu_{1}^{t_{s1}}\cdots\nu_{l}^{t_{ll}}}\widetilde{\varchi}^{\mu_{1}^{1-t_{s1}}\cdots\nu_{l}^{1-t_{ll}}} n^{-t} L^{l+s} \\
	=&~ \sum_{\substack{\substack{t_{l1}\cdots t_{ll}\\t_{s1}\cdots t_{ss}}=0\\t = \sum_{q=1}^{l} t_{l_q} + \sum_{q=1}^{s} t_{s_q}}}^{1} \varphi_{{l+s}_{\mu_{1}^{t_{l1}}\cdots\nu_{s}^{t_{ss}}}}[\widetilde{\mathcal{N}}(\widetilde{x})] \widetilde{\varchi}^{\mu_{1}^{1-t_{s1}}\cdots\nu_{l}^{1-t_{ll}}} n^{-t} L^{l+s} .
\end{align}
We observe that in higher dimensions, the form of the matrix elements and their scaling with $\widetilde{d}$ is much more complicated, however the scaling with respect to $h$ is identical and so the same conclusions as for $r_{\textrm{coef}}$ in \cref{eq:r_coef_scaling} for the one dimensional problem apply.

\noindent For example, if we assume that $\widetilde{d} \ll n$ for a local neighborhood around $\widetilde{x}$, then the Gram matrix elements can be written as
\begin{align}
	W_{l_{\mu_1\cdots\mu_{l}}s_{\nu_1\cdots\nu_{s}}} =&~ \varphi_{{l+s}_{\mu_{1}\cdots\nu_{s}}}[\widetilde{\mathcal{N}}(\widetilde{x})] h^{l+s},
\end{align}
the inverse Gram matrix can be written as
\begin{align}
	W_{l_{\mu_1\cdots\mu_{l}}s_{\nu_1\cdots\nu_{s}}}^{-1} =&~ \varphi_{{l+s}_{\mu_{1}\cdots\nu_{s}}}^{-1}[\widetilde{\mathcal{N}}(\widetilde{x})] h^{-(l+s)},
\end{align}
and identical conclusions about the scaling of the coefficients for local fitting neighbors can be made:
\begin{align}
	r_{\textrm{coef}_{l_{\mu_1\cdots\mu_{l}}}} &=~ \min{(k+1-l,r_{\textrm{der}_{l_{\mu_1\cdots\mu_{l}}}})}.
\end{align}


\subsection{Global model error} \label{app:error_analysis_global}

Given the scaling of the individual components of the error as $r_{\textrm{local}},~r_{\textrm{der}},~\textrm{and}~ r_{\textrm{coef}}$, and defining the local error as 
\begin{align}
	e(x | \widetilde{x}) =&~ \sum_{l=0}^{k}\frac{\gamma_l(\widetilde{x}) - 1}{l!} \uninderivative[l]{u(\widetilde{x})}{x} (x-\widetilde{x})^l \label{eq:localmodelerror}\\
	~+&~ \sum_{l=0}^{k} \frac{\gamma_l(\widetilde{x})}{l!} \varepsilon_l(\widetilde{x}) (x-\widetilde{x})^l \nonumber \\
	~-&~ \sum_{l=k+1}^{K}\frac{1}{l!} \uninderivative[l]{u(\widetilde{x})}{x} (x-\widetilde{x})^{l} \nonumber 
\end{align}
and assuming the model is consistent such that $\gamma_0(\widetilde{x}) = 1$ and $\varepsilon_0(\widetilde{x}) = 0$, we can observe the minimum scaling of the local model error as
\begin{align}
	e(x | \widetilde{x}) \sim&~ \sum_{l=1}^{k}\mathcal{O}(h^{r_{\textrm{coef}_l}+l}) + \sum_{l=1}^{k}(1 + \mathcal{O}(h^{r_{\textrm{coef}_l}}))\mathcal{O}(h^{r_{\textrm{der}_{l}}+l}) + \mathcal{O}(h^{k+1}) \\
	=&~ \mathcal{O}(h^{r_{\textrm{local}}}).
\end{align}
Given the scaling of the derivatives is
\begin{align}
	r_{\textrm{der}_{l}} =&~ r + 1 - l
\end{align}
where ideally $r \geq k$, and given the scaling of the coefficients is
\begin{align}
	r_{\textrm{coef}_l} =&~ \min {(k+1-l,r_{\textrm{der}_{l}})} = \min{(r,k)} + 1 - l
\end{align}
we can see the $\ith[l]$ derivative and coefficient scaling decreases linearly with $l$. Therefore the local model scaling is
\begin{align}
	r_{\textrm{local}}=&~ \min_{1\leq l \leq k}{(r_{\textrm{coef}_l}+l,r_{\textrm{der}_{l}}+l,k+1)} \\
	=&~ \min_{1 \leq l \leq k} \min{(k+1,r_{\textrm{der}_{l}}+l)} \nonumber
\end{align}	
which is seen to be independent of the bare coefficient scaling $\bar{r}_{\textrm{coef}_l}$ and so
\begin{align}
	r_{\textrm{local}} =&~ \min{(r,k)} + 1 \label{eq:local_scaling_bounds}.
\end{align} 

This result is true for both symmetric and unstructured neighborhoods due to having shown that $\bar{r}_{\textrm{coef}} + l \geq k+1$. Therefore for all local neighborhoods, the local error model is limited by the accuracy of either the derivatives, or the model itself relative to the function being modeled. For example, when $r = 2$, such as in a central difference finite difference scheme for the non-local derivatives, for $k>1$ the error is limited by the order of the model. We also note that the form of the Taylor series model means that since the error in the $\ith[l]$ coefficient and derivative terms scale as some constant $r-l$, and are products with terms that scale like $l$, there are fortuitous cancellations between these scalings with $l$, making the most terms in the local model error have a constant scaling. 

The global error when considering the finite set of testing points $x \in V$, given the training points $\widetilde{x} \in \widetilde{V}$ can be written as being proportional to the average in the $l$-norm of local errors:
\begin{align}
	\norm{e}_l^l \sim&~ \frac{1}{n}\sum_{V}\abs{e(x)}^l, \label{eq:totalloss}
\end{align}
and given the form of the local error
\begin{align}
	\norm{e}_l^l \sim&~ \frac{2h}{L}\sum_{j \in V}\abs{C_0(x_j|\widetilde{x}_{\widetilde{j}}) + C_1(x_j|\widetilde{x}_{\widetilde{j}})h + C_2(x_j|\widetilde{x}_{\widetilde{j}}) h^2 + \cdots}^l \\
	\leq&~ \frac{2h}{L}n\abs{C_0 + C_1h + C_2h^2 + \cdots}^l \\
	\sim&~ \abs{C_0 + C_1h + C_2h^2 + \cdots}^{l},
\end{align}
and therefore the global error scales as 
\begin{align}
	\norm{e}_l \sim \mathcal{O}(h^{r_{\textrm{local}}}) \\
	\intertext{and}
	r_{\textrm{global}} = r_{\textrm{local}}.
\end{align}
Given that we have found the scaling for the local error of the modified Taylor series model in \cref{eq:local_scaling_bounds}, the global error is
\begin{align}
	r_{\textrm{global}} =&~ \min{(r,k)} + 1\label{eq:global_scaling_bounds}.
\end{align}

\subsubsection{Computing global error}
To derive the form of the average error in \cref{eq:totalloss}, we approximate the total $l$-norm error $\norm{e}_l$ across the domain of interest. We assume there is a closed form for the local error $e(x|\widetilde{x}) \sim C_q(\widetilde{x})(x-\widetilde{x})^q \sim C_q(x|\widetilde{x})h^q$ that is $q$ order, and the data mesh has spacing $2h$ that splits the domain into $n$ intervals.

We approximate the total interpolation error of the model along the intervals of the entire domain by the error at the point $x_j \in [\widetilde{x}_{\widetilde{j}},\widetilde{x}_{i+1}]$ in the interior of the interval. The resulting form is what is used in \cref{eq:totalloss}:
\begin{align}
	\norm{e}_l^l =&~ \frac{1}{L} \int\limits_{0}^L dx \abs{e_{j}(x)}^l \\
	=&~ \frac{1}{L} \sum_{j=0}^{n-1}\int\limits_{\widetilde{x}_{i}}^{\widetilde{x}_{i+1}} dx~ \abs{e(x|\widetilde{x}_{\widetilde{j}})}^l \nonumber\\
	\approx&~ \frac{1}{L} \sum_{j=0}^{n-1} 2h \abs{e(x_j |\widetilde{x}_{\widetilde{j}})}^l \nonumber\\
	=&~ \frac{1}{n} \sum_{j=0}^{n-1} \abs{e(x_j |\widetilde{x}_{\widetilde{j}})}^l \nonumber \\	
	=&~ \frac{1}{n} \sum_{j=0}^{n-1} \abs{C_q(\widetilde{x}_{\widetilde{j}})}^l h^{lq} \nonumber \\
	\leq&~ \frac{1}{n} \sum_{j=0}^{n-1} C_{q}^l h^{lq} \nonumber \\
	\leq&~ C_{q}^l h^{lq} \nonumber \\
	\intertext{ and therefore the approximated global interpolation error scales identically to the local error:}
	\norm{e}_l \leq&~ \left[C_{q}\right] h^{q} \sim \mathcal{O}(h^{q}).
\end{align}


\subsection{Numerical error of modified Taylor series model in $p$ dimensions}\label{app:nonlocalerror_localderivative_results}
\begin{figure}[hpt]
\centering
\includegraphics[width=\textwidth]{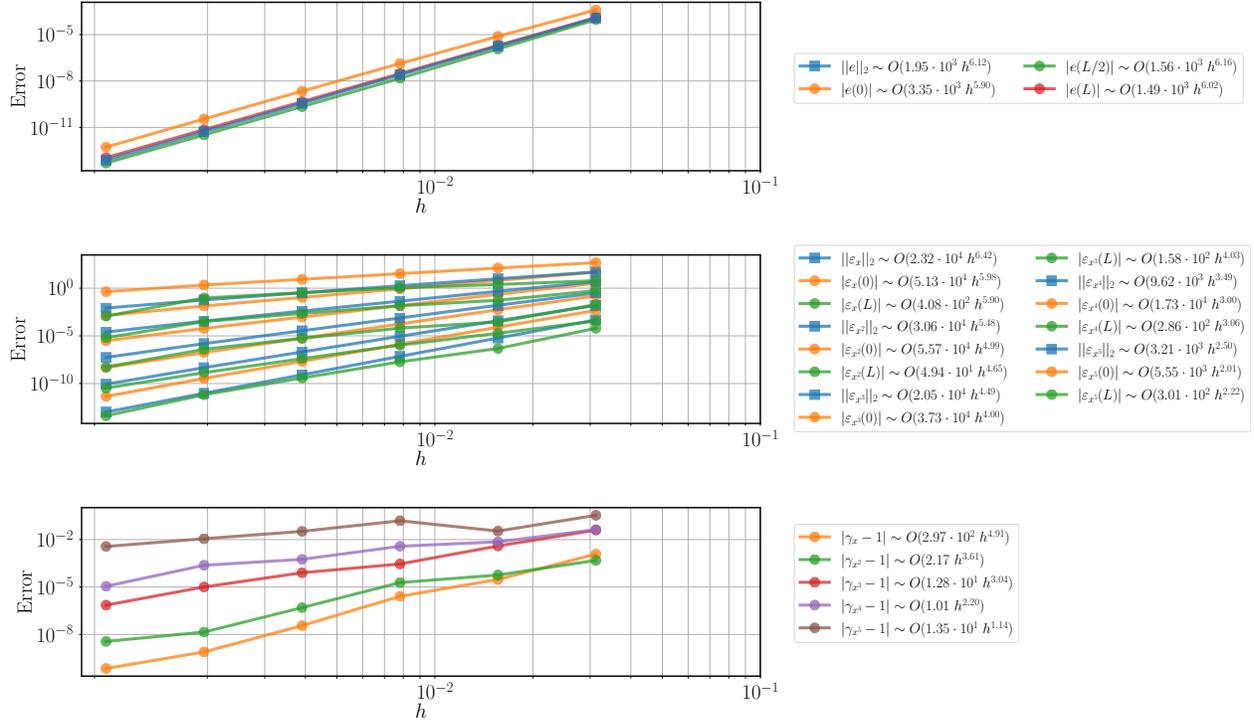}
\caption{Error scaling of components of a $r=k+1$ accurate, $k=5$-order Taylor series model for a $K=8$ order polynomial $u(x) = \sum_{\sum s \leq K}\alpha_{s} x^{s}$ in $p=1$ dimensions. Polynomial coefficients $\alpha_{s} \sim U[-1,-1]$ are sampled from a uniform distribution and all unique commuting Taylor series derivatives are included in the model up to $k$ order terms. Top: Global and local model error $e$. Middle: Global and local non-local derivative error $\varepsilon$. Bottom: Global and local error of $\gamma$. Leading order scaling with $h$ fits are shown in the legend for global and local values.}
\label{fig:errorscaling_1_app}
\end{figure}

We now present numerical results for $p=1,2$ dimensional error scaling of the local and global model error, the non-local derivatives, and the linear coefficients in \cref{fig:errorscaling_1_app,fig:errorscaling_2_app}. We find that the scaling behavior of the linear coefficients is highly dependent on the particular linear solver used, however the exact scaling does not affect the global model error. As predicted by the analysis, the plots confirm that this $k$-order Taylor series model, with $n = \mathcal{O}(1/h)$ local models and $r = k+1$-order accurate local stencils, has non-local derivative error of $\varepsilon_l = \mathcal{O}(h^{r+1-l})$, linear coefficient error of $\gamma_l - 1 = \mathcal{O}(h^{k+1-l})$, and global model error of $e = \mathcal{O}(h^{k+1})$.

\begin{figure}[hpt]
\centering
\includegraphics[width=\textwidth]{figures/scaling_2d.pdf}
\caption{Error scaling of components of a $r=k+1$ accurate, $k=3$-order Taylor series model for a $K=6$ order polynomial $u(x) = \sum_{\sum s \leq K}\alpha_{s} x^{s}$ in $p=2$ dimensions. Polynomial coefficients $\alpha_{s} \sim U[-1,-1]$ are sampled from a uniform distribution and all unique commuting Taylor series derivatives are included in the model up to $k$ order terms. Top: Global and local model error $e$. Middle: Global and local non-local derivative error $\varepsilon$. Bottom: Global and local error of $\gamma$. Leading order scaling with $h$ fits are shown in the legend for global and local values.}
\label{fig:errorscaling_2_app}
\end{figure}


\newpage
\section{Free energy functional representation} \label{sec:physicalsystems_psi}
In addition to modelling the phase evolution of microstructures via the coupled Cahn-Hilliard, gradient elasticity problem, we have applied the non-local calculus to represent the total free energy. In this, we are guided by the polynomial form of the free energy density, and propose $\Psi = \Psi(\bar{\varphi},{\vecbar{E}})$. The model therefore takes the form of a modified Taylor series at time $j$, about a base value at the time $i$:	
\begin{align}
	\Psi(\bar{\varphi}_{j},\vecbar{E}_{j}) = \Psi(\bar{\varphi}_{i},\vecbar{E}_{i})
	~&+~ \gamma^{\bar{\varphi}} \difference[1]{\Psi(\bar{\varphi}_{i},\vecbar{E}_{i})}{\bar{\varphi}} \Delta {{}{\bar{\varphi}}}_{ij} ~+~ \gamma^{\bar{E}_{\alpha \beta}} \difference[1]{\Psi(\bar{\varphi}_{i},\vecbar{E}_{i})}{\bar{E}_{\alpha \beta}} \Delta {{}\bar{E}_{\alpha \beta}}_{ij} \label{eq:model_psitaylor} \\
	~&+~ \gamma^{{\bar{\varphi} \bar{\varphi}}} \frac{1}{2!}\unidifference[2]{\Psi(\bar{\varphi}_{i},\vecbar{E}_{i})}{{\bar{\varphi}}} {\Delta {{}{\bar{\varphi}}}_{ij}}^2 \nonumber\\
	 ~&+~ \gamma^{{\bar{E}_{\alpha \beta}}{\bar{E}_{\xi \kappa}}} \frac{1}{2!}\difference[2]{\Psi(\bar{\varphi}_{i},\vecbar{E}_{i})}{{\bar{E}_{\alpha \beta}},{\bar{E}_{\xi \kappa}}} \Delta {{}{\bar{E}_{\alpha \beta}}}_{ij} \Delta {{}{\bar{E}_{\xi \kappa}}}_{ij} \nonumber\\
	 ~&+~ \gamma^{{\bar{\varphi}}{\bar{E}_{\alpha \beta}}} \difference[2]{\Psi(\bar{\varphi}_{i},\vecbar{E}_{i})}{{\bar{\varphi}},{\bar{E}_{\alpha \beta}}} \Delta {{}{\bar{\varphi}}}_{ij} \Delta {{}{\bar{E}_{\alpha \beta}}}_{ij} \nonumber \\
	 ~&+~ O(\Delta {{}{\bar{\varphi}}}_{ij}^3) + O(\Delta {{}{\bar{E}_{\alpha \beta}}}_{ij} ^3), \nonumber
\end{align}
where ${\Delta x_{ij}} = y_{} - x_{}$ represents the change in the $\varphi_{\alpha}$ or $\vecbar{E}$ state variable between times $j$ and $i$. 

\begin{figure}[hpt]
\centering
\includegraphics[width=0.5\textwidth]{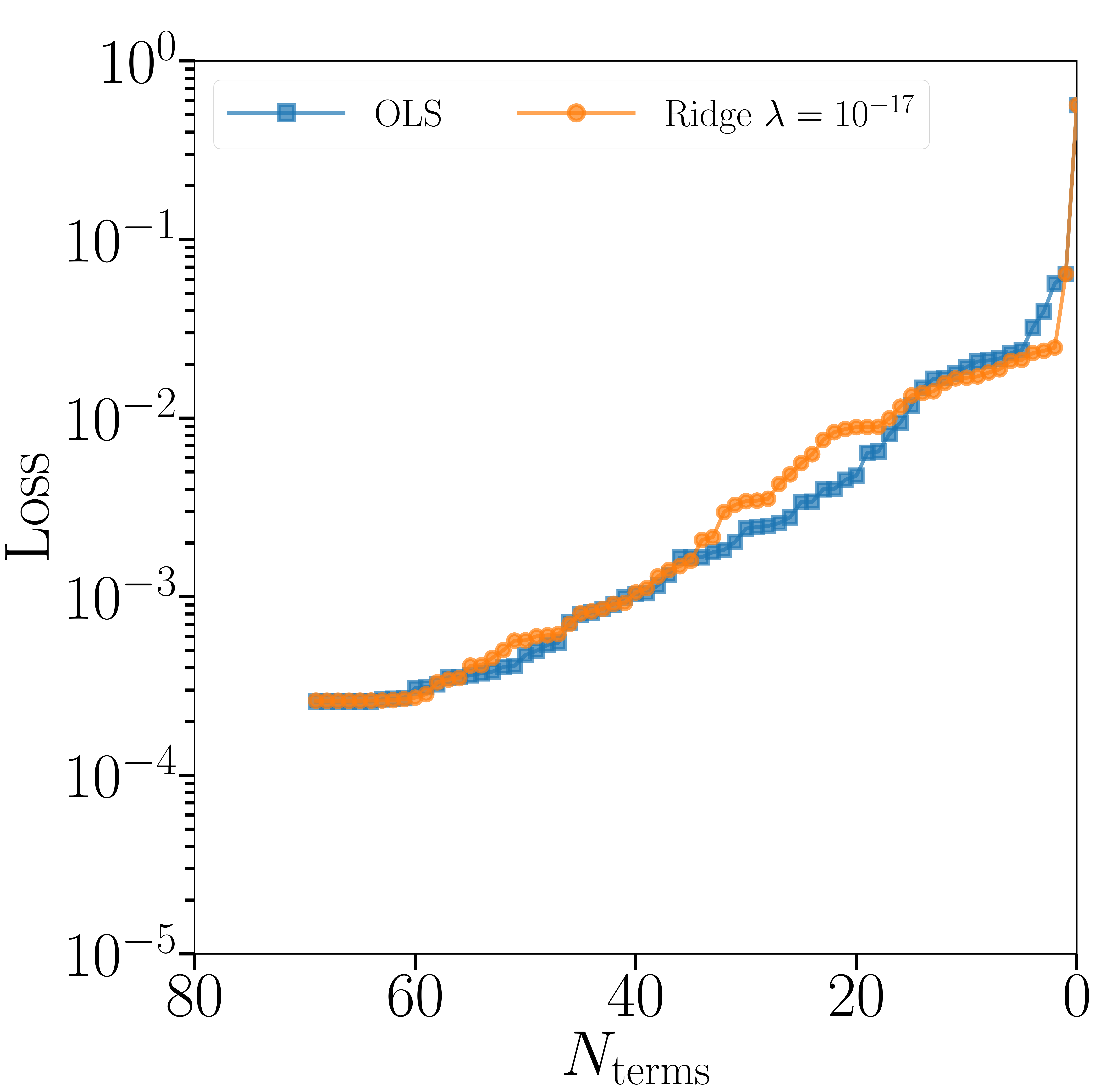}
\caption{Stepwise regression loss curve for the $\ith[4]$ order free energy Taylor series functional representation, using OLS and Ridge regression.}
\label{fig:psi_loss_OLS_Ridge17_local}
\end{figure}

\noindent Stepwise regression is performed with a $\ith[4]$-order Taylor series, using an $l_2$ loss function for this roughly monotonic function, with OLS and Ridge regression, where the optimal ridge parameter is found to be small $\lambda = 10^{-17}$. The first 10 terms from the stepwise regression for the $\ith[4]$ order Taylor series functional representation for $\Psi$ are found for OLS and Ridge regression to be:
\input{./figures/model_psi_OLS_Ridge17_local.tex}
\noindent The OLS and Ridge regression loss curves in \cref{fig:psi_loss_OLS_Ridge17_local}, unlike for the phase volume fraction dynamics, do not indicate a plateau where there is a clear distinction of the most relevant terms in the model. The loss increases smoothly over the regression iterations, likely due to the accuracy of Taylor series increasing monotonically with the order of the model. The full model has loss of order $\mathcal{O}(10^{-4})$, comparable to other methods. \cite{Zhang2020}

\begin{figure}[hpt]
\centering
\begin{subfigure}[t]{0.49\textwidth}
	\centering
	\includegraphics[width=\textwidth]{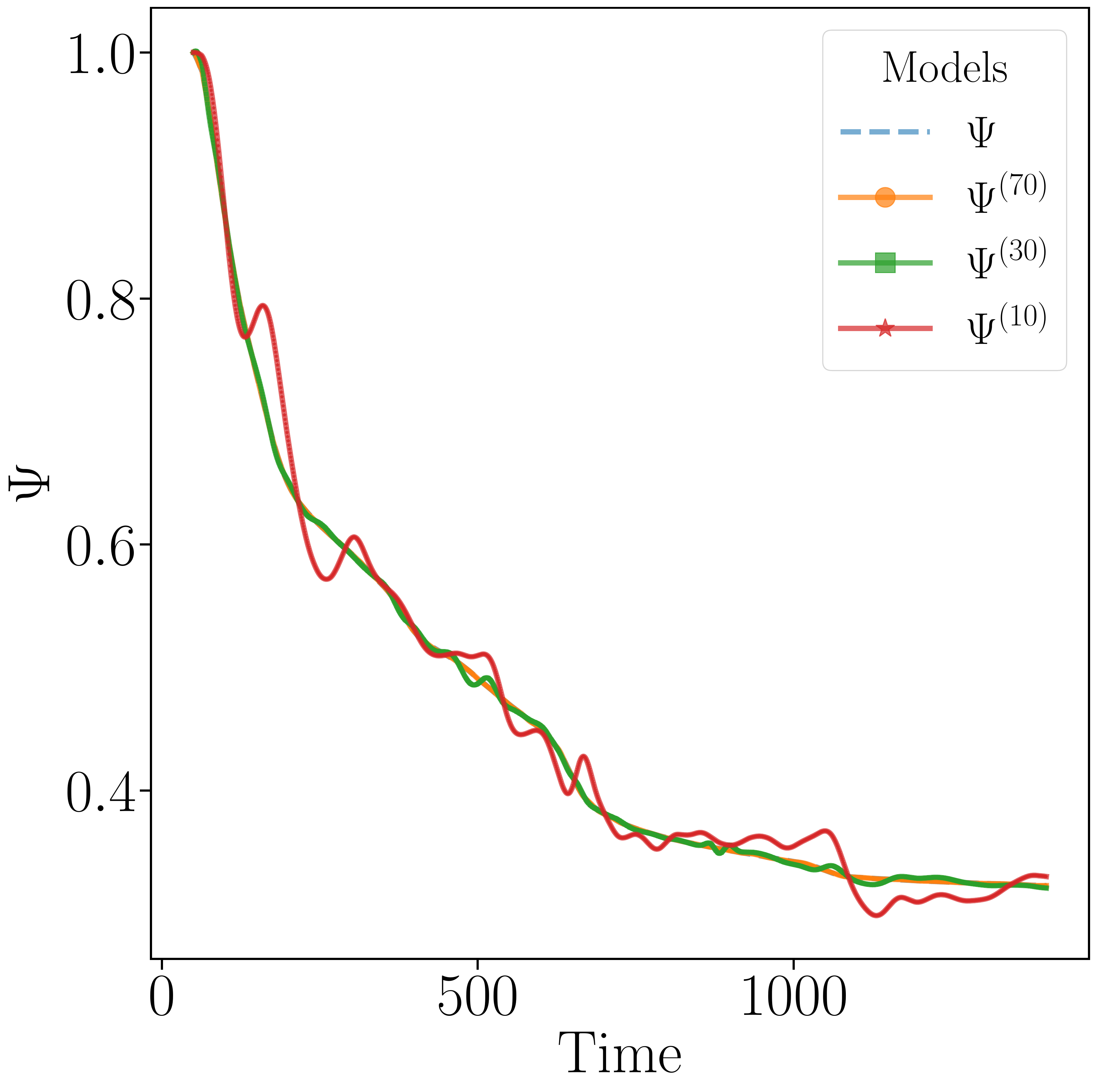}
	\subcaption{OLS regression.}
	\label{fig:psi_fit_OLS}
\end{subfigure}
\hfill
\begin{subfigure}[t]{0.49\textwidth}
	\centering
	\includegraphics[width=\textwidth]{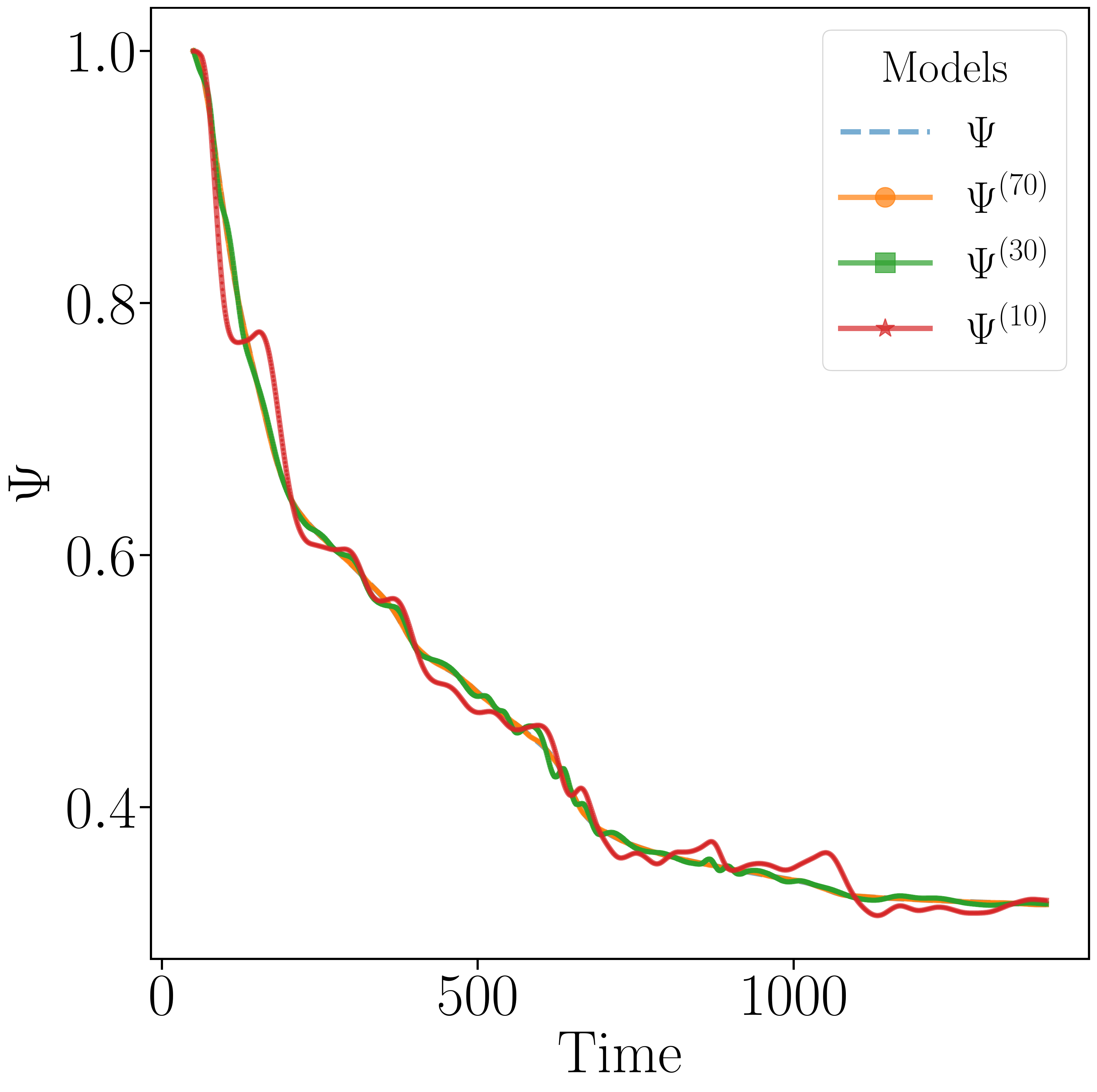}
	\subcaption{Ridge regression with $\lambda = 10^{-17}$.}
	\label{fig:psi_fit_Ridge17_local}
\end{subfigure}
\caption{Fitted curves for the $\ith[4]$ order free energy Taylor series functional representation, with $10$, $30$, and $70$ terms. DNS free energy data is shown with the dashed blue curve, and is essentially coincident with the higher order models.}
\label{fig:psi_fit_OLS_Ridge17_local}
\end{figure}

\noindent Regarding the fits for the free energy, Taylor series with $10$, $30$, and $70$ terms are shown in \cref{fig:psi_fit_OLS_Ridge17_local}. Free energies at earlier times also appear to be fit better by all models. When comparing the OLS and Ridge fits, the fits are almost identical for more complex models, however the OLS method has smaller oscillations for the most parsimonious models. This Taylor series approach is shown to be a logical and effective basis of terms to model this smooth free energy function. 


\end{appendices}


\end{document}

%% file: figures/model_diffusion_local.tex
\begin{subequations}
\renewcommand{\theequation}{\arabic{parentequation}~-~\the\numexpr6+6-1-\value{equation}\relax}

\begin{equation}
\begin{aligned}
{M_{\varphi}} =&~ 
{\gamma}^{{{{\varphi}_{\scriptscriptstyle{\nabla^2_{-}}}}^{}}}{{{\varphi}_{\scriptscriptstyle{\nabla^2_{-}}}}^{}} ~+~ 
{\gamma}^{{{\varphi_{\scriptscriptstyle{\nabla^2_{+}}}}^{}}}{{\varphi_{\scriptscriptstyle{\nabla^2_{+}}}}^{}} ~+~ 
{\gamma}^{{{{\varphi}_{4{-}}}^{}}}{{{\varphi}_{4{-}}}^{}} ~+~ 
{\gamma}^{{{{F_{+}}}^{}}}{{{F_{+}}}^{}}	 
\\
{\mathcal{E}_{\varphi}} =&~ {\gamma}^{{{\varphi_{{}}}^{}}}{{\varphi_{{}}}^{}} ~+~ 
{\gamma}^{{{\varphi_{3{+}}}^{}}}{{\varphi_{3{+}}}^{}} ~+~ 
{\gamma}^{{{{\varphi}_{\scriptscriptstyle{\nabla^2_{-}}}}^{}}}{{{\varphi}_{\scriptscriptstyle{\nabla^2_{-}}}}^{}} ~+~ 
{\gamma}^{{{{{F}^{\prime}_{-}}}^{}}}{{{{F}^{\prime}_{-}}}^{}} ~+~ 
{\gamma}^{{{{\varphi}_{5{-}}}^{}}}{{{\varphi}_{5{-}}}^{}} ~+~ 
{\gamma}^{{{\varphi_{5{+}}}^{}}}{{\varphi_{5{+}}}^{}}
\end{aligned}
\end{equation}

\begin{equation}
\begin{aligned}
{M_{\varphi}} =&~ {\gamma}^{{{{\varphi}_{\scriptscriptstyle{\nabla^2_{-}}}}^{}}}{{{\varphi}_{\scriptscriptstyle{\nabla^2_{-}}}}^{}} ~+~ 
{\gamma}^{{{\varphi_{\scriptscriptstyle{\nabla^2_{+}}}}^{}}}{{\varphi_{\scriptscriptstyle{\nabla^2_{+}}}}^{}} ~+~ 
{\gamma}^{{{{\varphi}_{4{-}}}^{}}}{{{\varphi}_{4{-}}}^{}} ~+~ 
{\gamma}^{{{{F_{+}}}^{}}}{{{F_{+}}}^{}}
\\
{\mathcal{E}_{\varphi}} =&~ {\gamma}^{{{\varphi_{{}}}^{}}}{{\varphi_{{}}}^{}} ~+~ 
{\gamma}^{{{\varphi_{3{+}}}^{}}}{{\varphi_{3{+}}}^{}} ~+~ 
{\gamma}^{{{{\varphi}_{\scriptscriptstyle{\nabla^2_{-}}}}^{}}}{{{\varphi}_{\scriptscriptstyle{\nabla^2_{-}}}}^{}} ~+~ 
{\gamma}^{{{{{F}^{\prime}_{-}}}^{}}}{{{{F}^{\prime}_{-}}}^{}} ~+~ 
{\gamma}^{{{{\varphi}_{5{-}}}^{}}}{{{\varphi}_{5{-}}}^{}}
\end{aligned}
\end{equation}

\begin{equation}
\begin{aligned}
{M_{\varphi}} =&~{\gamma}^{{{{\varphi}_{\scriptscriptstyle{\nabla^2_{-}}}}^{}}}{{{\varphi}_{\scriptscriptstyle{\nabla^2_{-}}}}^{}} ~+~ 
{\gamma}^{{{\varphi_{\scriptscriptstyle{\nabla^2_{+}}}}^{}}}{{\varphi_{\scriptscriptstyle{\nabla^2_{+}}}}^{}} ~+~ 
{\gamma}^{{{{\varphi}_{4{-}}}^{}}}{{{\varphi}_{4{-}}}^{}}
\\
{\mathcal{E}_{\varphi}} =&~ {\gamma}^{{{\varphi_{{}}}^{}}}{{\varphi_{{}}}^{}} ~+~ 
{\gamma}^{{{\varphi_{3{+}}}^{}}}{{\varphi_{3{+}}}^{}} ~+~ 
{\gamma}^{{{{\varphi}_{\scriptscriptstyle{\nabla^2_{-}}}}^{}}}{{{\varphi}_{\scriptscriptstyle{\nabla^2_{-}}}}^{}} ~+~ 
{\gamma}^{{{{{F}^{\prime}_{-}}}^{}}}{{{{F}^{\prime}_{-}}}^{}} ~+~ 
{\gamma}^{{{{\varphi}_{5{-}}}^{}}}{{{\varphi}_{5{-}}}^{}}
\end{aligned}
\end{equation}

\begin{equation}
\begin{aligned}
{M_{\varphi}} =&~{\gamma}^{{{{\varphi}_{\scriptscriptstyle{\nabla^2_{-}}}}^{}}}{{{\varphi}_{\scriptscriptstyle{\nabla^2_{-}}}}^{}} ~+~ 
{\gamma}^{{{\varphi_{\scriptscriptstyle{\nabla^2_{+}}}}^{}}}{{\varphi_{\scriptscriptstyle{\nabla^2_{+}}}}^{}} ~+~ 
\\
{\mathcal{E}_{\varphi}} =&~ {\gamma}^{{{\varphi_{{}}}^{}}}{{\varphi_{{}}}^{}} ~+~ 
{\gamma}^{{{\varphi_{3{+}}}^{}}}{{\varphi_{3{+}}}^{}} ~+~ 
{\gamma}^{{{{\varphi}_{\scriptscriptstyle{\nabla^2_{-}}}}^{}}}{{{\varphi}_{\scriptscriptstyle{\nabla^2_{-}}}}^{}} ~+~ 
{\gamma}^{{{{{F}^{\prime}_{-}}}^{}}}{{{{F}^{\prime}_{-}}}^{}} ~+~ 
{\gamma}^{{{{\varphi}_{5{-}}}^{}}}{{{\varphi}_{5{-}}}^{}}
\end{aligned}
\end{equation}

\begin{equation}
\begin{aligned}
{M_{\varphi}} =&~{\gamma}^{{{{\varphi}_{\scriptscriptstyle{\nabla^2_{-}}}}^{}}}{{{\varphi}_{\scriptscriptstyle{\nabla^2_{-}}}}^{}} 
\\
{\mathcal{E}_{\varphi}} =&~ {\gamma}^{{{\varphi_{{}}}^{}}}{{\varphi_{{}}}^{}} ~+~ 
{\gamma}^{{{\varphi_{3{+}}}^{}}}{{\varphi_{3{+}}}^{}} ~+~ 
{\gamma}^{{{{\varphi}_{\scriptscriptstyle{\nabla^2_{-}}}}^{}}}{{{\varphi}_{\scriptscriptstyle{\nabla^2_{-}}}}^{}} ~+~ 
{\gamma}^{{{{{F}^{\prime}_{-}}}^{}}}{{{{F}^{\prime}_{-}}}^{}} ~+~ 
{\gamma}^{{{{\varphi}_{5{-}}}^{}}}{{{\varphi}_{5{-}}}^{}}
\end{aligned}
\end{equation}

\begin{equation}
\begin{aligned}
{M_{\varphi}} =&~{\gamma}^{{{{\varphi}_{\scriptscriptstyle{\nabla^2_{-}}}}^{}}}{{{\varphi}_{\scriptscriptstyle{\nabla^2_{-}}}}^{}} 
\\
{\mathcal{E}_{\varphi}} =&~ {\gamma}^{{{\varphi_{{}}}^{}}}{{\varphi_{{}}}^{}} ~+~ 
{\gamma}^{{{\varphi_{3{+}}}^{}}}{{\varphi_{3{+}}}^{}} ~+~ 
{\gamma}^{{{{\varphi}_{\scriptscriptstyle{\nabla^2_{-}}}}^{}}}{{{\varphi}_{\scriptscriptstyle{\nabla^2_{-}}}}^{}} ~+~ 
{\gamma}^{{{{{F}^{\prime}_{-}}}^{}}}{{{{F}^{\prime}_{-}}}^{}}
\end{aligned}
\end{equation}

\begin{equation}
\begin{aligned}
{M_{\varphi}} =&~{\gamma}^{{{{\varphi}_{\scriptscriptstyle{\nabla^2_{-}}}}^{}}}{{{\varphi}_{\scriptscriptstyle{\nabla^2_{-}}}}^{}} 
\\
{\mathcal{E}_{\varphi}} =&~ {\gamma}^{{{\varphi_{{}}}^{}}}{{\varphi_{{}}}^{}} ~+~ 
{\gamma}^{{{\varphi_{3{+}}}^{}}}{{\varphi_{3{+}}}^{}} ~+~ 
{\gamma}^{{{{\varphi}_{\scriptscriptstyle{\nabla^2_{-}}}}^{}}}{{{\varphi}_{\scriptscriptstyle{\nabla^2_{-}}}}^{}}
\end{aligned}
\end{equation}

\begin{equation}
\begin{aligned}
{M_{\varphi}} =&~ 0
\\
{\mathcal{E}_{\varphi}} =&~ {\gamma}^{{{\varphi_{{}}}^{}}}{{\varphi_{{}}}^{}} ~+~ 
{\gamma}^{{{\varphi_{3{+}}}^{}}}{{\varphi_{3{+}}}^{}} ~+~ 
{\gamma}^{{{{\varphi}_{\scriptscriptstyle{\nabla^2_{-}}}}^{}}}{{{\varphi}_{\scriptscriptstyle{\nabla^2_{-}}}}^{}}
\end{aligned}
\end{equation}

\begin{equation}
\begin{aligned}
{M_{\varphi}} =&~ 0
\\
{\mathcal{E}_{\varphi}} =&~ {\gamma}^{{{\varphi_{{}}}^{}}}{{\varphi_{{}}}^{}} ~+~ 
{\gamma}^{{{\varphi_{3{+}}}^{}}}{{\varphi_{3{+}}}^{}}
\end{aligned}
\end{equation}

\begin{equation}
\begin{aligned}
{M_{\varphi}} =&~ 0
\\
{\mathcal{E}_{\varphi}} =&~ {\gamma}^{{{\varphi_{{}}}^{}}}{{\varphi_{{}}}^{}}
\end{aligned}
\end{equation}

\end{subequations}

%% file: figures/model_vol_OLS_Ridge17_local.tex
\begin{align}
\frac{{\delta}^{} {\bar{\varphi}}_{\textrm{OLS}}}{{\delta} {{{t}}_{}}^{}} =&~ {\gamma}^{{{{\bar{E}}_{12}}^{}{{\bar{E}}_{22}}^{}{{\bar{l}}_{1}}^{}}}{{{\bar{E}}_{12}}^{}{{\bar{E}}_{22}}^{}{{\bar{l}}_{1}}^{}} ~+~ {\gamma}^{{{{\bar{E}}_{12}}^{}{{\bar{\varphi}}_{}}^{}{{\bar{l}}_{1}}^{}}}{{{\bar{E}}_{12}}^{}{{\bar{\varphi}}_{}}^{}{{\bar{l}}_{1}}^{}} \label{eq:vol_model_OLS_local} \\
+&~ {\gamma}^{{{{\bar{E}}_{11}}^{}{{\bar{\varphi}}_{}}^{}}}{{{\bar{E}}_{11}}^{}{{\bar{\varphi}}_{}}^{}} ~+~ {\gamma}^{{{{\bar{E}}_{11}}^{}}}{{{\bar{E}}_{11}}^{}}  \nonumber \\ 
+&~   {\gamma}^{{{{\bar{E}}_{22}}^{}{{\bar{\varphi}}_{}}^{}{{\bar{N}}_{1}}^{}}}{{{\bar{E}}_{22}}^{}{{\bar{\varphi}}_{}}^{}{{\bar{N}}_{1}}^{}} ~+~ {\gamma}^{{{{\bar{E}}_{11}}^{}{{\bar{E}}_{22}}^{}{{\bar{l}}_{1}}^{}}}{{{\bar{E}}_{11}}^{}{{\bar{E}}_{22}}^{}{{\bar{l}}_{1}}^{}}  \nonumber \\ 
+&~    {\gamma}^{{{{\bar{E}}_{22}}^{2}{{\bar{l}}_{1}}^{}}}{{{\bar{E}}_{22}}^{2}{{\bar{l}}_{1}}^{}} ~+~ {\gamma}^{{{{\bar{E}}_{11}}^{}{{\bar{\varphi}}_{}}^{}{{\bar{l}}_{1}}^{}}}{{{\bar{E}}_{11}}^{}{{\bar{\varphi}}_{}}^{}{{\bar{l}}_{1}}^{}}   \nonumber \\ 
+&~  {\gamma}^{{{{\bar{E}}_{11}}^{}{{\bar{\varphi}}_{}}^{}{{\bar{l}}_{2}}^{}}}{{{\bar{E}}_{11}}^{}{{\bar{\varphi}}_{}}^{}{{\bar{l}}_{2}}^{}} ~+~ {\gamma}^{{{{\bar{E}}_{22}}^{}{{\bar{\varphi}}_{}}^{}{{\bar{l}}_{2}}^{}}}{{{\bar{E}}_{22}}^{}{{\bar{\varphi}}_{}}^{}{{\bar{l}}_{2}}^{}}  \nonumber \\ \nonumber \\
\frac{{\delta}^{} {\bar{\varphi}}_{\textrm{Ridge}}}{{\delta} {{{t}}_{}}^{}} =&~ {\gamma}^{{{{\bar{E}}_{12}}^{}{{\bar{E}}_{22}}^{}{{\bar{l}}_{1}}^{}}}{{{\bar{E}}_{12}}^{}{{\bar{E}}_{22}}^{}{{\bar{l}}_{1}}^{}} ~+~ {\gamma}^{{{{\bar{E}}_{12}}^{}{{\bar{\varphi}}_{}}^{}{{\bar{l}}_{1}}^{}}}{{{\bar{E}}_{12}}^{}{{\bar{\varphi}}_{}}^{}{{\bar{l}}_{1}}^{}}  \label{eq:vol_model_Ridge17_local} \\
+&~   {\gamma}^{{{{\bar{E}}_{11}}^{}{{\bar{\varphi}}_{}}^{}}}{{{\bar{E}}_{11}}^{}{{\bar{\varphi}}_{}}^{}} ~+~ {\gamma}^{{{{\bar{E}}_{11}}^{}}}{{{\bar{E}}_{11}}^{}}  \nonumber \\ 
+&~   {\gamma}^{{{{\bar{E}}_{22}}^{}{{\bar{\varphi}}_{}}^{}{{\bar{N}}_{1}}^{}}}{{{\bar{E}}_{22}}^{}{{\bar{\varphi}}_{}}^{}{{\bar{N}}_{1}}^{}} ~+~ {\gamma}^{{{{\bar{E}}_{11}}^{}{{\bar{E}}_{22}}^{}{{\bar{l}}_{1}}^{}}}{{{\bar{E}}_{11}}^{}{{\bar{E}}_{22}}^{}{{\bar{l}}_{1}}^{}}  \nonumber \\ 
+&~   {\gamma}^{{{{\bar{E}}_{22}}^{2}{{\bar{l}}_{1}}^{}}}{{{\bar{E}}_{22}}^{2}{{\bar{l}}_{1}}^{}} ~+~ {\gamma}^{{{{\bar{E}}_{11}}^{}{{\bar{\varphi}}_{}}^{}{{\bar{l}}_{1}}^{}}}{{{\bar{E}}_{11}}^{}{{\bar{\varphi}}_{}}^{}{{\bar{l}}_{1}}^{}}  \nonumber \\ 
+&~   {\gamma}^{{{{\bar{E}}_{11}}^{}{{\bar{\varphi}}_{}}^{}{{\bar{l}}_{2}}^{}}}{{{\bar{E}}_{11}}^{}{{\bar{\varphi}}_{}}^{}{{\bar{l}}_{2}}^{}} ~+~ {\gamma}^{{{{\bar{E}}_{22}}^{}{{\bar{\varphi}}_{}}^{}{{\bar{l}}_{2}}^{}}}{{{\bar{E}}_{22}}^{}{{\bar{\varphi}}_{}}^{}{{\bar{l}}_{2}}^{}} \nonumber
\end{align}

%% file: figures/model_psi_OLS_Ridge17_local.tex
\begin{align}
{\Psi}_{\textrm{OLS}} =&~ {\Psi}_{}({{{\bar{\varphi}}}}_{0},{{\mathbf{\bar{E}}}}_{0}) + {{\gamma}^{{{{\bar{\varphi}}_{}}^{2}}}\frac{1}{{2!}}{\frac{{\delta}^{2}{\Psi}_{}({{{\bar{\varphi}}}}_{0},{{\mathbf{\bar{E}}}}_{0})}{{\delta {\bar{\varphi}}_{}}^{2}}}{\Delta {\bar{\varphi}}_{}}^{2}}  \label{eq:psi_model_OLS_local}\\ 
+&~  {{\gamma}^{{{{\bar{\varphi}}_{}}^{}}{{{\bar{E}}_{12}}^{}}}{\frac{{\delta}^{2}{\Psi}_{}({{{\bar{\varphi}}}}_{0},{{\mathbf{\bar{E}}}}_{0})}{{\delta {\bar{\varphi}}_{}}^{} {\delta {\bar{E}}_{12}}^{}}}{\Delta {\bar{\varphi}}_{}}^{} {\Delta {\bar{E}}_{12}}^{}} + {{\gamma}^{{{{\bar{\varphi}}_{}}^{2}}{{{\bar{E}}_{12}}^{}}}\frac{1}{{2!}}{\frac{{\delta}^{3}{\Psi}_{}({{{\bar{\varphi}}}}_{0},{{\mathbf{\bar{E}}}}_{0})}{{\delta {\bar{\varphi}}_{}}^{2} {\delta {\bar{E}}_{12}}^{}}}{\Delta {\bar{\varphi}}_{}}^{2} {\Delta {\bar{E}}_{12}}^{}} \nonumber\\
+&~   {{\gamma}^{{{{\bar{E}}_{11}}^{}}{{{\bar{E}}_{22}}^{}}}{\frac{{\delta}^{2}{\Psi}_{}({{{\bar{\varphi}}}}_{0},{{\mathbf{\bar{E}}}}_{0})}{{\delta {\bar{E}}_{11}}^{} {\delta {\bar{E}}_{22}}^{}}}{\Delta {\bar{E}}_{11}}^{} {\Delta {\bar{E}}_{22}}^{}} ~+~ {{\gamma}^{{{{\bar{\varphi}}_{}}^{}}{{{\bar{E}}_{11}}^{}}{{{\bar{E}}_{22}}^{}}}{\frac{{\delta}^{3}{\Psi}_{}({{{\bar{\varphi}}}}_{0},{{\mathbf{\bar{E}}}}_{0})}{{\delta {\bar{\varphi}}_{}}^{} {\delta {\bar{E}}_{11}}^{} {\delta {\bar{E}}_{22}}^{}}}{\Delta {\bar{\varphi}}_{}}^{} {\Delta {\bar{E}}_{11}}^{} {\Delta {\bar{E}}_{22}}^{}} \nonumber\\
+&~   {{\gamma}^{{{{\bar{\varphi}}_{}}^{2}}{{{\bar{E}}_{11}}^{}}{{{\bar{E}}_{22}}^{}}}\frac{1}{{2!}}{\frac{{\delta}^{4}{\Psi}_{}({{{\bar{\varphi}}}}_{0},{{\mathbf{\bar{E}}}}_{0})}{{\delta {\bar{\varphi}}_{}}^{2} {\delta {\bar{E}}_{11}}^{} {\delta {\bar{E}}_{22}}^{}}}{\Delta {\bar{\varphi}}_{}}^{2} {\Delta {\bar{E}}_{11}}^{} {\Delta {\bar{E}}_{22}}^{}} + {{\gamma}^{{{{\bar{\varphi}}_{}}^{}}{{{\bar{E}}_{22}}^{}}}{\frac{{\delta}^{2}{\Psi}_{}({{{\bar{\varphi}}}}_{0},{{\mathbf{\bar{E}}}}_{0})}{{\delta {\bar{\varphi}}_{}}^{} {\delta {\bar{E}}_{22}}^{}}}{\Delta {\bar{\varphi}}_{}}^{} {\Delta {\bar{E}}_{22}}^{}} \nonumber\\
+&~   {{\gamma}^{{{{\bar{\varphi}}_{}}^{}}{{{\bar{E}}_{11}}^{}}}{\frac{{\delta}^{2}{\Psi}_{}({{{\bar{\varphi}}}}_{0},{{\mathbf{\bar{E}}}}_{0})}{{\delta {\bar{\varphi}}_{}}^{} {\delta {\bar{E}}_{11}}^{}}}{\Delta {\bar{\varphi}}_{}}^{} {\Delta {\bar{E}}_{11}}^{}} + {{\gamma}^{{{{\bar{\varphi}}_{}}^{2}}{{{\bar{E}}_{22}}^{}}}\frac{1}{{2!}}{\frac{{\delta}^{3}{\Psi}_{}({{{\bar{\varphi}}}}_{0},{{\mathbf{\bar{E}}}}_{0})}{{\delta {\bar{\varphi}}_{}}^{2} {\delta {\bar{E}}_{22}}^{}}}{\Delta {\bar{\varphi}}_{}}^{2} {\Delta {\bar{E}}_{22}}^{}} \nonumber \\ \nonumber \\
{\Psi}_{\textrm{Ridge}} =&~ {\Psi}_{}({{{\bar{\varphi}}}}_{0},{{\mathbf{\bar{E}}}}_{0}) + {{\gamma}^{{{{\bar{\varphi}}_{}}^{2}}}\frac{1}{{2!}}{\frac{{\delta}^{2}{\Psi}_{}({{{\bar{\varphi}}}}_{0},{{\mathbf{\bar{E}}}}_{0})}{{\delta {\bar{\varphi}}_{}}^{2}}}{\Delta {\bar{\varphi}}_{}}^{2}}  \label{eq:psi_model_Ridge17_local}\\
+&~ {{\gamma}^{{{{\bar{\varphi}}_{}}^{3}}}\frac{1}{{3!}}{\frac{{\delta}^{3}{\Psi}_{}({{{\bar{\varphi}}}}_{0},{{\mathbf{\bar{E}}}}_{0})}{{\delta {\bar{\varphi}}_{}}^{3}}}{\Delta {\bar{\varphi}}_{}}^{3}} + {{\gamma}^{{{{\bar{\varphi}}_{}}^{}}{{{\bar{E}}_{11}}^{}}{{{\bar{E}}_{22}}^{}}}{\frac{{\delta}^{3}{\Psi}_{}({{{\bar{\varphi}}}}_{0},{{\mathbf{\bar{E}}}}_{0})}{{\delta {\bar{\varphi}}_{}}^{} {\delta {\bar{E}}_{11}}^{} {\delta {\bar{E}}_{22}}^{}}}{\Delta {\bar{\varphi}}_{}}^{} {\Delta {\bar{E}}_{11}}^{} {\Delta {\bar{E}}_{22}}^{}} \nonumber\\
+&~   {{\gamma}^{{{{\bar{E}}_{11}}^{}}{{{\bar{E}}_{22}}^{}}}{\frac{{\delta}^{2}{\Psi}_{}({{{\bar{\varphi}}}}_{0},{{\mathbf{\bar{E}}}}_{0})}{{\delta {\bar{E}}_{11}}^{} {\delta {\bar{E}}_{22}}^{}}}{\Delta {\bar{E}}_{11}}^{} {\Delta {\bar{E}}_{22}}^{}} + {{\gamma}^{{{{\bar{\varphi}}_{}}^{}}{{{\bar{E}}_{11}}^{3}}}\frac{1}{{3!}}{\frac{{\delta}^{4}{\Psi}_{}({{{\bar{\varphi}}}}_{0},{{\mathbf{\bar{E}}}}_{0})}{{\delta {\bar{\varphi}}_{}}^{} {\delta {\bar{E}}_{11}}^{3}}}{\Delta {\bar{\varphi}}_{}}^{} {\Delta {\bar{E}}_{11}}^{3}}  \nonumber \\
+&~ {{\gamma}^{{{{\bar{\varphi}}_{}}^{}}{{{\bar{E}}_{11}}^{}}{{{\bar{E}}_{12}}^{}}}{\frac{{\delta}^{3}{\Psi}_{}({{{\bar{\varphi}}}}_{0},{{\mathbf{\bar{E}}}}_{0})}{{\delta {\bar{\varphi}}_{}}^{} {\delta {\bar{E}}_{11}}^{} {\delta {\bar{E}}_{12}}^{}}}{\Delta {\bar{\varphi}}_{}}^{} {\Delta {\bar{E}}_{11}}^{} {\Delta {\bar{E}}_{12}}^{}} + {{\gamma}^{{{{\bar{\varphi}}_{}}^{3}}{{{\bar{E}}_{11}}^{}}}\frac{1}{{3!}}{\frac{{\delta}^{4}{\Psi}_{}({{{\bar{\varphi}}}}_{0},{{\mathbf{\bar{E}}}}_{0})}{{\delta {\bar{\varphi}}_{}}^{3} {\delta {\bar{E}}_{11}}^{}}}{\Delta {\bar{\varphi}}_{}}^{3} {\Delta {\bar{E}}_{11}}^{}} \nonumber\\
+&~   {{\gamma}^{{{{\bar{\varphi}}_{}}^{2}}{{{\bar{E}}_{11}}^{}}}\frac{1}{{2!}}{\frac{{\delta}^{3}{\Psi}_{}({{{\bar{\varphi}}}}_{0},{{\mathbf{\bar{E}}}}_{0})}{{\delta {\bar{\varphi}}_{}}^{2} {\delta {\bar{E}}_{11}}^{}}}{\Delta {\bar{\varphi}}_{}}^{2} {\Delta {\bar{E}}_{11}}^{}} + {{\gamma}^{{{{\bar{\varphi}}_{}}^{}}{{{\bar{E}}_{11}}^{}}}{\frac{{\delta}^{2}{\Psi}_{}({{{\bar{\varphi}}}}_{0},{{\mathbf{\bar{E}}}}_{0})}{{\delta {\bar{\varphi}}_{}}^{} {\delta {\bar{E}}_{11}}^{}}}{\Delta {\bar{\varphi}}_{}}^{} {\Delta {\bar{E}}_{11}}^{}} \nonumber
\end{align}